\documentclass[CJK,11pt]{amsart}

\usepackage{geometry}
\usepackage{amsmath,amsfonts,amssymb,amsthm}
\usepackage{graphicx}
\usepackage{caption}
\usepackage{subcaption}
\usepackage{float}
\usepackage{dsfont}
\usepackage{epsfig}
\usepackage{esint}
\usepackage[colorlinks,linkcolor=blue]{hyperref}
\usepackage{verbatim}
\usepackage{algorithmicx}
\usepackage{algorithm}

\newcommand{\ignore}[1]{}

\newcommand{\del}{\delta\hspace{-0.2ex}}

\newtheorem{proposition}{Proposition}
\newtheorem{theorem}{Theorem}
\newtheorem{remark}{Remark}
\newtheorem{lemma}{Lemma}
\newtheorem{corollary}{Corollary}
\newtheorem{assumption}{Assumption}

\makeatletter
\newcommand\footnoteref[1]{\protected@xdef\@thefnmark{\ref{#1}}\@footnotemark}
\makeatother

\newcommand{\C}{\mathcal{C}}

\newcommand{\QQ}{\mathrm{Q}}
\newcommand{\R}{\mathbb{R}}

\newcommand{\cc}{\mathrm{C}}

\newcommand{\dd}{\mathrm{d}}

\newcommand{\bb}{\mathrm{B}}
\newcommand{\LL}{\mathrm{L}}
\newcommand{\rhs}{r.h.s.\,}
\newcommand{\lhs}{l.h.s.\,}
\newcommand{\calF}{\mathcal{F}}
\newcommand{\fintstar}{\fint_{\bb_{r_\star(x)}(x)}}
\newcommand{\rss}{r_{\star\star}}
\newcommand{\ahom}{a_{\textrm{hom}}}

\usepackage{color}
\definecolor{darkred}{rgb}{0.9,0.1,0.1}
\definecolor{darkblue}{rgb}{0,0,0.7}
\definecolor{darkgreen}{rgb}{0,0.5,0}

\begin{document}
\title[Artificial boundary conditions for random elliptic systems]{{\bf Artificial boundary conditions for random elliptic systems with correlated coefficient field}}
\author{Nicolas Clozeau \quad Lihan Wang}
\thanks{Institute of Science and Technology Austria, Klosterneuburg, Austria (nicolas.clozeau@ist.ac.at) \and Carnegie Mellon University, Pittsburgh, PA, USA (lihanw@andrew.cmu.edu)}
\thanks{NC has received funding from the European Research Council (ERC) under the Eu\-ropean Union’s Horizon 2020 research and innovation programme (grant agreement No 948819) \includegraphics[height=\fontcharht\font`\B]{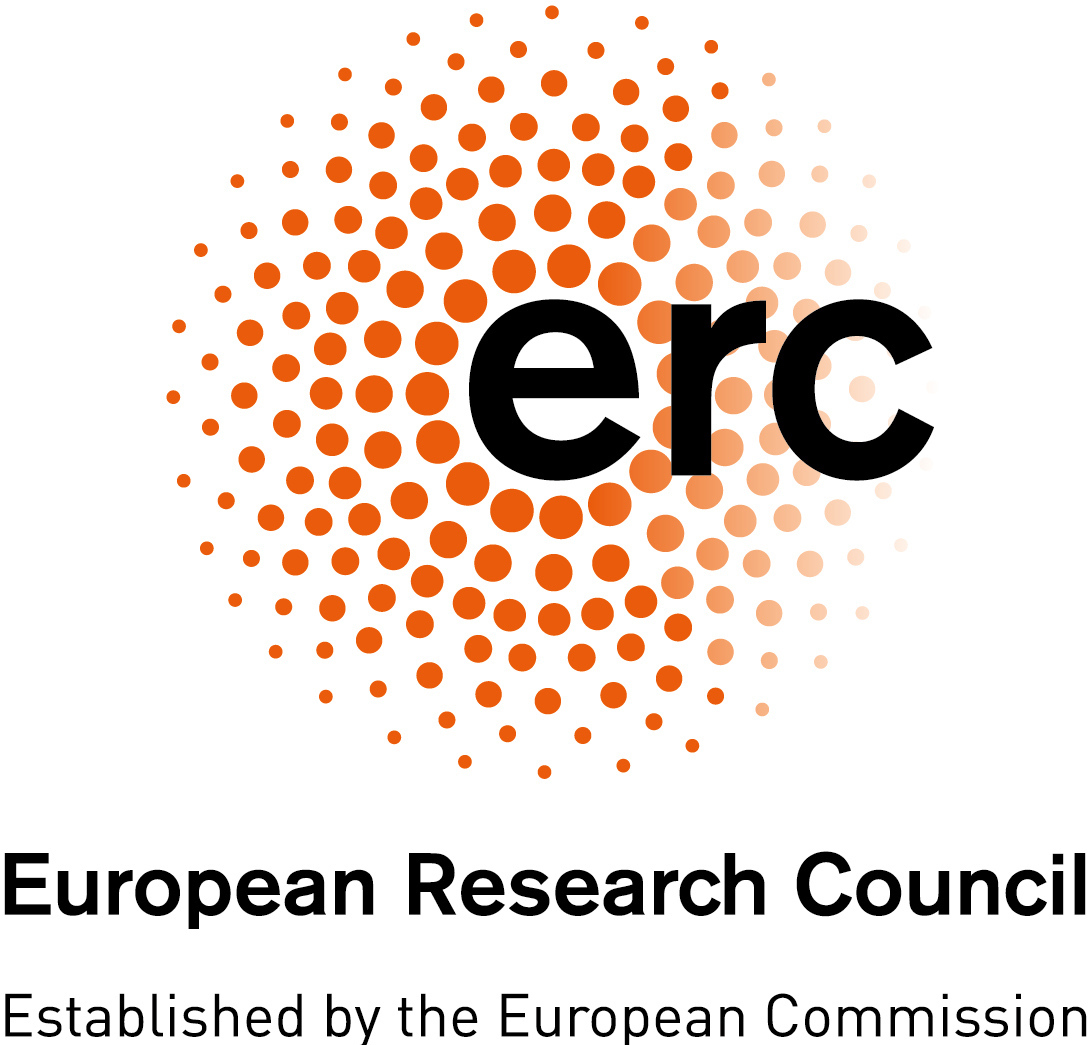}\,\includegraphics[height=\fontcharht\font`\B]{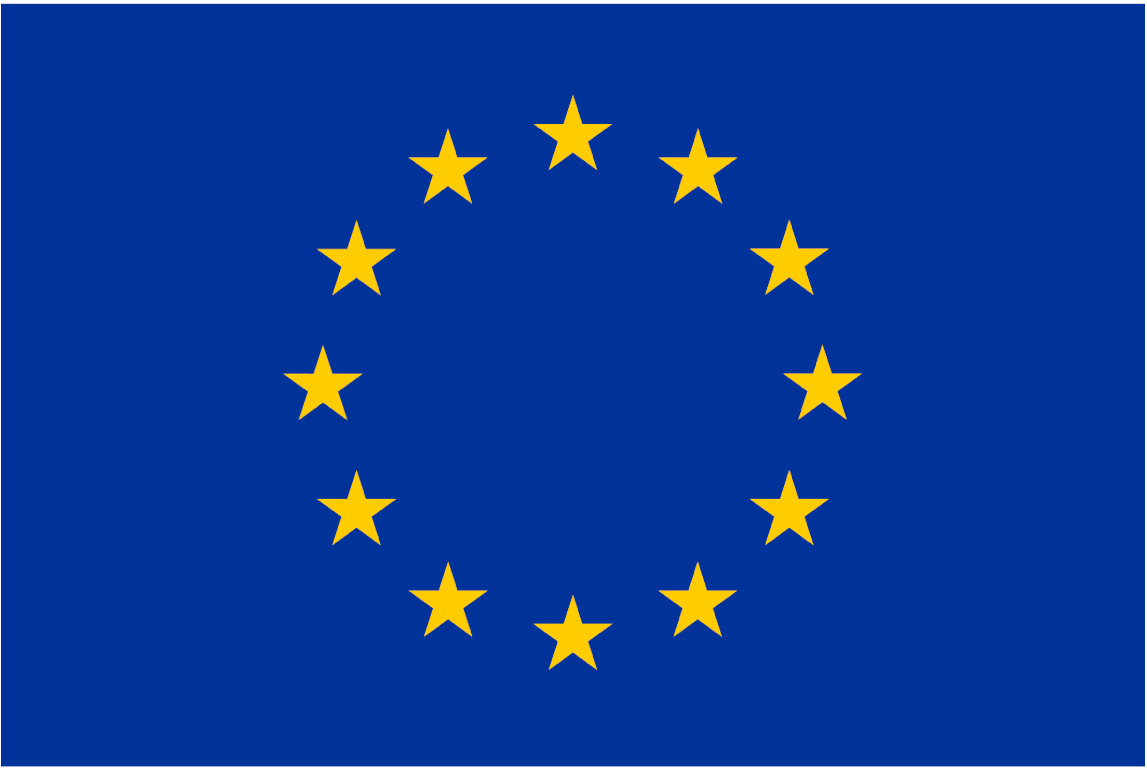}}
\begin{abstract}
We are interested in numerical algorithms for computing the electrical field generated by a charge distribution localized on scale $\ell$ in an infinite heterogeneous correlated random medium, in a situation where the medium is only known in a box of diameter $L\gg\ell$ around the support of the charge. We show that the algorithm in \cite{lu2021optimal}, suggesting optimal Dirichlet boundary conditions motivated by the multipole expansion \cite{bella2020effective}, still performs well in correlated media. With overwhelming probability, we obtain a convergence rate in terms of $\ell$, $L$ and the size of the correlations for which optimality is supported with numerical simulations. These estimates are provided for ensembles which satisfy a multi-scale logarithmic Sobolev inequality, where our main tool is an extension of the semi-group estimates in \cite{clozeau2021optimal}. As part of our strategy, we construct sub-linear second-order correctors in this correlated setting which is of independent interest.
\end{abstract}
\maketitle

\tableofcontents
%\section*{To Do List}
%%
%\begin{itemize}
%	%
%	\item Say somewhere the small regularity results we need for the quantities such that we refer it to.
%	\item Consistency of notations: min \& max or $\wedge \& \vee$? I'm fine with the latter (it's nicer in the exponent), but we need to define it at the very beginning if we go that way.
%	\item add $d,\beta>2$ in assumptions of lemmas
%	%
%\end{itemize}
\section{Introduction and main results}
\subsection{Random conducting media}
We consider a conducting medium described by a symmetric $\lambda$-uniformly elliptic coefficient field $a$ in $d\geq 2$ dimensional space, that is
\begin{equation}\label{UniformElliptiIntro}
\lambda\vert\zeta\vert^2\leq \zeta\cdot a(x)\zeta\leq \vert\zeta\vert^2\quad\text{for any $\zeta,x\in\mathbb{R}^d$.}
\end{equation}
We are interested in computing  $\nabla u$, where $u$ is the decaying solution of the elliptic divergence-form equation
\begin{equation}\label{EquationIntro}
-\nabla\cdot a\nabla u=\nabla\cdot h\quad\text{in $\mathbb{R}^d$,}
\end{equation}
interpreted as the electric field generated by the neutral charge distribution $\nabla\cdot h$. Furthermore, we assume that the charge is localized at a scale $\ell\geq 1$ : there exists $\widehat{h}$ compactly supported in the unit ball $\bb$ such that
\begin{equation}\label{StructureAssumptionh}
h=\widehat{h}(\tfrac{\cdot}{\ell}).
\end{equation}
From the two assumptions \eqref{UniformElliptiIntro} and \eqref{StructureAssumptionh}, we define $u$ as the energy solution of \eqref{EquationIntro} in $\dot{\mathrm{H}}^1:=\mathrm{H}^1\slash \mathbb{R}$ where
$$\mathrm{H}^1:=\bigg\{u: \mathbb{R}^d\rightarrow\mathbb{R}\,\text{ such that }\, \int_{\mathbb{R}^d} \vert \nabla u\vert^2<\infty\bigg\}.$$
In addition to the medium $a$ satisfying \eqref{UniformElliptiIntro}, we assume that the conductivity information is only known statistically. Concretely, it amounts to assume that $a$ is sampled from a stationary and ergodic probability measure $\mathbb{P}$. By the latter we mean a probability measure on the space of tensor fields $a$ satisfying \eqref{UniformElliptiIntro}. Stationarity means that for any shift vector $z\in\mathbb{R}^d$, the shifted random field $x\mapsto a(x+z)$ has the same (joint) distribution as $a$. Ergodicity is a qualitative assumption that encodes the decorrelation of the values of $a$ over large distances. Here, we are deliberately vague on the $\sigma$-algebra and on the notion of ergodicity because we will consider a very explicit class in this paper so that quantitative results can be obtained, see Section \ref{sec:assump}.
\subsection{Optimal artificial boundary conditions}\label{subsec:ABC}
For numerical purposes, we analyze the approximation of \eqref{StructureAssumptionh} in a finite domain $\QQ_L:=[-L,L)^d$ for $L\gg1$, that is 
\begin{equation}\label{IntroBC}
\left\{
    \begin{array}{ll}
        -\nabla\cdot a\nabla u^{(L)}=\nabla\cdot h & \text{in $\QQ_L$,} \\
        u^{(L)}=u^{(L)}_{\mathrm{bc}} & \text{on $\partial\QQ_L$,}
    \end{array}
\right.
\end{equation}
for some $u^{(L)}_{\text{bc}}$ to be determined. We follow the strategy of Lu, Otto and the second author in \cite{lu2018optimal,lu2021optimal} where they establish an expression of $u^{(L)}_{\text{bc}}$ allowing for stochastic cancellations. In the case where the probability distribution $\mathbb{P}$ has a \textbf{finite range of dependence}\footnote{meaning that there exists $r\geq 1$ such that $a\vert_U$ and $a\vert_V$ are independent for any open sets $U,V$ satisfying $\text{dist}(U,V)\geq r$}, they show that, compared to the na\"ive approach of setting $u^{(L)}_{\text{bc}}=0$, a suitable choice of the boundary condition improves the rate of convergence by the order of the fluctuation scaling $O(L^{-\frac{d}{2}})$ (see \cite[Theorem $1.2$]{lu2021optimal}). Moreover, optimality in this setting is proven in \cite[Theorem $2$]{lu2018optimal}.

%In the recent contribution \cite{lu2018optimal,lu2021optimal}, the authors proposed an algorithm to calculate the value of $\nabla u$ using only the information of the coefficient field $a$ on a large box $Q_{2L}:=[-2L,2L)^d$ to design suitable artificial Dirichlet boundary conditions and solve \eqref{EquationIntro} in $Q_L$. They show that in the case where $\mathbb{P}$ has a \textbf{finite range of dependence}, stochastic cancellations occur and improve the rate of convergence, compared to a na\"ive approach, by the order of the Central Limit Theorem (CLT) scaling $O((\frac{1}{L})^\frac{d}{2})$. Moreover, the optimality of the method in the finite range of dependence setting is proven in \cite[Theorem $2$]{lu2018optimal}.

%
\medskip

Our goal in the present contribution is to go beyond the independent case and generalize the analysis in \cite{lu2021optimal} to the setting where the probability distribution $\mathbb{P}$ may possess \textbf{long-range correlations}. Such correlated coefficient fields can be, for instance, generated from diverse point processes or Gaussian fields which are known to not satisfy the finite range of dependence assumption (see for instance \cite{duerinckx2017weighted}) and are used in practical models of interest in the applied sciences (see \cite{torquato2002random}). For the purpose of the paper, we illustrate the effect of the correlations for a specific class of Gaussian type coefficient fields, see Assumption \ref{Gaussian} and Assumption \ref{DecayCor}.  Our class includes for instance log-normal random coefficients, that is of the type
$$a(x)=\frac{b_1+e^{-\tilde{\kappa}(g(x)-m)}}{b_2+e^{-\kappa(g(x)-m)}}\text{Id}\quad\text{for any $x\in\mathbb{R}^d$},$$
where $b_1,b_2>0$, $\tilde{\kappa},\kappa,m\in\mathbb{R}$ and $g$ is a stationary mean-zero (scalar) Gaussian field with decaying correlations %\lw{Since $\beta$ already appears at this point, I would prefer to discuss the critical index $2$ earlier in section 1.2 instead of 1.3. Hence it would be good to emphasize this below formula.}
\begin{equation}\label{eq:cx}
	|c(x)|\lesssim (1+\vert x\vert)^{-\beta}, \ \text{ where }c(x):=\mathbb{E}[g(x)g(0)]
\end{equation} for some $\beta>0$ (possibly very small). We show that depending on the size $\beta$ of the correlations and the dimension we can construct the first-order boundary condition as in \cite{lu2018optimal} or second-order boundary condition as in \cite{lu2021optimal}, which allows to capture stochastic cancellations of the order of the fluctuation scaling\footnote{We use $a\wedge b:=\min\{a,b\}$ and $a\vee b := \max\{a,b\}$ throughout this work.} $O(L^{-\frac{\beta\wedge d}{2}})$ in both cases. The critical values $\beta\wedge d=2$ directly corresponds to the borderline case where the gradient of second-order corrector  cannot be constructed in a $\LL^2$-stationary sense (this will be discussed in more details in Section \ref{DiscussTwoScale}). We state in the following a post-processed version of our main Theorem \ref{thm:EffBdryAlg} showing what we achieve in this paper.
\begin{corollary}\label{CoroIntro}
We define $u^{(L)}$ the solution of \eqref{IntroBC} with the boundary condition $u^{(L)}_{\mathrm{bc}}$ as defined in \eqref{eqn:algapproxbdry} when $d,\beta >2$ or \eqref{eqn:2dproxbdry} when $\beta \wedge d\le 2$. There exists an exponent $\gamma>0$ depending on $d$ and $\beta$ such that for any $\varepsilon>0$ there exists a constant $C>0$ for which for any $0<R\leq L$ and $\frac{L}{C}\geq \ell\geq C$, it holds
\begin{equation}\label{IntroMainResult}
\mathbb{P}\bigg(\|\nabla u-\nabla u^{(L)}\|_{\LL^{2}(\bb_R)}\leq C\Big(\frac{\ell}{L}\Big)^d L^{-\frac{\beta\wedge d}{2}(1-\varepsilon)}\bigg)\geq \big(1-C(\exp(-\ell^{\gamma})+\exp(-R^\gamma)\big)\big)\big(1-\exp(-L^{\frac{1}{C}})\big).
\end{equation}
\end{corollary}
We refer to Theorem \ref{thm:EffBdryAlg} for a more precise statement and Appendix \ref{appendix:Alg} for the precise description of the algorithm. Despite we prove our main result only for this specific class of coefficient fields, our methods adapt to various other type of randomness and we refer to Section \ref{ExtensionOtherType} for further discussions. Additionally to the theoretical result in Corollary \ref{CoroIntro}, we perform numerical simulations in Section \ref{NumericalSimulations} which supports the scaling $O((\frac{\ell}{L})^d L^{-\frac{\beta\wedge d}{2}})$. 

\medskip

Our proof strategy is in the vein of the series of works \cite{lu2021optimal,clozeau2021optimal,gloria2015corrector,armstrong2019quantitative}. As in \cite{lu2021optimal}, we combine the semigroup strategy first developed in \cite{gloria2015corrector} and \cite[Chapter $9$]{armstrong2019quantitative} with the multipole expansion theory developed in \cite{bella2020effective}. Our first main contribution are the semigroup estimates in Section \ref{TimeDecaySemigroup}, where our proof techniques differ highly from \cite{lu2021optimal} and rather relies on sensitivity calculus, functional inequalities and large-scale regularity in the spirit of \cite{clozeau2021optimal}. As a by-product, our second main contribution is the construction of sub-linear second-order correctors in our correlated setting in Section \ref{ConstructionSecondCor} which is of independent interest. Combined together with the strategy in \cite{lu2021optimal} based on the multipole expansion in \cite{bella2020effective}, we deduce \eqref{IntroMainResult}.

\medskip

Before we present the technical assumptions and the rigorous formulation of our main results, we briefly review the construction of $u^{(L)}_{\text{bc}}$ in \cite{lu2021optimal} and introduce the quantities therin, where the theory of stochastic homogenization and the theory of multipole expansion in random media kick in. 
\subsubsection{Large-scale approximation via two-scale expansion}\label{DiscussTwoScale}
The (qualitative) homogenization theory, which started with the pioneer works of Kozlov \cite{kozlov1979averaging} and Papanicolaou and Varadhan \cite{papanicolaou1979boundary}, states that the large-scale behavior of $(-\nabla\cdot a\nabla)^{-1}$ is captured by $(-\nabla\cdot a_{\mathrm{hom}}\nabla)^{-1}$ where the constant and deterministic homogenized coefficient $a_{\text{hom}}$ is given by\footnote{where $\{e_i\}_{1\leq i\leq d}$ denotes the canonical basis of $\mathbb{R}^d$}, 
\begin{equation}\label{intro:ah}
	\ahom e_i = \lim_{R\uparrow\infty} \fint_{\mathrm{Q}_R}  q_i^{(1)}\quad\text{with }q_i^{(1)}:= a  (e_i+\nabla \phi_i^{(1)})\quad \text{for any $i\in\{1,\cdots,d\}$}.
\end{equation}
In \eqref{intro:ah}, $\phi^{(1)}=\{\phi^{(1)}_i\}_{1\leq i\leq d}$ is the so-called first-order corrector, unique (up to an additive constant) random field with stationary gradient\footnote{that is, for any $x\in\mathbb{R}^d$, $\nabla\phi^{(1)}_i(a,\cdot+x)=\nabla\phi^{(1)}_i(a(\cdot+x),\cdot)$ almost-surely} that solves in the distributional sense 
\begin{equation}\label{intro:phi1}
-\nabla \cdot a (e_i+\nabla \phi_i^{(1)})=0\quad\text{in $\mathbb{R}^d$}\quad\text{with }\mathbb{E}\big[\nabla\phi^{(1)}_i\big]=0\text{ and }\mathbb{E}\big[\vert\nabla\phi^{(1)}_i\vert^2\big]<\infty.
\end{equation}
With the help of $a_\text{hom}$, we can compute cheaply the homogenized solution $u_\text{hom}$ of 
\begin{equation}\label{EquationHom}
-\nabla\cdot \ahom\nabla u_{\text{hom}}= \nabla \cdot h,
\end{equation}
which will play an important role in our algorithm. A quantitative description of this large-scale behavior is provided by the so-called two-scale expansion and we refer to the works of Gloria, Neukamm and Otto \cite{gloria2014regularity,gloria2015quantification,gloria2017quantitative} as well as Armstrong, Kuusi, Mourrat and Smart \cite{armstrong2016quantitative,armstrong2019quantitative,armstrong2016lipschitz} for a complete overview on quantitative stochastic homogenization. 
	In particular, the first-order flux corrector $\sigma^{(1)}$ introduced in \cite{gloria2014regularity} in the stochastic setting (see \cite[Proposition 7.2]{jikov2012homogenization} for periodic homogenization) plays an important role and jointly $(\phi^{(1)},\sigma^{(1)})$ governs the homogenization error. For any indices $i\in\{1,\cdots,d\}$, $\sigma^{(1)}_i=\{\sigma^{(1)}_{ijk}\}_{1\leq j,k\leq d}$ is the unique (up to an additive constant) skew-symmetric tensor with stationary gradient which solves in the distributional sense\footnote{we use the notation $(\nabla \cdot \sigma^{(1)}_i)_j=\sum_{k}\partial_k \sigma^{(1)}_{ijk}$}
\begin{equation}\label{intro:sig1}
	\left\{ \begin{aligned}
	- \Delta \sigma_{ijk}^{(1)} =& \nabla \cdot (q_{ik}^{(1)}e_j - q_{ij}^{(1)}e_k), \\  
	 \nabla \cdot \sigma^{(1)}_i =&\,q_{i}
	\end{aligned}\right.\quad\text{with $\mathbb{E}\big[\nabla\sigma^{(1)}_i\big]=0$ and $\mathbb{E}\big[\vert\nabla\sigma^{(1)}_i\vert^2\big]<\infty$.}
\end{equation}
For earlier works on quantitative stochastic homogenization, we refer the readers to \cite{yurinskii1986averaging, naddaf1997homogenization, naddaf1998estimates}. When $d>2$, higher-order accuracy may be obtained by using second-order corrector $\phi^{(2)}=\{\phi^{(2)}_{ij}\}_{1\leq i,j\leq d}$ which is, for any indices $i,j\in\{1,\cdots, d\}$, the distributional solution with stationary gradient of 
\begin{equation}\label{intro:phi2}
-\nabla\cdot a \nabla \phi_{ij}^{(2)} = \nabla \cdot (a\phi_i^{(1)}-\sigma_i^{(1)})e_j\quad\text{in $\mathbb{R}^d$}\quad\text{with $\mathbb{E}[\nabla\phi^{(2)}_{ij}]=0$ and $\mathbb{E}\big[\vert\nabla\phi^{(2)}_{ij}\vert^2\big]<\infty$.}
\end{equation} 
More precisely, we can upgrade the two-scale expansion to second-order, that is\footnote{we use Einstein's notation for repeated indices}
\begin{equation}\label{TwoScaleIntro}
u^{\text{2sc}}:=(1+\phi^{(1)}_i\partial_i+\phi^{(2)}_{ij}\partial_{ij})u_{\text{hom}},
\end{equation}
and gain higher-order accuracy in the large-scale approximation (see for instance \cite{bella2017stochastic}). The existence of $\phi^{(2)}$ is not guaranteed for general stationary and ergodic distribution of coefficient fields and requires the existence of stationary first-order correctors (as opposed to only their gradients), which itself asks for stronger ergodicity assumptions. It is established when the probability distribution satisfies a logarithmic Sobolev inequality in \cite{clozeau2021optimal, gloria2019quantitative} or a finite range of dependence assumption in \cite{lu2021optimal}. With this consideration, to be able to define $\phi^{(2)}$ rigorously in our correlated setting \eqref{eq:cx}, we not only require $d>2$ but we also need to strengthen the quantitative ergodicity conditions to $\beta>2$, identified as the regime where stationary first-order correctors can be constructed (see \cite[Corollary $3$]{clozeau2021optimal}). For the construction of second- and higher-order correctors, we refer the readers to the works \cite{fischer2016higher, gu2017high, duerinckx2019higher}.

\medskip

Analogously to the first-order case, we need second-order flux correctors $\sigma^{(2)}$, which were first introduced in \cite{bella2017stochastic} in the stochastic case, since the spatial growth of $(\phi^{(2)},\sigma^{(2)})$ governs the homogenization error at the second-order. When $(\phi^{(1)},\sigma^{(1)})$ are stationary, we can define for each pair of indices $i,j \in \{1,\cdots,d\}$ a skey-symmetric tensor $\sigma_{ij}^{(2)} = \{\sigma_{ijkn}^{(2)}\}_{1\le k,n \le d}$, which is unique (up to an additive constant), has a stationary gradient, and solves in the distributional sense\footnote{we use the notation $(\nabla\times q_{ij}^{(2)})_{kn}:=\partial_k q_{ijn}^{(2)}-\partial_n q_{ijk}^{(2)}$}
\begin{equation*}
\left\{
\begin{array}{ll}
	\nabla\cdot \sigma^{(2)}_{ij}=q_{ij}^{(2)}:=a\nabla\phi^{(2)}_{ij}+(\phi^{(1)}_i a-\sigma_i^{(1)})e_j, &  \\
	-\Delta \sigma^{(2)}_{ij}=\nabla\times q_{ij}^{(2)},& 
\end{array}
\right. \quad\text{with $\mathbb{E}\big[\nabla\sigma^{(2)}_{ij}\big]=0$ and $\mathbb{E}\big[\vert\nabla\sigma^{(2)}_{ij}\vert^2\big]<\infty$.}
\end{equation*}

\medskip

As one may notice, the quantities $(\phi^{(1)},\sigma^{(1)},\phi^{(2)})$ still satisfy the whole-space equations \eqref{intro:phi1},\eqref{intro:sig1} and \eqref{intro:phi2} which are not exactly computable from a single realization $a|_{\QQ_{2L}}$ restricted in a finite domain. Furthermore, direct Dirichlet approximations with homogeneous boundary conditions of these quantities may introduce parasite boundary effects and impact the convergence rate by a surface factor $L^{d-1}$. To bypass this issue, we introduce the so-called massive correctors $(\phi^{(1)}_{M},\sigma^{(1)}_{M},\phi^{(2)}_{M})$ belonging to the space
$$\mathrm{H}^1_{\text{uloc}}(\mathbb{R}^d):=\bigg\{\psi\in \mathrm{H}^1_{\text{loc}}(\mathbb{R}^d)\Big\vert \sup_{x\in\mathbb{R}^d}\int_{\bb(x)}\vert\psi\vert^2+\vert\nabla\psi\vert^2<\infty\bigg\},$$
where, for $M\geq 1$, $\phi^{(1)}_M=\{\phi^{(1)}_{i,M}\}_{1\leq i\leq d}$ satisfies
\begin{equation}\label{eq:1stOrderMassCor} 
\tfrac{1}{M}\phi^{(1)}_{i,M}-\nabla\cdot a(e_i+\nabla\phi^{(1)}_{i,M})=0\quad\text{in $\mathbb{R}^d$},
\end{equation}
$\sigma^{(1)}_M=\{\sigma^{(1)}_{ijk,M}\}_{1\leq i,j,k\leq d}$ satisfies
\begin{equation}\label{eq:1stOrderMassFlux}
	\frac{1}{M}\sigma_{ijk,M}^{(1)} - \Delta \sigma_{ijk,M}^{(1)} = \nabla \cdot (q_{ik,M}^{(1)}e_j - q_{ij,M}^{(1)}e_k)\quad\text{with }q_{i,M}^{(1)} = a(e_i+\nabla \phi_{i,M}^{(1)}),
\end{equation}
and $\phi^{(2)}_M=\{\phi^{(2)}_{ij,M}\}_{1\leq i,j\leq d}$ satisfies
\begin{equation}\label{eq:2ndOrderMassCor} 
\tfrac{1}{M}\phi^{(2)}_{ij,M}-\nabla\cdot a\nabla\phi^{(2)}_{ij,M}= \nabla \cdot (a\phi^{(1)}_{i,M}-\sigma^{(1)}_{i,M})e_j\quad\text{in $\mathbb{R}^d$}.
\end{equation} 
%
%Sincedirect Dirichlet approximations of these quantities may be insufficient to provide good approximation rates, in this work we are considering massive approximations of all these quantities. Namely, for any $M\ge 1$, we define $(\phi_{i,M}^{(1)})_i$ to be the first-order massive corrector
%
%
%and $\sigma^{(1)}_{i,M} = (\sigma_{ijk,M})_{j,k}$ to denote the first-order massive flux corrector
%
%with first-order flux $q_{i,M}^{(1)} = a(e_i+\nabla \phi_{e_i,M}^{(1)})$. Similarly, we may construct $(\phi_{ij,M}^{(2)})_{ij}$ to be the solution of 
%
On the one hand, using probabilistic tools, we can prove that $(\phi^{(1)}_M,\sigma^{(1)}_M,\phi^{(2)}_M)$ provide good approximations of $(\phi^{(1)},\sigma^{(1)},\phi^{(2)})$ as $M\uparrow\infty$, which we optimally quantify in Lemma \ref{lem:1stCorApprox} and Corollary \ref{Cor:Approx2ndCor}. On the other hand, using deterministic tools, it is proven in \cite[Proposition 2.8]{lu2021optimal} that $(\phi^{(1)}_M,\sigma^{(1)}_M,\phi^{(2)}_M)$ can be further approximated by $(\phi^{(1)}_{M,L},\sigma^{(1)}_{M,L},\phi^{(2)}_{M,L})$, defined in \eqref{eqn:phiML}, \eqref{eqn:algsigma} and \eqref{eqn:2ndcorapprox} in Appendix \ref{appendix:Alg}, which solve the same equations but with zero Dirichlet boundary conditions on boxes of size $O(L)$, and the errors are subalgebraic in $L$ for the suitable choice $M \sim L^{2-}$.
%The quantities $(\phi^{(1)}_M,\sigma^{(1)}_M,\phi^{(2)}_M)$ are shift-covariant themselves, so we may characterize their approximation error with their non-massive counterparts using probabilistic tools. Moreover, as is proved in \cite{lu2021optimal} using only deterministic tools, they can be further approximated by imposing zero Dirichlet boundary conditions with a subalgebraic error, which are realized in \cite[Algorithm 1]{lu2021optimal}. A description of the algorithm including the roles of dipole and quadpole correction terms is given in more detail in Appendix \ref{appendix:Alg}.
%
\subsubsection{Multipole expansion in random media}
As observed in \cite{bella2015quantitative, bella2020effective}, na\"ively solving 
\begin{equation*}
\left\{\begin{array}{ll}
		-\nabla\cdot a\nabla u_1=\nabla\cdot h & \text{in $\QQ_L$,} \\
		u_1=(1+\phi_i^{(1)}\partial_i + \phi_{ij}^{(2)}\partial_{ij}) u_{\text{hom}} & \text{on $\partial\QQ_L$}
	\end{array}
	\right.\end{equation*}
will not provide a better approximation of $\nabla u(x)$ for large $|x|$, since it fails to capture the correct \emph{multipole behavior} of the solution, which is the far-field effect generated by the intrinsic moments of a localized right-hand side $\nabla\cdot h$. Instead of rigorous explanations, for which we refer the readers to \cite{bella2020effective,lu2021optimal}, here we provide a heuristic argument: the solution $u_{\text{hom}}$ of \eqref{EquationHom} has the following Green's function representation 
\begin{equation}\label{Greenrepuh}
	u_{\text{hom}}(x) = \int_{\R^d} \dd y \, G_{\text{hom}}(x-y) \nabla \cdot h(y),
\end{equation}
where $G_{\text{hom}}$ is the fundamental solution of $-\nabla\cdot a_{\text{hom}}\nabla$. By assumption \eqref{StructureAssumptionh}, in the above integral \eqref{Greenrepuh}, $y$ is only supported in $\bb_{\ell}$, hence for $x\in \partial\QQ_L$ we have $|x|\gg |y|$, which means that we can perform a Taylor expansion on \eqref{Greenrepuh} and obtain at leading order
\begin{equation*}
	u_{\text{hom}}(x) \approx  \int_{\R^d} \dd y \, \Big(G_{\text{hom}}(x) - y_i \partial_iG_{\text{hom}}(x) + \frac{y_iy_j}{2}\partial_{ij}G_{\text{hom}}(x)   \Big) \nabla \cdot h(y).
\end{equation*}
Since the $\nabla\cdot h$ is neutral, the contribution of the first term vanishes, and we can therefore rewrite
\begin{equation}\label{uhmultpole}
u_{\text{hom}}(x) \approx - \underbrace{\partial_iG_{\text{hom}}(x) \int_{\R^d} \dd y \,  y_i  \nabla \cdot h(y)}_{\text{Dipole term}}+ \underbrace{\partial_{ij}G_{\text{hom}}(x) \int_{\R^d} \dd y \, \frac{y_iy_j}{2}   \nabla \cdot h(y)}_{\text{Quadrupole term}}.
\end{equation}
In other words, $u_\text{hom}$ admits a (second-order) multi-pole expansion where the coefficients are given by moments of the neutral charge against harmonic polynomials. For the random heterogeneous operator $-\nabla\cdot a\nabla$, the harmonic coordinates are given by $y\mapsto y_i+\phi^{(1)}_i$, and thus by similar arguments we obtain the leading order expression (up to dipole behaviour)
\begin{equation}\label{umultpole}
	u(x) \approx -\partial_iG(x) \int_{\R^d} \dd y \,  (y_i +\phi_i^{(1)} )\nabla \cdot h(y),
\end{equation} 
where $G$ denotes the fundamental solution of $-\nabla\cdot a \nabla$. Hence, comparing \eqref{uhmultpole} with \eqref{umultpole}, neglecting the small error between $G$ and $G_{\text{hom}}$, we need to correct $u_\text{hom}$ by adding the \emph{Dipole correction term}
\begin{equation}\label{MultiPoleAnsatz}
u_{\text{hom}} \mapsto u_{\text{hom}} - \partial_iG_{\text{hom}}(x) \int_{\R^d} \dd y \,  \phi_i^{(1)} \nabla \cdot h(y).
\end{equation}
In the regime $d=2$ or $0<\beta\le 2$, as already mentioned before, second-order ``harmonic polynomials'' for the random operator $-\nabla\cdot a\nabla$ cannot be explicitly constructed and thus one cannot improve upon the first-order expansion. In this case, we therefore take for $\xi_i^{(1)} := \int_{\R^d} \dd y \,  \phi_i^{(1)} \nabla \cdot h(y),$
\begin{equation}\label{EquationBCIntro2D}
	u_{\text{bc}}:=\underbrace{(1+\phi^{(1)}_{i}\partial_i)}_{\text{First-order two-scale expansion}}\Big(u_{\text{hom}}\quad-\quad\underbrace{\xi^{(1)}_{i}\partial_iG_{\text{hom}}(x)}_{\text{Dipole correction}}\Big).
\end{equation}
In the regime $d>2$ and $\beta>2$, second-order harmonic polynomials are well-defined and we can thus expand to second-order, resulting in a \emph{Quadrupole correction term} that is of the similar form $\xi_{ij}^{(2)}\partial_{ij} G_{\text{hom}}(x)$. We refer the reader to \cite[(39)]{lu2021optimal} for a derivation of the expression of $\xi_{ij}^{(2)}$.
%\nc{Detail about the differences and difficulties in the approach of the correlated case compared to the finite range case. Describe here the algorithm, etc... Here introduce first and second order correctors as well as massive. Point out that stationary correctors and sub-linear second order correctors are required so that too slowly decaying correlations is an issue.}  
%
Thus, combining the second-order two-scale expansion \eqref{TwoScaleIntro} with \eqref{MultiPoleAnsatz}, we build the boundary condition for \eqref{IntroBC} as
\begin{equation}\label{EquationBCIntro}
u_{\text{bc}}:=\underbrace{(1+\phi^{(1)}_{i}\partial_i+\phi^{(2)}_{ij}\partial_{ij})}_{\text{Second-order two-scale expansion}}\Big(u_{\text{hom}}\quad-\quad\underbrace{\xi^{(1)}_{i}\partial_iG_{\text{hom}}(x)}_{\text{Dipole correction}}\quad+\quad\underbrace{\xi^{(2)}_{ij}\partial_{ij}G_{\text{hom}}(x)}_{\text{Quadrupole correction}}\Big).
\end{equation}
The algorithm in Appendix \ref{appendix:Alg} is then established by replacing the true correctors $(\phi^{(1)},\sigma^{(1)},\phi^{(2)})$ with their proxies $(\phi^{(1)}_{M,L},\sigma^{(1)}_{M,L},\phi^{(2)}_{M,L})$, that are computable from a realization $a\vert_{\text{Q}_{2L}}$, with the scaling $M\sim L^{2-}$.
%
%Our main contribution is the construction of sub-linear second-order correctors in our correlated setting, see Assumption \ref{Gaussian} and Assumption \ref{DecayCor} for more precisions. Our approach follows the strategy developed in \cite{clozeau2021optimal}, where the original idea comes from \cite{gloria2015corrector,gloria2015quantification,armstrong2019quantitative}, based on the semigroup associated to the corrector equation \eqref{intro:phi2}. It is in the vein of the series of work \cite{gloria2014regularity,gloria2015quantification,gloria2017quantitative} which are based on a sensitivity calculus and functional inequalities. {\color{red}detail more}.
%
%We now state our assumptions on the coefficient field and the rigorous statements of our main results.
%
\subsection{Assumptions and formulation of the main results}\label{sec:assump}
We now introduce the class of probability measures $\mathbb{P}$ that we consider which generates $\lambda$-uniformly elliptic coefficient fields $a$. We assume that $a$ is given by a non-linear transformation of a Gaussian field $g$ generated from a standard white noise $\zeta$ in $\LL^2(\mathbb{R}^d)$. 
\begin{assumption}[Gaussian-type coefficient]\label{Gaussian}
We fix the space dimension $d\ge 2$ and consider a centred stationary scalar\footnote{That is only made for notational simplicity, we may consider values in any finite-dimensional linear space.} white noise $\zeta$ on $(\Omega,\mathcal{A},\mathbb{P}_0)$. We define $\mathbb{P}$ (with expectation denoted by $\mathbb{E}[\cdot]$) as the push-forward of $\mathbb{P}_0$ under the map
\begin{equation}\label{Ass:Coeff}
\Omega\ni \omega\mapsto a:=\big(x\mapsto A(g(\omega,x))\big)\quad\text{with } g:=m\star\zeta\text{ and }m\in\LL^2(\mathbb{R}^d),
\end{equation}
where $A$ is a smooth map which takes values into the set of $\lambda$-uniformly elliptic, symmetric and bounded matrices (for some parameter $\lambda>0$ fixed) satisfying
\begin{equation}\label{RegA}
\sup_{g\in\mathbb{R}}\vert A'(g)\vert+\vert A''(g)\vert<\infty.
\end{equation}
\end{assumption}
We further assume that $g$ is ergodic with a quantitative description of the decay of correlations allowing for long-range interactions.
\begin{assumption}[Decay of correlations]\label{DecayCor}
We assume that the model $m$ satisfies
\begin{equation}\label{Regm}
m\in\cc^{2}(\mathbb{R}^d)\cap \mathrm{H}^2(\mathbb{R}^d),
\end{equation}
and that there exists a parameter $\beta>0$ such that 
\begin{equation}\label{Ass:DecayCor}
\sup_{x\in\mathbb{R}^d}(1+\vert x\vert)^{\frac{1}{2}(d+\beta)}\big(\mathds{1}_{\beta\neq d}+\log^{\frac{1}{2}}(1+|x|)\mathds{1}_{\beta=d}\big)\vert m(x)\vert <\infty.
\end{equation}
\end{assumption}
We now comment on some direct consequences of Assumptions \ref{Gaussian} and  \ref{DecayCor}. On the one hand, \eqref{Ass:DecayCor} implies an algebraic decay on the covariance $c:=\mathbb{E}_0[g(\cdot)g(0)]=m\star m$, namely
\begin{equation}\label{DecayCovarianceFunction}
\sup_{x\in\R^d}(1+|x|)^{\beta}\vert c(x)\vert<\infty.
\end{equation}
The assertions \eqref{DecayCovarianceFunction} and \eqref{Ass:Coeff} ensure that $\mathbb{P}$ satisfies a \emph{multiscale logarithmic Sobolev inequality} (MSLSI), see \cite[Theorem 3.1]{duerinckx2017weighted}, that is: For any square-integrable random variable $\mathcal{F}$ of the coefficient field $a$, it holds 
\begin{equation}\label{msLSI}
	\mathbb{E}\Big[\calF^2 \log \tfrac{\calF^2}{\mathbb{E}[ \calF^2]}\Big] \lesssim \mathbb{E}\bigg[ \int_1^\infty \dd \ell\, \ell^{-d-\beta-1} \int_{\R^d}\dd x\, |\partial_{x,\ell}^{{\rm fct}} \calF|^2 \bigg],
\end{equation}
where 
\begin{equation}\label{DefDerivative}
\partial_{x,\ell}^{{\rm fct}} \calF(a):=\sup\bigg\{\limsup_{h\downarrow 0}\frac{\mathcal{F}(a+h\delta a)-\mathcal{F}(a)}{h}, \sup_{\bb_\ell(x)}\vert\delta a\vert\leq 1, \delta a=0\text{ outside }\bb_{\ell}(x)\bigg\}.
\end{equation}
It is known, see \cite[Proposition 1.10]{duerinckx2017weighted2}, that \eqref{msLSI} allows to control arbitrary algebraic moments, that is for any $p<\infty$ and any random variable $\mathcal{F}$ with $\mathbb{E}[\vert \mathcal{F}\vert^p]<\infty$,
\begin{equation}\label{msPIp}
\mathbb{E}\big[ |\calF-\mathbb{E}[\calF]|^p\big]^\frac{1}{p} \lesssim \sqrt{p}\,\mathbb{E}\bigg[ \Big( \int_1^\infty \dd \ell\, \ell^{-d-\beta-1} \int_{\R^d}\dd x\,|\partial_{x,\ell}^{fct} \calF|^2 \Big)^\frac{p}{2}\bigg]^\frac{1}{p}.
\end{equation}
The control of moments \eqref{msPIp} implies fine concentration properties for non-linear functions $\mathcal{F}=\mathcal{F}(a)$ and will be extensively used for stochastic estimates. 

\medskip

On the other hand, the assertion \eqref{Regm} ensures that almost surely, the realizations of $a$ belong to $\cc^{1,\alpha}_\text{loc}(\mathbb{R}^d)$ for any $\alpha\in (0,1)$, namely
\begin{equation}\label{SmoothnessCoef}
\sup_{x\in\mathbb{R}^d}\mathbb{E}\Big[\|a\|^p_{\cc^{1,\alpha}(\bb_1(x))}\Big]^{\frac{1}{p}}\lesssim_{\alpha}\sqrt{p}\quad \text{for any $\alpha\in (0,1)$ and $p<\infty$.}
\end{equation}
We shortly repeat the standard Kolmogorov argument for \eqref{SmoothnessCoef}. Defining the (vector-valued) Gaussian field $\nabla g=\nabla m\star \zeta$ with covariance
$$\tilde{c}(z)=\mathbb{E}\big[\nabla g(z)\otimes \nabla g(0)\big]=\int_{\mathbb{R}^d}\dd x\, \nabla m(z-x)\otimes \nabla m(x),$$
we have from \eqref{Regm} $\sup_{x,y}\vert x-y\vert^{-2}\mathbb{E}[\vert\nabla g(x)-\nabla g(y)\vert^2]=2\sup_z\vert z\vert^{-2}\text{tr}(\tilde{c}(0)-\tilde{c}(z))<\infty$. By Gaussianity, this extends to arbitrary moments : $\sup_{x,y}\vert x-y\vert^{-1}\mathbb{E}[\vert\nabla g(x)-\nabla g(y)\vert^p]^{\frac{1}{p}}\lesssim \sqrt{p}$. Estimating the H\"older semi-norm $[\nabla g]^p_{\alpha,\bb_1}$ by the Besov norm $\int_{\bb_1}\dd z\, \vert z\vert^{-d-p\alpha}\int_{\bb_1}\dd x\,\vert \nabla g(x+z)-\nabla g(x)\vert^p$, one derives $\sup_{x}\mathbb{E}\Big[[\nabla g]^p_{\alpha,\bb_1(x)}\Big]^{\frac{1}{p}}\lesssim \sqrt{p}$ for any $\alpha<1$ and $p<\infty$. By \eqref{RegA}, this transmits to the coefficient field $a$, that is \eqref{SmoothnessCoef}. The assertion \eqref{SmoothnessCoef} is ensured to simplify the proofs as it will be extensively used to apply Schauder's theory for linear elliptic and parabolic equations. We believe that this is not essential and that the result of Theorem \ref{thm:EffBdryAlg} still holds if \eqref{SmoothnessCoef} is dropped.

\smallskip
%

%Finally, if $d>2$ and we strengthen our assumption \eqref{DecayCovarianceFunction} to $\beta>2$, then the first-order correctors $(\phi^{(1)},\sigma^{(1)})$ themselves are stationary, see for instance \cite[Corollary $3$]{clozeau2021optimal}, which is required to build in Section \ref{ConstructionSecondCor} the sub-linear second-order correctors $(\phi^{(2)},\sigma^{(2)})$. 
%

%, the second assertion of \eqref{Ass:DecayCor} ensures local regularity on the coefficient field in form of 
%
%\begin{equation}\label{RegCoef}
%\sup_{x\in\mathbb{R}^d}\mathbb{E}\Big[\|a\|^p_{\cc^{0,\alpha}(\bb_1(x))}\Big]<\infty \quad \text{ for any $\alpha<1$ and $p<\infty$.}
%\end{equation}
%
%
%One important motivation of this paper arises from a numerical analysis perspective \cite{lu2021optimal}. Consider the field $\nabla u$ which is the decaying solution of
%the elliptic divergence-form equation 
%\begin{align}\label{eqn:intrbaseq}
	%\nabla\cdot(a\nabla u+g)=0\quad\mbox{in}\;\mathbb{R}^d.
%\end{align}
%where the \rhs describes a compactly supported and overall neutral charge distribution. Let us denote the characteristic scale of $g$ by
%$\ell$, assuming that it is of the form
%\begin{align}\label{eqn:conditionrhs}
	%g(x)=\hat g(\frac{x}{\ell})
%\end{align}
%for some sufficiently smooth $\hat g$ supported in the unit ball. We want to give the best possible approximation of the solution $\nabla u$ given a realization of $a|_{Q_{2L}}$. 

\medskip

%Due to the correlation of the coefficient field $a$, the optimal error that an algorithm based on artificial boundary conditions can achieve is unclear; indeed, in the case $\beta \le d$, due to the correlations, one may ``obtain more information'' of $a$ from its realization in $a|_{Q_{2L}}$ compared to the finite-range scenario considered in \cite{lu2021optimal}. \lw{We may further illustrate this in numerical examples. Depending on our numerical results, we may want to expand this further.} On the other hand, however, the behavior and approximation of $\nabla\phi^{(2)}$ may become worse, as is shown in Corollaries \ref{Cor:Growth2ndCor} and \ref{Cor:Approx2ndCor}. Nevertheless, the algorithm in \cite{lu2021optimal} still provides an educated prediction of $\nabla u$. 
%
Our main result shows an estimate on the error for the approximation of $\nabla u$ by $\nabla u^{(L)}$ solution of \eqref{IntroBC} with the boundary condition \eqref{EquationBCIntro}.
\begin{theorem}[Effective boundary conditions]\label{thm:EffBdryAlg}
We fix $R\leq L$, the dimension $d\in\{2,3\}$ and consider a coefficient field $a$ sampled from $\mathbb{P}$ as defined in Assumption \ref{Gaussian} and Assumption \ref{DecayCor}. Consider $u$ and $h$ related via \eqref{EquationIntro} where $h$ satisfies \eqref{StructureAssumptionh}. We denote by $u^{(L)}$ the output of \cite[Algorithm 1]{lu2021optimal}, recalled in Appendix \ref{appendix:Alg}. 

\medskip

For any $\varepsilon \in (0,\frac{1}{4})$, there exist a constant $C$ depending on $\lambda$, $\hat g$, $\beta$, $\varepsilon$ and a random radius $\rss$, such that conditioning on $\ell\ge \rss$ and $R\geq \rss$, with probability at least $1 - \exp\big(-L^{\frac{1}{C}}\big)$, we have
\begin{equation}\label{eq:AlgApprox}
		\|\nabla u-\nabla u^{(L)}\|_{\LL^2(\bb_R)}
		\leq C \Big(\frac{\ell}{L}\Big)^d\bigg(\frac{\rss}{L}\bigg)^{\frac{\beta\wedge d}{2}(1-\varepsilon)} \quad\mbox{ provided } \;
		\frac{L}{C}\ge  \ell \ge C.
	\end{equation}
Moreover, there exists $\gamma>0$ such that
\begin{align}\label{eqn:thm1rss}
\mathbb{E}\big[\exp\big(\rss^\gamma\big)\big]\le C.
\end{align}
\end{theorem}
\begin{remark}
We obtain Corollary \ref{CoroIntro} as a combination of the conditional probability \eqref{eq:AlgApprox} and the moment bound \eqref{eqn:thm1rss}. More precisely, given the two events 
$$\mathcal{E}_1:=\Big\{\|\nabla u-\nabla u^{(L)}\|_{\LL^2(\bb_R)}
		\leq C \big(\tfrac{\ell}{L}\big)^d\big(\tfrac{\rss}{L}\big)^{\frac{\beta\wedge d}{2}(1-\varepsilon)}\Big\}\quad\text{and}\quad \mathcal{E}_2:=\big\{\ell\geq r_{\star\star}\big\}\cap\big\{R\geq r_{\star\star}\big\},$$
we have 
$$\mathbb{P}(\mathcal{E}_1)\geq \mathbb{P}(\mathcal{E}_2)\mathbb{P}(\mathcal{E}_1\vert\mathcal{E}_2).$$
We deduce \eqref{IntroMainResult} using that \eqref{eq:AlgApprox} provides
$$\mathbb{P}(\mathcal{E}_1\vert\mathcal{E}_2)\geq 1-\exp(-L^{\frac{1}{C}}),$$
and \eqref{eqn:thm1rss} provides by Markov's inequality
$$\mathbb{P}(\mathcal{E}_2)\geq 1-C(\exp(-\ell^\gamma)+\exp(-R^{\gamma})).$$
\end{remark}
 Throughout the paper, we focus on the proof of the second-order case when $d>2$ and $\beta>2$, as the proofs for the first-order case $\beta \wedge d\leq 2$ are simpler and rely only on estimates of the first-order quantitative homogenization theory established in \cite{clozeau2021optimal}. 

\medskip

Our next result concerns optimality of Theorem \ref{thm:EffBdryAlg}. Our goal is to heuristically argue that in the regime of weak correlations, that is $\beta \gg d$, the scaling in \eqref{eq:AlgApprox} is essentially optimal in the sense that there exists a class of Gaussian coefficient field for which we can expect at most a gain of the order of $O(L^{-\frac{d}{2}})$ from random fluctuations. For simplicity, our example is constructed in the setting of the discrete lattice $\mathbb{Z}^d$ and for media with  small contrast, that is for some $\eta\in (0,1)$, our coefficient field $a$ is of the form
 \begin{equation}
 	\label{intro:forma}
a:=(1+A(\eta g))\mathrm{Id},
\end{equation} 
where $A(x) \approx x$ near $x=0$ and $g=\{\mathcal{G}_n\}_{n\in\mathbb{Z}^d}$ is a family of Gaussian random variables that will be specified below. Therefore, instead of working on $\nabla u$ we consider its small contrast approximation $\nabla \bar{u}$ defined as the unique decaying solution of
\begin{equation}\label{SmallContrastEquation}
	-\Delta \bar{u}=\nabla\cdot g\nabla v\quad\text{with}\quad -\Delta v=\nabla\cdot h.
\end{equation}
The equation \eqref{SmallContrastEquation} is motivated as follows : Combining \eqref{SmallContrastEquation} and \eqref{EquationIntro}, the difference $e:=u-(v+\eta\bar{u})$ solves
$$-\nabla\cdot a\nabla e=\eta\nabla\cdot (a-1)\nabla\bar{u}+\nabla\cdot (a-1-\eta g)\nabla v.$$
Thus, one may notice that from energy estimates the difference between $\nabla u$ and $\nabla (v+\eta \bar u)$ is $O(\eta^2)$. Moreover, since $v$ is deterministic, we also see that the scaling (in terms of $\ell, L$) of the variance of $\nabla u$, conditioned on $a|_{Q_{L}}$, should be reflected by that of $\nabla \bar{u}$, which is easier to study since its dependence on the randomness is essentially linear.
%We show optimality of Theorem \ref{thm:EffBdryAlg}, in the sense that we can expect at most a gain of the order of $L^{-\frac{\beta\wedge d}{2}}$ from random fluctuations, for a specific class of coefficient fields. That is done in the regime of small ellipticity contrast (that is $\lambda\ll 1$ in \eqref{UniformElliptiIntro}) for simplifying the presentation and computations. For a given scalar field $a(x)=A(g(x))\text{Id}$ satisfying Assumption \ref{Gaussian}, we introduce by Gaussianity the notation
%
%\begin{equation*}
%
%\mathbb{E}\big[A(g(x))A(g(0))\big]=\mathcal{A}(c(x))\text{Id},
%
%\end{equation*}
%
%with
%
%\begin{equation}\label{CorrelationCoefField}
%
%\begin{aligned}
%
%\mathcal{A}(s):=(2\pi)^{-\frac{d}{2}}\int_{\mathbb{R}^d\times\mathbb{R}^d} &A(z_1)A(z_2)\frac{1}{\sqrt{\text{det}(\Sigma_s)}}\exp\Big(-\frac{1}{2}(z_1,z_2)^{T}\Sigma^{-1}_s(z_1,z_2)\Big)\\
%&\text{and $\Sigma_s:=
   %\begin{pmatrix} 
        % c(0)&s \\ s&c(0) 
   %\end{pmatrix}$}.
%
%\end{aligned}
%\end{equation}
%
%
\begin{proposition}\label{OptimalScaling}
Let $g:=\{\mathcal{G}_n\}_{n\in\mathbb{Z}^d}$ be a Gaussian sequence of identically distributed and centred one-dimensional Gaussian random variables such that there exists $\beta>0$ for which 
\begin{equation}\label{correlationExample}
c(n-n'):=\mathbb{E}[\mathcal{G}_n\mathcal{G}_{n'}]=(1+\vert n-n'\vert)^{-\beta}\quad\text{for any $n,n'\in\mathbb{Z}^d$.}
\end{equation}
%
%Let the coefficient field $a$ be of the form \eqref{intro:forma} with $A : \mathbb{R}\rightarrow \mathbb{R}$ satisfying 
%
%\begin{equation}\label{AssumptionAheuristic}
%
%\sup_{g\in\mathbb{R}}\vert A(g)\vert<1,\quad \sup_{g\in\mathbb{R}}\vert A'(g)\vert+\vert A''(g)\vert<\infty,\quad A(0)=0\text{ and $A'(0)=1$}.
%
%\end{equation}
%
Let $\bar{u}$ be defined as in \eqref{SmallContrastEquation} %be the unique decaying solution of the discrete partial differential equation on $\mathbb{Z}^d$
%
%\begin{equation}\label{DiscreteEquation}
%
%-\nabla\cdot a\nabla u=\nabla\cdot h,
%
%\end{equation}
%
where $h$ satisfies \eqref{StructureAssumptionh} with the normalisation $\sum_{n\in\mathbb{Z}^d} h(n)=\ell^d\,\nu$ for some $\nu\in\mathbb{R}^d$.

\medskip

Then, there exists $\beta_d>d$ such that for any $\beta\geq \beta_d$ we have
\begin{equation}\label{Optimality}
\begin{aligned}
\mathbb{E}\Big[\big\vert\nabla \bar{u}(0)-\mathbb{E}\big[\nabla \bar{u}(0)\big\vert \mathcal{F}_L\big]\big\vert^2\Big]^{\frac{1}{2}}&\gtrsim_{\beta} \Big(\frac{\ell}{L}\Big)^{d}L^{-\frac{d}{2}},%-\eta^2\|h\|\\
% &\text{with $\|h\|:=\|h\|_{\ell^{2}(\mathbb{Z}^d)}+\Big(\sum_{n\in\mathbb{Z}^d}(1+\vert n\vert)^{2\alpha}\vert h(n)\vert^4\Big)^{\frac{1}{4}}$,}
\end{aligned}
\end{equation}
where $\mathbb{E}[\cdot\vert\mathcal{F}_L]$ denotes the conditional expectation with respect to the family $\mathcal{F}_L:=\{\mathcal{G}_n,\vert n\vert_\infty\leq L\}$.
%where $\mathbb{E}[\cdot\vert \mathcal{F}_L]$ denotes the conditional expectation with respect to $\{\mathcal{G}_n, \vert n\vert_\infty\leq L\}$.
%
\end{proposition}
Extensions of Proposition \ref{OptimalScaling} to the regime of high contrast or in the setting $\beta \lesssim d$ would require more involved arguments, and constitute an interesting open question to investigate. In particular, when $\beta<d$ the upper bound \eqref{eq:AlgApprox} worsen as $\beta$ decreases (mainly due to the growth of $\phi^{(2)}$, see Corollary \ref{Cor:Growth2ndCor}) while, on the other hand, stronger correlations of the random field may provide more information on $a|_{\QQ_{2L}\backslash \QQ_L}$ which may lower the conditional variance of the solution. 
\subsection{Extension to other type of coefficient fields}\label{ExtensionOtherType}
The result in Theorem \ref{thm:EffBdryAlg} is not limited to the Gaussian setting and holds as soon as the coefficient field $a$ satisfies a MSLSI of type \eqref{msLSI}. For instance, for various coefficient generated from a point process $\mathcal{P}$, namely $a=A(\mathcal{P})$ where $A$ satisfies \eqref{RegA}, it holds for some positive weight $\pi$
\begin{equation}\label{ExampleOscillation}
	\mathbb{E}\Big[\calF^2 \log \tfrac{\calF^2}{\mathbb{E}[ \calF^2]}\Big] \lesssim \mathbb{E}\bigg[ \int_1^\infty \dd \ell\, \ell^d\pi(\ell) \int_{\R^d}\dd x\, |\partial_{x,\ell}^{{\rm osc}} \calF|^2 \bigg], 
\end{equation}
where 
$$\partial_{x,\ell}^{{\rm osc}} \calF:=\sup_{A':A'\vert_{\mathbb{R}^d\backslash \bb_\ell(x)}=A\vert_{\mathbb{R}^d\backslash \bb_\ell(x)}}\calF(A')-\inf_{A':A'\vert_{\mathbb{R}^d\backslash \bb_\ell(x)}=A\vert_{\mathbb{R}^d\backslash \bb_\ell(x)}}\calF(A').$$
The MSLSI \eqref{ExampleOscillation} is satisfied for instance for Poisson tessellations and random parking measures with $\pi(\ell)=\exp(-\tfrac{1}{C}\ell^d)$ or for Poisson inclusions with random radii with $\pi(\ell)=\ell^d\,\mathbb{P}(\ell-1\leq \mathcal{R}\leq \ell)$ where $\mathcal{R}$ denotes the law of the radius. For more details, we refer to \cite{duerinckx2017weighted}. For some comments on the proof techniques in this case, we refer to \cite[Appendix D]{clozeau2021optimal}.
\subsection{Numerical Simulations}\label{NumericalSimulations}
%
%Similar to \cite{lu2021optimal}
To avoid any discretization errors, we consider the equation \eqref{EquationIntro} in the discrete setting $\mathbb{Z}^3$ and the coefficients are evaluated at the edges of $\mathbb{Z}^3$. We also consider only coefficient fields generated from a Gaussian random field that satisfy Assumption \ref{Gaussian} and \eqref{DecayCovarianceFunction}. We define the map $A:\mathbb{R} \to \mathbb{R}^{d\times d}$ into the space of symmetric, $4$-uniformly elliptic and bounded matrices by 
\begin{equation}\label{numerics:Amap}
	A(g) = \left(1 + \frac{3}{1+e^{-g}}\right) \mathrm{Id}. 
\end{equation}
The Gaussian field is generated using the software M\'erope \cite{Josien_Merope}. %The rest of the setting are identical to that in \cite{lu2021optimal}. 
For the right hand side, we take some function $f$ compactly supported in the box $\{-1, 0, 1\}^3$ with average zero, so that there exists some vector valued function $h$ such that $f=\nabla \cdot h$ for a function $h$ supported in the slightly larger box $Q_2$. The algorithm is compared against 
\begin{enumerate} 
	\item Solving the equation \eqref{IntroBC} with homogeneous (zero) Dirichlet boundary condition.
	\item Solving the equation \eqref{IntroBC} with modified correctors but without dipole or quadruple corrections, i.e. the boundary condition is given by \begin{equation}\label{eqn:numnopole}
		u_{\text{nc}}^{(L)}=(1+\phi_{i,M,L}^{(1)}\partial_i+\phi_{ij,M,L}^{(2)}\partial_{ij})u_{\text{hom},L} \ \mbox{ on } \partial Q_L.
	\end{equation}
	\item Solving the equation \eqref{IntroBC} with the boundary condition corrected up to first-order correctors and dipoles: \begin{equation}\label{eqn:numdipole}
		u_{\text{dp}}^{(L)}=(1+\phi_{k,M,L}^{(1)}\partial_k)\bigg(u_{\text{hom},L} +\Big(\int_{\R^d} h \cdot \nabla \phi_{i,M,L}^{(1)}\Big) \partial_i G_{\text{hom},L}\bigg) \ \mbox{ on } \partial Q_L.
\end{equation}\end{enumerate} 
We compare the numerical rate of $\lvert \nabla (u^{(2L)}- u^{(L)})(\frac{L}{2},\frac{L}{2},\frac{L}{2})\rvert$ and plot it for various $L$ for all four algorithms, and we show the numerical error plots for different choices of covariance functions.

\begin{figure}[ht]
	\begin{minipage}{3.3in}
	\includegraphics[width=3.1in]{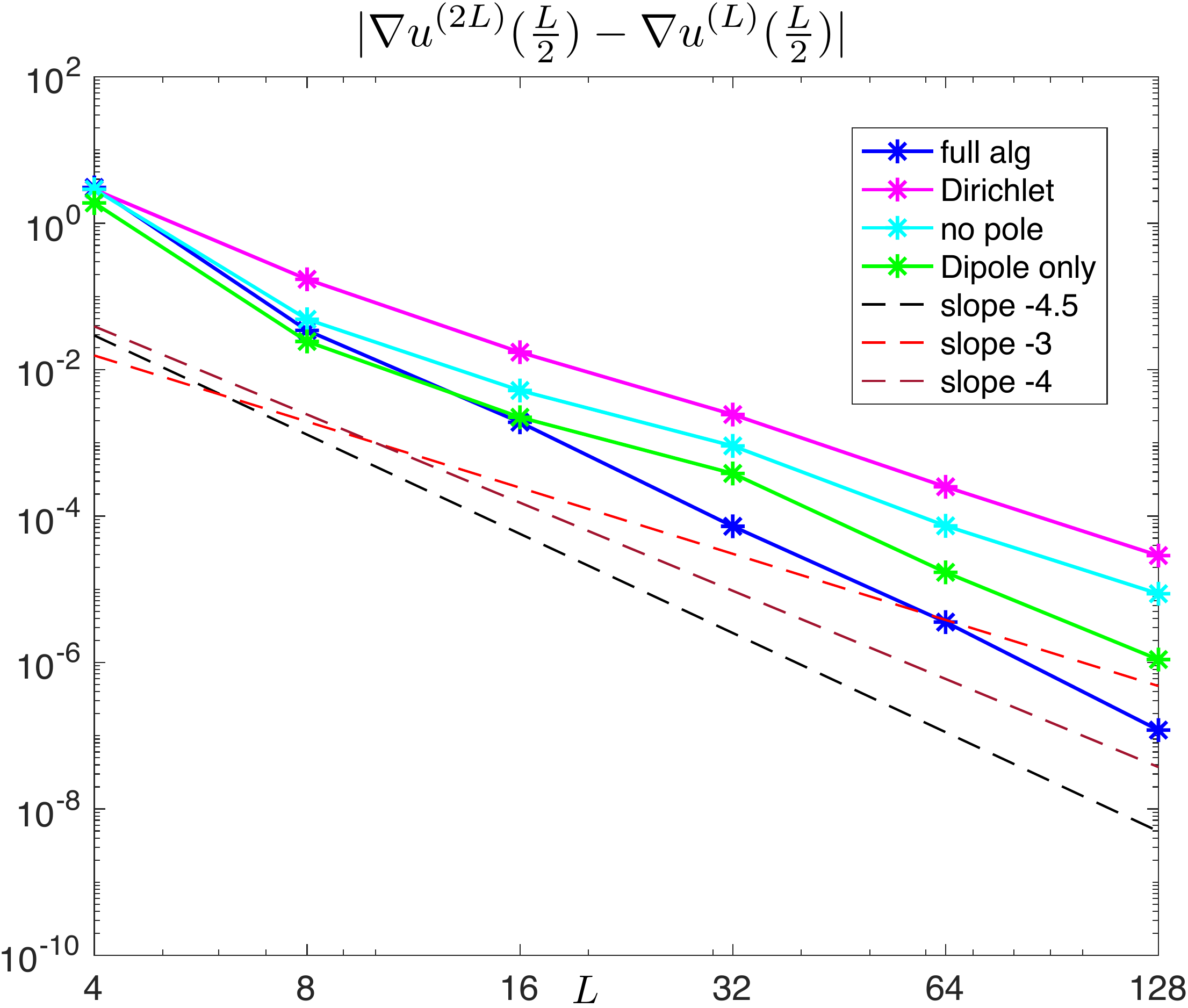}
	\end{minipage}
\begin{minipage}{3.3in}
	\includegraphics[width=3.25in]{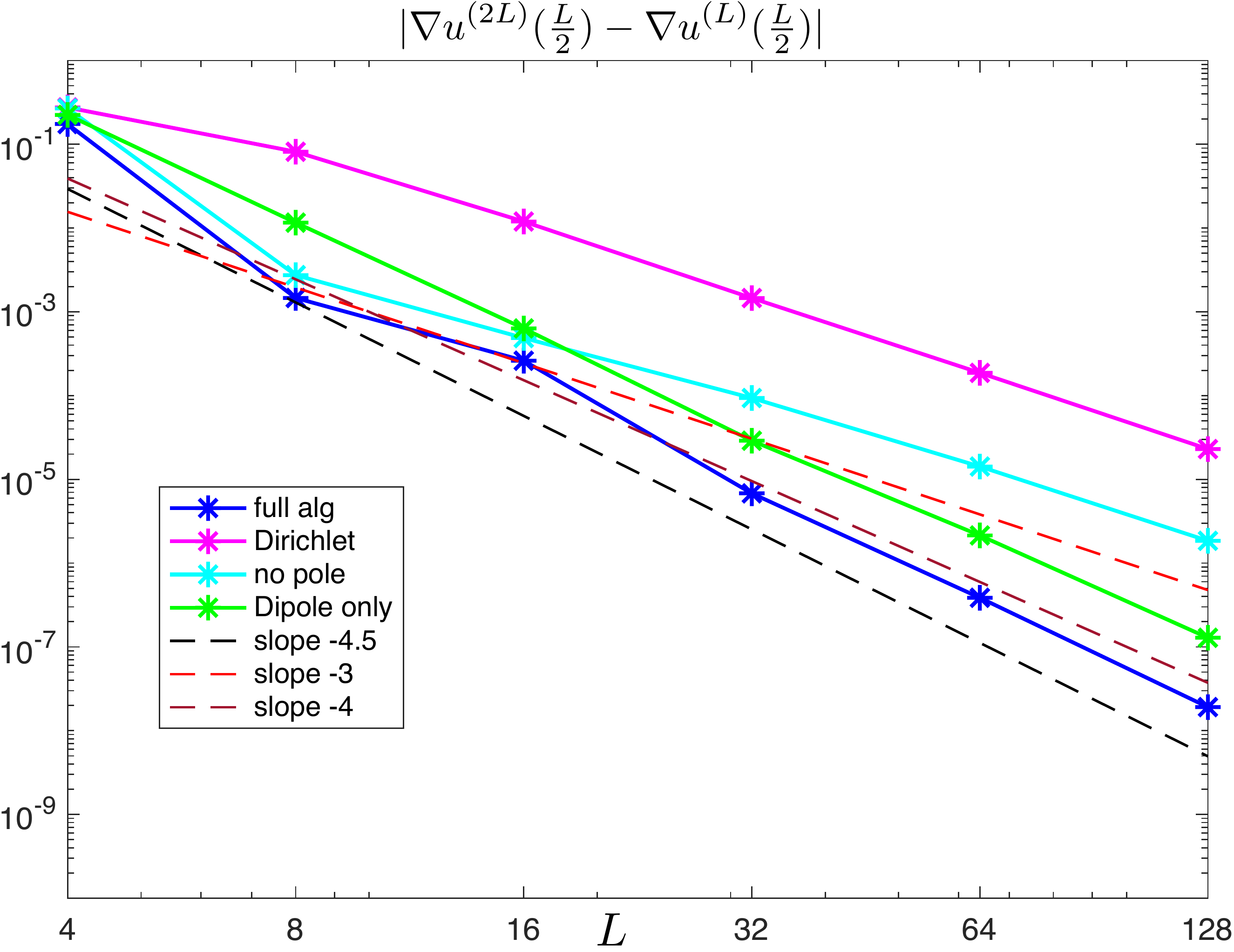}
\end{minipage}
\caption{Numerical convergence rates of $|\nabla u^{(2L)}(\tfrac{L}{2})- \nabla u^{(L)}(\tfrac{L}{2})|$ for four algorithms with covariance functions: left: $c(|x-y|)=\exp(-\tfrac{|x-y|^2}{8})$; right: $c(|x-y|)=(1+|x-y|/2)^{-5}$.}
\label{fig:goodfig}
\end{figure}

\medskip

The two figures in Figure \ref{fig:goodfig} show results that are similar to those in \cite{lu2021optimal} and consistent with what we present in Theorem \ref{thm:EffBdryAlg}, and our Algorithm \ref{alg:truealg} produces an error rate $O(L^{-4.5})$ that is optimal in the sense of Proposition \ref{OptimalScaling}. Meanwhile, the algorithm with dipoles only has a convergence rate of $O(L^{-4})$ while the other two algorithms behave as $O(L^{-3})$. 

\medskip

We would like to comment here, however, that the scenario with correlated coefficients is numerically more complicated and due to memory limitations it is difficult to generate large Gaussian fields and compute solutions of elliptic equations on large scales in dimension 3. Therefore, the increase in small scale correlations of the random field makes it difficult for computations to reach the asymptotic regime in which we observe the desired theoretical rates. This is manifested in Figure \ref{fig:notasyfig} below. We can observe that in both figures, the asymptotic regime is likely only reached around $L=128$, which is almost at our computational limit. If we further change our covariance functions to $\exp(-\tfrac{|x-y|^2}{32})$ or $(1+|x-y|/4)^{-5}$, which displays even worse numerical behavior than that in Figure \ref{fig:notasyfig}, then we are unable to observe the asymptotic regime even when $L=128$, as the full and Dipole only algorithms behave almost identically. Similar computational issue can also be observed in \cite[Figure 7(a)]{schneider2022representative}, where, despite obtaining a significantly small error, the bias for periodizing the ensemble does not appear to converge at the desired rate.
\begin{figure}[!ht]
	\begin{minipage}{3.3in}
		\includegraphics[width=3.1in]{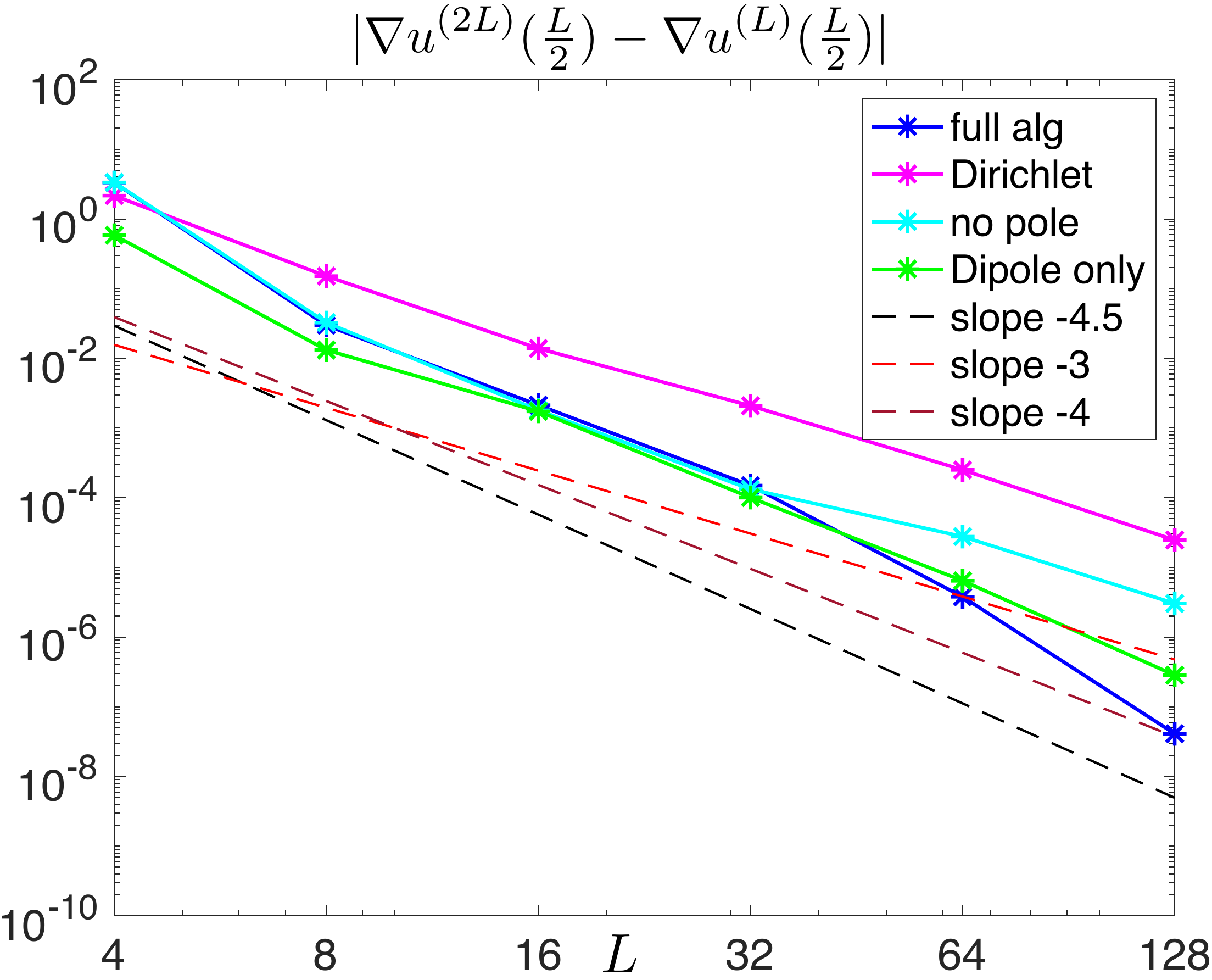}
	\end{minipage}
	\begin{minipage}{3.3in}
		\includegraphics[width=3.2in]{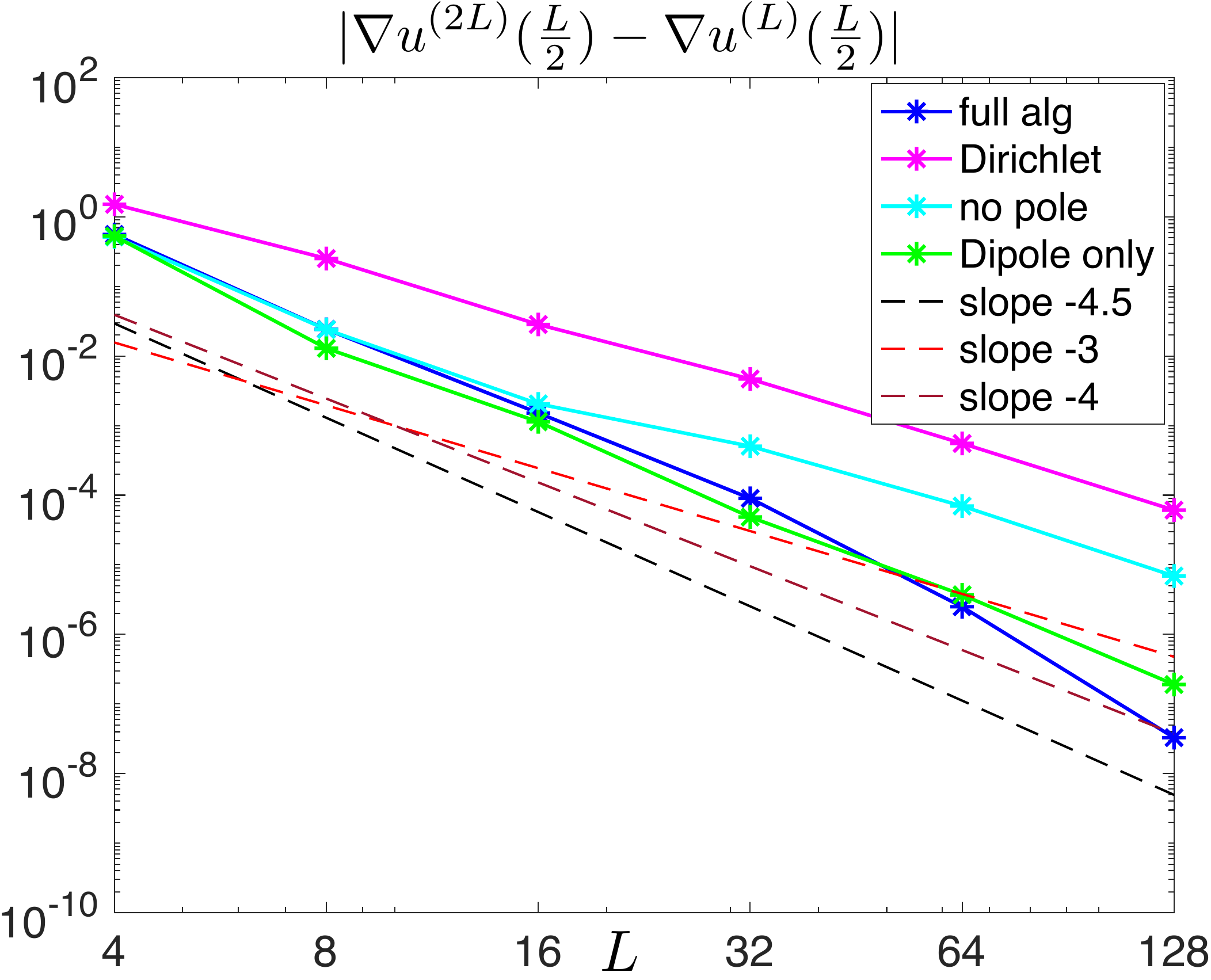}
	\end{minipage}
	\caption{Numerical convergence rates with covariance functions: left: $c(|x-y|)=\exp(-\tfrac{|x-y|^2}{16})$; right: $c(|x-y|)=(1+|x-y|/3)^{-5}$.}
	\label{fig:notasyfig}
\end{figure}

\medskip

\begin{figure}[!ht]
	\begin{minipage}{3.3in}
		\includegraphics[width=3.1in]{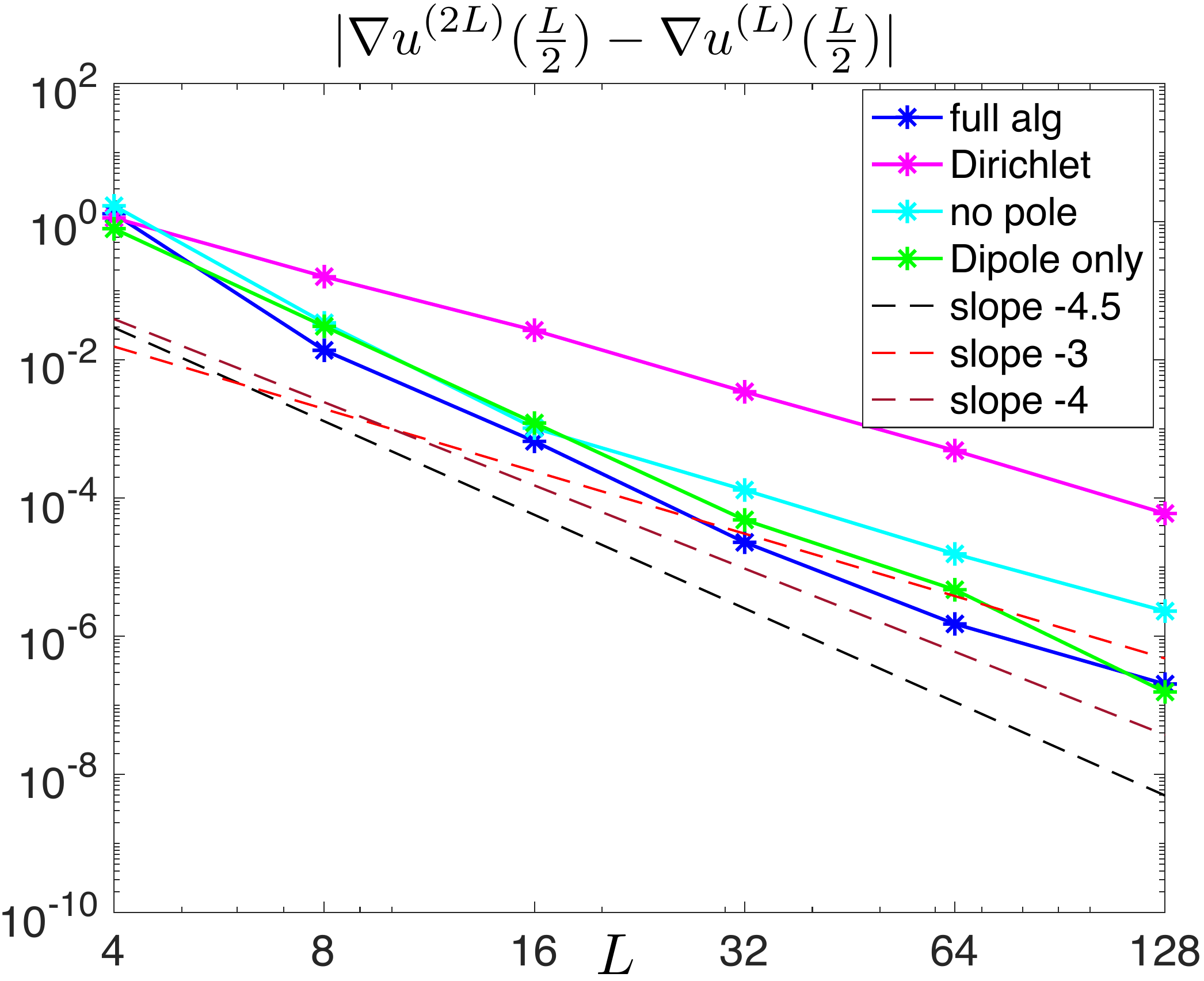}
	\end{minipage}
	\begin{minipage}{3.3in}
		\includegraphics[width=3.2in]{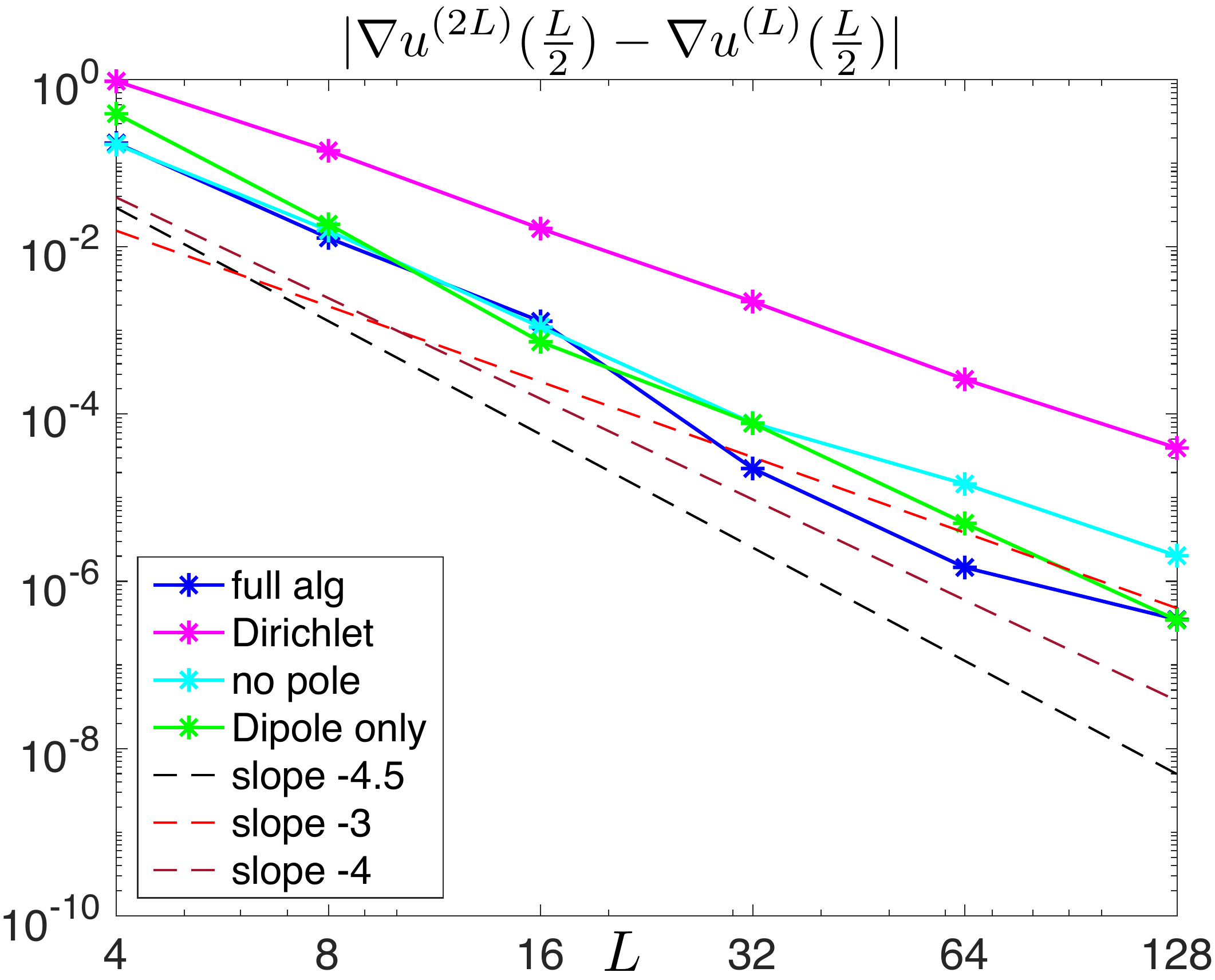}
	\end{minipage}
	\caption{Numerical convergence rates with covariance functions: left: $c(|x-y|)=(1+|x-y|)^{-\frac{5}{2}}$; right: $c(|x-y|)=(1+2|x-y|)^{-1.3}$.}
	\label{fig:badfig}
\end{figure}

We also present in Figure \ref{fig:badfig} numerical results when correlations are strong and standard LSI no longer holds. Based on above observations, however, we understand that due to failure of reaching asymptotic regime, we do not fully observe the behavior of numerical plots in the theoretical decay rate $\beta$ of $c$. In particular, error rates for the dipole algorithm and the full algorithm are both close to $O(L^{-4})$ in both scenarios. The results do suggest, however, that the benefit of adding dipole corrections diminish with a decreasing $\beta$.

\section{Structure of the proof}
\subsection{Strategy of the proof}\label{StrategySection}

It is reasonable to speculate that our proof strategy is a combination of \cite{lu2021optimal} and \cite{clozeau2021optimal}. Indeed, one may observe from \cite{lu2021optimal} that, thanks to the effective multipole observations in \cite{bella2020effective}, most of the remaining technical work lie in estimating the magnitude and the approximation errors of the second-order corrector $\phi^{(2)}$, which has been carried out in \cite{lu2021optimal} for distribution satisfying a finite-range of dependence. Therefore, our work generalizes the estimates of $\phi^{(2)}$ to the setting of MSLSI \eqref{msLSI}, by following the strategy based on optimal time decay estimates of the second-order semigroup. More precisely, we consider the parabolic equation 
\begin{equation}\label{strategy:naiveo2}
	\left\{
	\begin{array}{ll}
		\partial_{\tau} u^{(2)}_{ij}-\nabla\cdot a\nabla u^{(2)}_{ij}=0 & \text{ in $\mathbb{R}^d\times (0,+\infty)$}, \\
		u^{(2)}_{ij}(0)=\nabla\cdot(a\phi^{(1)}_{i}-\sigma^{(1)}_{i})e_j,& 
	\end{array}
	\right.
\end{equation} 
referred to as second-order semigroup, in line that it allows to disintegrate \eqref{intro:phi2} in form of 
\begin{equation}\label{StratDisinPhi2}
\nabla\phi^{(2)}_{ij}=\int_0^\infty \dd t\,\nabla u^{(2)}_{ij}(t,\cdot).
\end{equation}
Optimal time decay estimates on $t\mapsto \nabla u^{(2)}(t,\cdot)$ then allow to construct and deduce from \eqref{StratDisinPhi2} quantitative estimates on $(\phi^{(2)},\nabla\phi^{(2)})$ and error estimates for their massive approximation, as established in \cite{lu2021optimal}. Despite this similarity in strategy, our work also significantly differs from simply combining the works of \cite{lu2021optimal} and \cite{clozeau2021optimal} together. For distribution with a finite range of dependence \cite{lu2021optimal}, since the scale of locality, and as a result, CLT-type cancellations are spread through successive propagations of the parabolic equation \eqref{strategy:naiveo2}, the authors further decomposed $(\phi^{(1)},\sigma^{(1)})$ into the first-order semigroup and studied parabolic equations with a general approximately local initial condition to capture this sequential propagation behaviour. That strategy cannot be applied here, as functional inequalities of type \eqref{msPIp} cannot be easily iterated, and we fully exploit the structure of the initial data, in particular the corrector equations \eqref{intro:phi1} and \eqref{intro:sig1}. In the spirit of \cite{clozeau2021optimal}, we rather control the fluctuations of the second-order time dependent flux in Theorem \ref{SecondOrderDecayFlux} by monitoring the $\LL^2$-norm of its functional derivative that we compute in Lemma \ref{lem:funcderiv}. 
%While extending the strategy in \cite{clozeau2021optimal}, where the main difference lies in the structure of the initial data, further large-scale regularity results are needed which is the object of Section \ref{LargeScaleResultSec}. In addition, 
However, as opposed to \cite{clozeau2021optimal}, the initial data of \eqref{strategy:naiveo2} is no more given by the divergence of a stationary random field in $\LL^\infty$, which requires us to appeal to the massive term approximations $(\phi_M^{(1)},\sigma_M^{(1)})$ to make sense of the functional derivative in a pathwise way. For these results, we need to establish preliminary large-scale regularity estimates for elliptic and parabolic systems that we develop in Section \ref{LargeScaleResultSec}. The control of the fluctuations of the time dependent second-order flux then allows us to deduce optimal time decay estimates on $(u^{(2)},\nabla u^{(2)})$ in Corollary \ref{SecondOrderDecay}, which are themselves the key to construct the second-order correctors in Section \ref{ConstructionSecondCor} and to establish error estimates for their massive approximation in Section \ref{MassiveApproxSection}. We repeat here that the results in this section are stated for the second-order regime $\beta\wedge d >2$; in the first-order regime $\beta\wedge d \le 2$ only results in \cite{clozeau2021optimal} will be needed.

%, while at the same time, generalizes the approach of \cite{clozeau2021optimal} to right hand-sides that are not given by the divergence of a stationary uniformly bounded random field in $\LL^\infty$ but with exponential moment bounds and depending only implicitly on $a$.

%\smallskip

%However, . Indeed, if we consider the parabolic equation

%then we have the decomposition $\nabla\phi^{(2)}_{ij}=\int_{0}^{+\infty}\dd t\,\nabla u^{(2)}_{ij}(t,\cdot)$. In the finite range case of , In our setting, thanks to the more direct probabilistic tools \eqref{msLSI}, one can work with a \rhs that is only stationary with suitable moment bounds. On the other hand, when computing functional derivatives, which involves the derivatives of $(\phi^{(1)},\sigma^{(1)})$, there are substantial issues with making sense of the dual equations. To encounter this problem, and considering that we also need estimates for $(\phi_M^{(2)},\sigma_M^{(2)})$, we replace the correctors by their massive counterparts $(\phi^{(1)}_M,\sigma_M^{(1)})$, leading to the equation \eqref{2NDOrder:Eq1} that we mainly work with. 

%
\subsection{Large-scale regularity of elliptic and parabolic systems}\label{LargeScaleResultSec}
We state in this section large-scale regularity results for linear random elliptic and parabolic operators with coefficients satisfying Assumptions \ref{Gaussian} and \ref{DecayCor}. In the following results and the rest of the paper, we will generically use the notation $\mathcal{C}$ to denote a random field having uniform exponential moments in its argument, namely if $\mathcal{C}$ depends on parameters $\alpha_1,\cdots,\alpha_N$,
\begin{equation}\label{MomentBoundGeneric}
\exists \text{ deterministic $C,\gamma$ depending on $d, \beta$ and \eqref{RegA} s.t. $\sup_{\alpha_1,\cdots,\alpha_N}\mathbb{E}\big[\exp(\tfrac{1}{C}\mathcal{C}^\gamma(\alpha_1,\cdots,\alpha_N))\big]<\infty$.}
\end{equation}
Equivalently, we can check \eqref{MomentBoundGeneric} by controlling algebraic moments, that is 
\begin{equation}\label{BoundGenericAlgebraic}
\sup_{\alpha_1,\cdots,\alpha_N}\mathbb{E}\big[\mathcal{C}^p(\alpha_1,\cdots,\alpha_N)\big]^{\frac{1}{p}}\leq C p^{\frac{1}{\gamma}}\quad\text{for any $1\le p<\infty$,}
\end{equation}
where a proof of this equivalence can be found in \cite[Lemma 3.7]{ledoux2013probability}. Large-scale regularity have been considered since the work \cite{avellanedaLp} (see also the monologue \cite{shen2018periodic}) and was first considered in the stochastic setting in \cite{armstrong2016quantitative, gloria2014regularity}. The first result below, a generalization of \cite[Corollary 3]{gloria2014regularity}, states the large-scale $\cc^{0,1}$-regularity for random parabolic operators.  A proof can be found in \cite[Corollary $7$]{clozeau2021optimal} for the gradient estimate and a simple adaptation of \cite[Lemma 4.18]{lu2021optimal} provides the estimate for the solution itself. 
\begin{proposition}[Large-scale $\cc^{0,1}$-estimates for parabolic systems]\label{LSregSoutionItself}
%
%There exists a $\frac{1}{8}$-Lipschitz stationary random field $r_{\star}$ satisfying for a generic constant $C>0$
%
%\begin{equation}\label{MomentBoundrstar}
%
%\sup_{x\in\mathbb{R}^d}\mathbb{E}\big[\exp(\tfrac{1}{C}r^{\beta\wedge d}_\star(x))\big]<\infty,
%
%\end{equation}
%
Fix $(s,x)\in \mathbb{R}^{d+1}$. There exists a stationary random field $r_\star$ satisfying for any $\varepsilon>0$
\begin{equation}\label{Momentrstar}
\mathbb{E}\big[r^{dp}_\star\big]^{\frac{1}{p}}\lesssim_\varepsilon p^{\frac{1}{\gamma_{\beta,\varepsilon}}}\quad\text{with }\gamma_{\beta,\varepsilon}:=(1+\varepsilon\mathds{1}_{\beta=d})(\beta \wedge d),
\end{equation}
%\begin{itemize}
%
%\item[(i)] 
such that, given $r_\star(x)\leq r\leq R$ and (time independent) $f,g,h\in\LL^{2}(\mathbb{R}^d)$, the weak solution $u$ of
$$\partial_\tau u-\nabla\cdot a(\nabla u +g)=0 \quad\text{on $(s-R^2,s)\times \bb_R(x)$,}$$
satisfies 
\begin{equation}\label{LargeScaleC01Para}
\begin{aligned}
\fint_{s-r^2}^{s}\dd\tau\fint_{\bb_r(x)}\vert \nabla u(\tau,\cdot)+g\vert^2  \lesssim \fint_{s-R^2}^{s}\dd\tau\,\fint_{\bb_R(x)}\vert\nabla u(\tau,\cdot)+g\vert^2+ \sup_{\rho\in [r_\star,R]} \Big(\frac{R}{\rho}\Big)^{2}\fint_{\bb_\rho(x)}\Big|g-\fint_{\bb_\rho(x)}g\Big|^2. 
\end{aligned}
\end{equation}

In addition, if $g=0$ it holds
\begin{equation}\label{LargeScaleC01Paranondiv}
\fint_{s-\rho^2}^{s}\dd\tau\fint_{\bb_\rho(x)}\vert u(\tau,\cdot)\vert^2 \lesssim \fint_{s-R^2}^{s}\dd\tau\,\fint_{\bb_R(x)}\vert u(\tau,\cdot)\vert^2. %\\ &+ R^2\sup_{\rho\in [r,R]} \bigg[\fint_{\bb_\rho(x)}|g|^2 + R^2 \fint_{\bb_\rho(x)}|\nabla g|^2 + \Big|h-\fint_{\bb_\rho(x)}h\Big|^2\biggr]. 
\end{equation}
%

%\item[(ii)]For $M\geq 1$, let $v$ be the weak solution of
%$$(\tfrac{1}{M} - \nabla \cdot a \nabla)v = 0 \quad \text{in $\bb_R(x)$ for $R\leq \sqrt{M}$.}$$
%%\end{itemize}
%In the case where $r_\star(x)\leq \sqrt{M}$ and for any $r_\star(x)\leq \rho\leq R$, it holds
%%
%\begin{equation}\label{LargeScaleC01Massive}
%%
%\fint_{\bb_{\rho}}\vert(\tfrac{1}{\sqrt{M}}v,\nabla v)\vert^2\lesssim \fint_{\bb_R} \big\vert(\tfrac{1}{\sqrt{M}}v, \nabla v)\big\vert^2.
%%
%\end{equation}

%
\end{proposition}
%
% whereas a proof of \eqref{LargeScaleC01Massive} cannot be found directly in the literature in this form and a proof in provided in Section \nc{cite the proof}. %\lw{We probably don't need (ii) anymore. Come back to here later.}

%
\medskip

The next result states decay estimates of solutions of parabolic systems that we obtain from the large-scale regularity in Proposition \ref{LSregSoutionItself}. These estimates are in the spirit of \cite[Lemma 4]{clozeau2021optimal} and \cite[Lemma 2 \& 3 \& 4]{gloria2019quantitative}.
\begin{proposition}[Decay estimates for parabolic systems]\label{RegResultParabolic}
%\nc{I don't understand the scalings here. First, $(51)$ implies the second-item of $(47)$. Second, in view of its application to $w_2$, $f_r$ will be $\nabla\cdot v^{T}(0,\cdot)a e_j$ and $v^{T}$ satisfies $(50)$. So the bound on the gradient is unchanged but the bound on $F_r$ is weakened by a factor $\sqrt{T}$.} \lw{$(51)$ actually doesn't imply $(47)$, because when $|y|\ge \sqrt{T}\ge r$, $(51)$ is actually weaker. Also, for divergence form rhs with $F_r =\frac{1}{\sqrt{T}} v^T(0,\cdot)e_j$, it doesn't quite satisfy the second part of $(47)$ but only $(51)$, and that's why we need $(52)$. Finally will you explain why you use variable $y$ for assumptions and $x$ for conclusions? }
%
%\begin{equation}\label{AssumptionDecayParaLarge}
%
%\sup_{x\in\mathbb{R}^d}(\vert x\vert+r)^d\big(\vert f_r(x)\vert+(\vert x\vert+r)\vert\nabla f_r(x)\vert\big)<\infty.
%
%\end{equation}
%
Fix $r\geq 1$ and $f_r$ satisfying one of the following two assumptions 
\begin{itemize}
\item There exists $\delta\geq 0$ such that 
\begin{equation}\label{LSreg:Eq1}
		\vert f_r(x)\vert\leq \mathcal{C}(r,x)\frac{\log^\del(2+\frac{\vert x\vert}{r})}{(\vert x\vert+r)^d} \quad\text{for any $x\in\mathbb{R}^d$}.
	\end{equation}
\item There exists $\delta\geq 0$ and a vector field $F_r$ such that $f_r=\nabla\cdot F_r$ or $f_r=\nabla\cdot a F_r$ and satisfying
\begin{equation}\label{LSreg:Eq1p}
\vert F_r(x)\vert\leq \mathcal{C}(r,x)\frac{\log^\del(2+\tfrac{\vert x\vert}{r})}{(\vert x\vert+r)^d} \quad\text{and}\quad |\nabla F_r(x)| \le \C(r,x)\frac{\log^\del(2+\tfrac{\vert x\vert}{r})}{(\vert x\vert+r)^{d+1}}\quad\text{for any $x\in\mathbb{R}^d$}.
\end{equation}
\end{itemize}
Consider $v_r$ the weak solution of the backward-parabolic equation
\begin{equation}\label{EquationParaLargeScale}
\left\{
    \begin{array}{ll}
        \partial_{\tau}v_r+\nabla\cdot a\nabla v_r=f_r & \text{ in $(-\infty,0)\times \mathbb{R}^d$,} \\
        v_r(0,\cdot)\equiv 0.& 
    \end{array}
\right.
\end{equation}
For any $x\in\mathbb{R}^d$, we have :

\medskip

\begin{itemize}
	\item[(i)] If \eqref{LSreg:Eq1} holds, then 
\begin{equation}\label{LSreg:PointwiseParabolicnondiv}
	\big\vert\nabla v_r(t,x) \big\vert\leq \mathcal{C}(r,x)\frac{\log^{1+\del}(2+\tfrac{\vert x\vert}{r})}{(\vert x\vert+r)^{d-1}}\quad\text{for all $\sqrt{-t}\geq 2r_*(x)$}.
\end{equation}
\item[(ii)]If \eqref{LSreg:Eq1p} holds, then
\begin{equation}\label{LSreg:PointwiseParabolic}
|(\nabla v_r(t,x),\tfrac{1}{\sqrt{|t|}}v_r(t,x))| \le \C(r,x)  \frac{\log^{1+\del}(2+\tfrac{|x|}{r})}{(|x|+r)^d}\quad\text{for all $\sqrt{-t}\geq 2r_*(x)$}.
\end{equation}
In addition, if the second-item of \eqref{LSreg:Eq1p} is replaced by
\begin{equation}\label{LSreg:Eq1pq}
	|\nabla F_r(x)| \le \C(r,x) \frac{\log^\del (2+\tfrac{|x|}{r})}{\sqrt{T}(|x|+r)^d}\quad\text{for some $\sqrt{T}\ge r$ and any $x\in\mathbb{R}^d$ ,}
\end{equation} 
then it holds
\begin{equation}\label{LSreg:PointwiseParabolicalt}
	|\nabla v_r(t,x)| \le \C(r,x)  \frac{\log^{1+\del}(2+\tfrac{|x|}{r})}{(|x|+r)^d}\left(\frac{|x|+r}{\sqrt{T}}\vee 1 \right).
\end{equation}
\end{itemize}
We recall that $\mathcal{C}$ denotes a generic random constant satisfying \eqref{MomentBoundGeneric}.
\end{proposition}
An immediate consequence of Proposition \ref{RegResultParabolic} is the following decay estimates for solutions of elliptic systems with a massive term.
\begin{corollary}[Decay estimates for elliptic systems]\label{RegResultElliptic}
Let $M,r\geq 1$, $f_r$ satisfying either \eqref{LSreg:Eq1} or \eqref{LSreg:Eq1p} and $u_r$ be the weak solution of
\begin{equation}\label{LSreg:Eq3}
(\tfrac{1}{M}-\nabla\cdot a\nabla) u_r=f_r.
\end{equation}
For any $x\in\mathbb{R}^d$, we have :
\begin{itemize}
\item[(i)]If \eqref{LSreg:Eq1} holds, then 
\begin{equation}\label{LSreg:Eq4}
  \vert \nabla u_r(x)\vert\leq \mathcal{C}(r,x)\frac{\log^{1+\del}(2+\tfrac{\vert x\vert}{r})}{(\vert x\vert+r)^{d-1}}.
\end{equation}
\item[(ii)]If \eqref{LSreg:Eq1p} holds, then 
%
%\begin{equation}\label{LSreg:Eq17}
%%
% \big\vert\big((\vert x\vert+r)\nabla h_r(x),h_r(x)\big)\big\vert\leq \mathcal{C}(r,x)\frac{\log(2+\tfrac{\vert x\vert}{r})}{(\vert x\vert+r)^{d-1}}.
%%
%\end{equation}
%
%
\begin{equation}\label{LSreg:Eq18}
  \big\vert (\tfrac{1}{\sqrt{M}}u_r,\nabla u_r)(x)\big\vert \leq \mathcal{C}(r,x)\frac{\log^{1+\delta}(2+\tfrac{\vert x\vert}{r})}{(\vert x\vert+r)^d}.
\end{equation}
In addition, if the second-item of \eqref{LSreg:Eq1p} is replaced by \eqref{LSreg:Eq1pq}, then it holds
\begin{equation}\label{LSreg:Eq17}
	|\nabla u_r(x)| \le \C(r,x)\frac{\log^{1+\del}(2+\tfrac{|x|}{r})}{(|x|+r)^d}\left(\frac{|x|+r}{\sqrt{T}}\vee 1 \right).
\end{equation}
\end{itemize}
We recall that $\mathcal{C}$ denotes a general random constant satisfying \eqref{MomentBoundGeneric}.
\end{corollary}

\subsection{Optimal time decay estimates of the  $2^{{\rm nd}}$-order semigroup }\label{TimeDecaySemigroup}
%
%\subsection{Approximation estimates}
%
%\subsection{Semigroup estimates}
%
We state in this section the optimal time decay estimate of the second-order semigroup $u^{(2)}$ defined in \eqref{strategy:naiveo2}, which constitutes our main tool for the construction of sub-linear second-order correctors in Section \ref{ConstructionSecondCor} and for the massive approximation of second-order correctors in Section \ref{MassiveApproxSection}. As already mentioned in Section \ref{StrategySection}, our strategy follows from \cite{clozeau2021optimal} and first establishes moment bounds of averaging of the second-order time dependent flux. However, we do not work directly on $u^{(2)}$ but we appeal to a massive approximation of the first-order correctors, namely we define the second-order massive semigroup $u^{(2)}_M:=\{u^{(2)}_{ij,M}\}_{ij}$ for any indices $1\leq i,j\leq d$ and $M\geq 1$ as the weak solution of
\begin{equation}\label{2NDOrder:Eq1}
\left\{
    \begin{array}{ll}
        \partial_{\tau} u^{(2)}_{ij,M}-\nabla\cdot a\nabla u^{(2)}_{ij,M}=0 & \text{ in $(0,\infty)\times \mathbb{R}^d$}, \\
        u^{(2)}_{ij,M}(0)=\nabla\cdot(a\phi^{(1)}_{i,M}-\sigma^{(1)}_{i,M})e_j,& 
    \end{array}
\right.
\end{equation}
where we recall that $(\phi^{(1)}_M,\sigma^{(1)}_M)$ are defined in \eqref{eq:1stOrderMassCor} and \eqref{eq:1stOrderMassFlux}. From $u^{(2)}_M$ we define the second-order time-dependent flux $q^{(2)}_M:=\{q^{(2)}_{ij,M}\}_{ij}$ for any indices $1\leq i,j\leq d$ and $M\geq 1$
\begin{equation}\label{2NDOrder:Eq8}
q^{(2)}_{ij,M}(T,\cdot):=(a\phi^{(1)}_{i,M}-\sigma^{(1)}_{i,M})e_j+\int_{0}^T \dd s\, a\nabla u^{(2)}_{ij,M}(s,\cdot)\text{\quad for any $T<\infty$,}
\end{equation}
and the second-order time-dependent corrector 
\begin{equation}\label{TimeDepSecondCor}
\phi^{(2)}_{ij,M}(T,\cdot)=\int_{0}^T u^{(2)}_{ij,M}(s,\cdot)\quad\text{for any $T<\infty$.}
\end{equation}
Considering \eqref{2NDOrder:Eq1} rather than \eqref{strategy:naiveo2} has two advantages. First, the massive correctors $(\phi^{(1)}_M,\sigma^{(1)}_M)$, as opposed to their massive-less counterparts, are well-defined for every realization $a$ so that we can compute their functional derivatives defined in \eqref{DefDerivative}. Second, $(\phi^{(1)}_M,\sigma^{(1)}_M)$ are more regular and satisfy the energy estimate \eqref{LocalizedMassiveTerm} which, combined with \eqref{GOlm1} and \eqref{EnergyEstiAppendix:Eq1}, controls the energy of $u^{(2)}_M$ and will be crucially used to justify the functional derivative of $q^{(2)}_M$ in Lemma \ref{lem:funcderiv}. The optimal time decay estimate of $u^{(2)}$ are then deduced from the optimal time decay estimate of $u^{(2)}_{M}$ by passing to the limit as $M\uparrow\infty$, making use of the convergence of the massive correctors $(\phi^{(1)}_M,\sigma^{(1)}_M)$ towards their massive-less counterparts $(\phi^{(1)},\sigma^{(1)})$ (see Lemma \ref{lem:1stCorApprox}).

\medskip

We now state the main result of this section, namely optimal fluctuation estimates for $(q^{(2)}_M(T,\cdot),\nabla\phi^{(2)}_M(T,\cdot))$ in the spirit of \cite[Proposition $1$]{clozeau2021optimal}. In the following result and the rest of the paper, we define for any $r\geq 1$ and function $f$,
\begin{equation}\label{GaussianKernel}
g_r:=(2\pi r^2)^{-\frac{d}{2}}\exp(-\tfrac{\vert\cdot\vert^2}{r^2})\quad\text{and}\quad f_r:=f\star g_r.
\end{equation}
%
%We state in this section the optimal time decay estimates on the $2^{\text{nd}}$-semigroup, which is our main tool to prove the massive approximation results in Corollary \ref{Cor:Approx2ndCor} and Lemma \nc{Give a name} as well as for the construction of sub-linear $2^{\text{nd}}$-order correctors in Corollary \ref{Cor:Growth2ndCor}. 

%
%\medskip
%

%\nc{We probably want to define it before while explaining the strategy}. As mentioned in Section \nc{strategy of the proof}, we appeal to the massive approximation and define the $2^{\text{nd}}$-order massive-semigroup 
%

%
%where we recall that $\phi^{(1)}_{e,M}$ and $\sigma^{(1)}_{e_i,M}$ denote the first-order massive corrector and flux corrector, defined as the unique stationary distributional solution of \eqref{eq:1stOrderMassCor} and \eqref{eq:1stOrderMassFlux} respectively.
%
%The semigroup $u^{(2)}_{e,e'}$ defined in \eqref{strategy:naiveo2} is obtained by replacing $\phi^{(1)}_{e,M}$ and $\sigma^{(1)}_{e,M}$ in \eqref{2NDOrder:Eq1} by $\phi^{(1)}_e$ and $\sigma^{(1)}_e$. \nc{but should be defined earlier}

%
%\medskip
%

%Our strategy follows \cite{clozeau2021optimal} that we extend to control the $2^{nd}$-order semigroup. As in \cite[Theorem 1]{clozeau2021optimal}, we first control the fluctuation of averages of the massive second-order time-dependent flux $q^{(2)}_{M}$ defined as

%and the time dependent massive second-order corrector
%
\begin{theorem}[Fluctuations of the $2^{\text{nd}}$-order time-dependent flux]\label{SecondOrderDecayFlux}
Let $1\le r<\sqrt{T}<\sqrt{M}$. We have for any $x\in\mathbb{R}^d$ 
\begin{equation}\label{SecondOrderDecayFlux:Eq}
\big\vert \big(q^{(2)}_{ij,M}(T,\cdot)-\mathbb{E}[q^{(2)}_{ij,M}(T,\cdot)] , \nabla \phi^{(2)}_{ij,M}(T,\cdot)\big)_r(x)\big\vert\leq\mathcal{C}(r,T,x)\log^2(\tfrac{\sqrt{T}}{r}) r^{-\frac{d}{2}}\mu_{\beta}(T),
\end{equation}
with
\begin{equation}
	\label{eq:defmubeta}
\mu_\beta(T) := \left\{
    \begin{array}{ll}
        T^{\frac{1}{2}+\frac{d}{4}-\frac{\beta}{4}} & \text{if $\beta<d$}, \\
        \sqrt{T}\log^{\frac{1}{2}}(T) & \text{if $\beta=d$,}\\
        \sqrt{T} & \text{if $\beta>d$,}
    \end{array}
\right.\end{equation}
where we recall that $\mathcal{C}$ satisfies \eqref{MomentBoundGeneric}.
\end{theorem}
As a direct consequence, we deduce the optimal time decay estimates of the second-order semigroup.
\begin{corollary}[Optimal time decay of $2^{\text{nd}}$-order semigroups]\label{SecondOrderDecay}
For any $\sqrt{M}\geq \sqrt{T} \ge 1$, we have
\begin{equation}\label{2NDOrder:Eq2}
T^{-\frac{1}{2}}\big\vert (u^{(2)}_{ij}(T,x),u^{(2)}_{ij,M}(T,x))\big\vert+\big\vert(\nabla u^{(2)}_{ij}(T,x),\nabla u^{(2)}_{ij,M}(T,x))\big\vert\leq \mathcal{C}(T,x) T^{-1-\frac{d}{4}}\mu_{\beta}(T),
\end{equation}
where we recall that $\mathcal{C}$ satisfies \eqref{MomentBoundGeneric} and $u^{(2)}$ as well as $u^{(2)}_M$ are respectively defined in \eqref{strategy:naiveo2} and \eqref{2NDOrder:Eq1}.
\end{corollary}
The proof of Theorem \ref{SecondOrderDecayFlux} is based on a representation formula for the functional derivative of $q^{(2)}_{M}$ defined in \eqref{DefDerivative}, treating apart the contribution in $\sigma^{(1)}_{M}$ in \eqref{2NDOrder:Eq8} with more direct techniques. Due to the lack of regularity of the perturbations $\delta a$ in \eqref{DefDerivative}, we establish the representation formula for a regularized version of $q^{(2)}_{M}$, see \eqref{RegularizedQ2}, where the initial data is smoothed uniformly with respect to the randomness by convolving $a$ with a smooth kernel $\chi_\varepsilon$. This procedure avoids singularities at time $t=0$ and makes all integrals converge in a Lebesgue sense. This regularization step does not deteriorate the estimate \eqref{SecondOrderDecayFlux:Eq} thanks to the regularity assumption \eqref{SmoothnessCoef} which allows us to pass to the limit as $\varepsilon\downarrow 0$. 
\begin{lemma}[Functional derivative of the $2^{\text{nd}}$-order time dependent flux]\label{lem:funcderiv}
Let $M\geq 1$, $\varepsilon\in (0,1)$ and define the kernel $\chi_\varepsilon:=(2\pi\varepsilon)^{-\frac{d}{2}}\exp(-\tfrac{\vert\cdot\vert^{2}}{\varepsilon})$. We consider the regularized second-order time depend flux 
\begin{equation}\label{RegularizedQ2}
q^{(2)}_{ij,M,\varepsilon}(T,\cdot):=(a\phi^{(1)}_{i,M}-\sigma^{(1)}_{i,M})e_j+\int_{0}^{T}\dd s\,a\nabla u^{(2)}_{ij,M,\varepsilon}(s,\cdot)\quad\text{for any $T<\infty$,}
\end{equation}
with $a_\varepsilon:=\chi_\varepsilon\star a$ and
\begin{equation}\label{2NDOrder:Eq1Regularized}
\left\{
    \begin{array}{ll}
        \partial_{\tau} u^{(2)}_{ij,M,\varepsilon}-\nabla\cdot a\nabla u^{(2)}_{ij,M,\varepsilon}=0 & \text{ in $(0,+\infty)\times \mathbb{R}^d$}, \\
        u^{(2)}_{ij,M,\varepsilon}(0)=\nabla\cdot(a_\varepsilon\,\phi^{(1)}_{i,M}-\sigma^{(1)}_{i,M})e_j.& 
    \end{array}
\right.
\end{equation}
For any $\ell\geq 1$, $1\leq r<\sqrt{T}$ and unit vector $\upsilon\in \mathbb{R}^d$, we have
\begin{equation}\label{eqn:funcderiv}
		\begin{aligned}
			&\partial^{\mathrm{fct}}_{x,\ell}\big(q^{(2)}_{ij,M,\varepsilon}(T,\cdot)+\sigma^{(1)}_{i,M}e_j\big)_r(0)\cdot \upsilon\\
			=& \int_{\bb_\ell(x)} \vert\phi_{i,M}^{(1)} g_r e_j\otimes \upsilon\vert +\int_{\bb_{\ell}(x)} g_r\bigg\vert\int_{0}^T\dd s\,\nabla u^{(2)}_{ij,M,\varepsilon}(s,\cdot)\bigg\vert+\int_{\bb_\ell(x)}\,\big|(\nabla v^T(0,\cdot) \otimes\phi^{(1)}_{i,M}e_j)\star \chi_\varepsilon\big|\\
			&+ \int_{\bb_\ell(x)}\bigg |\int_0^T \dd s\,\nabla v^T(s,\cdot) \otimes \nabla u^{(2)}_{ij,M,\varepsilon}(s,\cdot) \bigg|+\int_{\bb_\ell(x)}\big|\nabla w_2\otimes (e_i+\nabla \phi^{(1)}_{i,M})\big|\\
			& + \int_{\bb_\ell(x)}\big| (v^T(0,\cdot)e_j-\tfrac{1}{M}w_1)\otimes (e_i+\nabla \phi^{(1)}_{i,M}) \big|,
		\end{aligned}
	\end{equation}
where $v^T$ solves the backward parabolic equation
\begin{equation}\label{eq:vT}
		\left\{ \begin{aligned}
			& \partial_\tau  v^T +\nabla \cdot a \nabla v^T= \nabla \cdot a g_r \upsilon \text{\quad in $(-\infty,T)\times \mathbb{R}^d$},\\
			& v^T(T,\cdot) \equiv 0 ,
		\end{aligned}   \right.
	\end{equation}
$w_1$ solves
\begin{equation}\label{eqn:phi1dual}
(\tfrac{1}{M}-\Delta) w_1 =v^{T}(0,\cdot)e_j,
	\end{equation}
and $w_2$ solves
\begin{equation}\label{eqn:2ndphi1dual}
(\tfrac{1}{M}-\nabla \cdot a \nabla) w _2= g_r a \upsilon\cdot e_j -\nabla v^{T}(0,\cdot)\cdot a_\varepsilon\,e_j+\nabla\cdot v^T(0,\cdot)ae_j-\tfrac{1}{M}\nabla\cdot \big((a-\mathrm{Id})w_1\big).
\end{equation}
\end{lemma}
Finally, we derive pointwise estimates on $w_1$ satisfying  \eqref{eqn:phi1dual} using the explicit formulation in terms of the whole space Green function of the operator $\frac{1}{M}-\Delta$. The below estimate on $w_1$ allows us to estimate $w_2$, which in turn is used to establish the fluctuation bounds \eqref{SecondOrderDecayFlux:Eq}.
\begin{lemma}[Pointwise bounds for the massive Laplacian operator]\label{lem:ptwisew1}
 Let $1\le r\le \sqrt{M \wedge T}$ and $w_1$ be the solution of \eqref{eqn:phi1dual}. It holds %\lw{with improvements we can again remove the 1 term in max. By a careful examination of the proof, we can also remove the $\sqrt{M}$ factor and replace it by $|x|+r$ - which will justify that $\tfrac{1}{M}w_1$ is really a perturbation term.}
	\begin{equation}\label{eqn:bdw1}
		\big|(\nabla w_1, \tfrac{1}{\sqrt{M}}w_1)(x)\big| \le \mathcal{C}(r,x) \sqrt{MT}\,\frac{\log(2+\frac{|x|}{r})}{(|x|+r)^d}\quad\text{for any $x\in\mathbb{R}^d$}.
	\end{equation}
\end{lemma}

\subsection{Construction of sub-linear $2^{{\rm nd}}$-order correctors}\label{ConstructionSecondCor}
We state in this section the existence of sub-linear second-order correctors $(\phi^{(2)},\sigma^{(2)})$, defined in \eqref{intro:phi2}, together with optimal estimates. The construction is based on the optimal time decay estimates of the second-order semigroup in Corollary \ref{SecondOrderDecay} together with the representation formula for any indices $1\leq i,j\leq d$
\begin{equation}\label{2NDOrder:Eq3}
\nabla\phi^{(2)}_{ij}=\int_0^\infty\dd t\,\nabla u^{(2)}_{ij}(t,\cdot),
\end{equation}
that it is formally obtained by integrating \eqref{strategy:naiveo2} in time and the uniqueness in \eqref{intro:phi2}. 
%We state in this section the existence and properties of second-order correctors and convergence rate for their massive approximation. As a consequence, we construct the second order correctors by time integration. Indeed, from \eqref{2NDOrder:Eq2}, $\int_{0}^{\infty} \dd t\, \nabla u^{(2)}_{e_i,e_j}(t,\cdot)$ in well defined in $\LL^2(\Omega\times \mathbb{R}^d)$ and curl free so that there exists a potential $\phi^{(2)}_{ij}\in  \LL^2(\Omega\times\mathbb{R}^d)$ such that
%
%\begin{equation}\label{2NDOrder:Eq3}
%\nabla\phi^{(2)}_{ij}=\int_{0}^{+\infty}\dd t\,\nabla u^{(2)}_{ij}(t,\cdot).
%\end{equation}
%
%Therefore, after taking spatial derivatives, and using the equation \eqref{2NDOrder:Eq1}, we can show that the equation \eqref{intro:phi2} holds in the distributional sense.
%
%\begin{equation}\label{2NDOrder:Eq4}
%-\nabla\cdot a\nabla\phi^{(2)}_{ij}=\nabla\cdot (a\phi^{(1)}_i-\sigma^{(1)}_{i})e_j.
%\end{equation}
%The argument is rigorous when we apply \eqref{2NDOrder:Eq2}. As a byproduct of Theorem \ref{SecondOrderDecayFlux}, we can also show that $\phi^{(2)}_{ij}$ is the unique (up to an additive constant) solution of \eqref{intro:phi2} and we further establish optimal stochastic estimate on $\phi^{(2)}_{ij}$ directly deduced from \eqref{2NDOrder:Eq2}.
%
\begin{corollary}[Sub-linear $2^{\text{nd}}$-order correctors] \label{Cor:Growth2ndCor}
For any indices $1\leq i,j\leq d$, there exist random fields $\phi^{(2)}_{ij}$ and $\sigma^{(2)}_{ij}=\{\sigma^{(2)}_{ijkn}\}_{k,n}$ belonging to $\mathrm{H}^1_{\mathrm{loc}}(\mathbb{R}^d)$ such that
\begin{itemize}
\item[(i)]$\phi^{(2)}_{ij}$ and $\sigma^{(2)}_{ij}$ are the unique (up to an additive constant) random fields with stationary gradient, $\mathbb{E}[\nabla \phi^{(2)}_{ij}]=\mathbb{E}[\nabla \sigma^{(2)}_{ij}]=0$ and satisfy in the distributional sense
\begin{equation}\label{EqutionPhi2Theorem}
-\nabla\cdot a\nabla\phi^{(2)}_{ij}=\nabla\cdot (a\phi^{(1)}_i-\sigma^{(1)}_i)e_j\quad\text{and}\quad\left\{
    \begin{array}{ll}
      \nabla\cdot \sigma^{(2)}_{ij}=q_{ij}^{(2)}:=a\nabla\phi^{(2)}_{ij}+(\phi^{(1)}_i a-\sigma_i^{(1)})e_j, &  \\
         -\Delta \sigma^{(2)}_{ij}=\nabla\times q_{ij}^{(2)},& 
    \end{array}
\right.
\quad\text{a.~s.}
\end{equation}
with for any $x\in\mathbb{R}^d$,
\begin{equation}\label{BoundPhi2Grad}
\vert\nabla (\phi^{(2)}_{ij}, \sigma_{ij}^{(2)})(x)\vert \le \C(x).
\end{equation}
\item[(ii)] For any $r>0$ and $x\in\mathbb{R}^d$
\begin{equation}\label{BoundPhi2Fluctu}
\vert\nabla(\phi^{(2)}_{ij},\sigma_{ij}^{(2)})_r(x)\vert \le \C(x,r) r^{-\frac{d}{2}}\mu_{\beta}(r^2),
\end{equation}
where we recall that $\mu_\beta$ is defined in \eqref{eq:defmubeta} and the subscript $r$ stands for the convolution \eqref{GaussianKernel}.

\medskip

\item[(iii)] For any $x\in\mathbb{R}^d$
\begin{equation}\label{BoundPhi2} \begin{aligned}
\bigl \vert \bigl(\phi^{(2)}_{ij}(x)-\phi^{(2)}_{ij}(0), \sigma^{(2)}_{ij}(x)-\sigma^{(2)}_{ij}(0) \bigr) \bigr\vert &\le \C(x)\nu_{d,\beta}(\vert x\vert), 
\end{aligned}\end{equation}
with 
\begin{equation}\label{DefNuBeta}
 \nu_{d,\beta}(\vert x\vert):=\left\{
    \begin{array}{ll}
        (\vert x\vert+1)^{2-\frac{\beta}{2}}& \text{ for $2<\beta<4$ and $\beta<d$}, \\
        \sqrt{\vert x\vert+1} & \text{ for $d=3$ and $\beta>3$,}\\
        \log^{\frac{1}{2}}(\vert x\vert+2)\sqrt{\vert x\vert +1} & \text{ for $\beta=d=3$,}\\
        \log(\vert x\vert+2) & \text{ for $\beta=4$, $d>4$ and $\beta>4$, $d=4$,}\\
        \log^{\frac{3}{2}}(\vert x\vert+2) & \text{ for $\beta=d=4$,}\\
        1 &\text{ for $\beta>4$ and $d>4$.}
    \end{array}
\right.
\end{equation}
\end{itemize}
We recall that $\mathcal{C}$ denotes a random constant that satisfies \eqref{MomentBoundGeneric}.
\end{corollary}

\subsection{Massive approximation of correctors}\label{MassiveApproxSection}
We state in this section quantitative estimates on the approximation of first- and second-order correctors by their massive counterparts. The following results are based on the representation formula using the first- and second-order semigroups.
%The semigroup estimates in Theorem \ref{SecondOrderDecayFlux} also allows us to construct approximations for the second-order correctors. Indeed, $\phi^{(2)}_M$, the solution of \eqref{eq:2ndOrderMassCor}, also allows the semigroup decomposition 
%
%begin{equation}\label{eq:2ndMassCorParaDecom}
%
%\phi_{ij,M}^{(2)} = \int_0^\infty \dd t \, e^{-\frac{t}{M}}u_{e_i,e_j,M}^{(2)}(t,\cdot).
%	
%\end{equation}
%
%The next results state the approximation error between $(\phi^{(1)}_i,\sigma^{(1)}_{ijk},\phi^{(2)}_{ij})$ with their massive counterparts $(\phi^{(1)}_{i,M},\sigma^{(1)}_{ijk,M},\phi^{(2)}_{ij,M})$. 
We specify our results for $d=3$, which is the numerically pertinent setting, the approximation rates get better if we increase $d$ and use repeated Richardson extrapolation as in \cite{gloria2015corrector, clozeau2021optimal}.
\begin{corollary}[Massive approximation of $2^{\text{nd}}$-order correctors]\label{Cor:Approx2ndCor}
For any $R,M\geq 1$, it holds
\begin{equation}\label{eq:phi2ptws}
	|\phi_M^{(2)}(x)| \le \C(x)M^{1-\frac{\beta\wedge 3}{4}}(1+\mathds{1}_{\beta=3} \log^\frac{1}{2}M) ,
\end{equation}
and
\begin{equation}\label{eq:2ndCorApproxBd}
\bigg(\int_{\R^d} \eta_R(x-\cdot) \big|\nabla \big(\phi_{M}^{(2)}-\phi^{(2)} \big)\big|^2\bigg)^\frac{1}{2} \le \C(R,M,x) M^{\frac{1}{2}-\frac{\beta\wedge 3}{4}}\big(1+ \mathds{1}_{\beta=3}\log^\frac{1}{2} M\big)\log M,
\end{equation} 
We recall that $\eta_{R}:=R^{-d}\exp(-\tfrac{\vert\cdot\vert}{R})$ and $\mathcal{C}$ denotes a random constant that satisfies \eqref{MomentBoundGeneric}.
\end{corollary}
Below we present an auxiliary lemma, which shows when $\beta\wedge d >2$, the proxy correctors $(\phi_M^{(1)},\sigma_M^{(1)})$ are stationary themselves and are good approximations of $(\phi^{(1)},\sigma^{(1)})$. We drop the proof as the ideas are essentially contained in \cite[Corollary $5$]{clozeau2021optimal} and in the proof of Corollary \ref{Cor:Approx2ndCor}.
\begin{lemma}[Massive approximation of $1^{\text{st}}$-order correctors]\label{lem:1stCorApprox}
For any $R,M\ge 1$, the first-order massive corrector satisfies 
\begin{equation}\label{eq:phiptws}
	|(\phi_M^{(1)},\nabla\phi^{(1)}_M,\sigma^{(1)}_M,\nabla\sigma^{(1)}_M)(x)| \leq \C(x), 
\end{equation}
and 
\begin{equation}
 \begin{aligned}\label{eq:1stCorApprox}
		&\bigg(\int_{\R^d} \eta_R(x-\cdot) \big|\big( \phi_M^{(1)}-\phi^{(1)}, \sigma_M^{(1)}-\sigma^{(1)}, \sqrt{M}\nabla (\phi_M^{(1)}-\phi^{(1)}, \sigma_M^{(1)}-\sigma^{(1)} )\big)\big|^2\big| \bigg)^\frac{1}{2} \\ 
		&\le \C(R,x)\,M^{{\frac{1}{2}-\frac{\beta\wedge 3}{4}}}(1+\mathds{1}_{\beta=3}\log^\frac{1}{2} M)
	\end{aligned}
\end{equation}
We recall that $\eta_{R}:=R^{-d}\exp(-\tfrac{\vert\cdot\vert}{R})$ and $\mathcal{C}$ denotes a random constant that satisfies \eqref{MomentBoundGeneric}.
\end{lemma}
\section{Proofs}
\subsection{Proof of Theorem \ref{thm:EffBdryAlg}: Artificial boundary approximation for whole-field problem}
We only provide a sketch of the proof here with the main ideas, for more detailed arguments we refer the reader to the proof of \cite[Theorem 1.2]{lu2021optimal}. We split the proof into four steps. In the first step, we define the random radius $r_{\star\star}$ from the definition of the first and second-order correctors as well as the moment bounds in Corollary \ref{Cor:Growth2ndCor}. In the second step, we appeal to the effective multipoles expansion theory in \cite{bella2020effective}. In the third step, we derive error estimates for the approximation of the correctors by their massive counterparts satisfying homogeneous Dirichlet boundary conditions. In the final step, we conclude by combining the whole.

\medskip

{\sc Step 1. Definition of $\rss$. }The second-order correctors $(\phi^{(2)},\sigma^{(2)})$ are not stationary themselves, but rather grow at a rate slightly worse than $2-\gamma_\beta$ (where $\gamma_\beta := \frac{\beta\wedge 3}{2}$) away from the origin, see \eqref{BoundPhi2}. Meanwhile, $(\phi^{(1)},\sigma^{(1)})$ are stationary but not uniformly bounded. Hence, for any fixed $\varepsilon>0$ small such that $\gamma_{\beta}(1-\varepsilon)>1$, we may define a random radius $\rss \ge 1$, starting from which the correctors have the desired growth rate :
\begin{equation}\label{rStarStarSecond}
	\frac{1}{r^2}\bigg(\fint_{\bb_r}\Big\vert (\phi^{(2)},\sigma^{(2)})-\fint_{\bb_r}(\phi^{(2)},\sigma^{(2)})\Big\vert^2\bigg)^{\frac{1}{2}}\leq \Big(\frac{\rss}{r}\Big)^{\gamma_\beta(1-\varepsilon)}\quad \text{for $r\geq \rss$,}
\end{equation}
and
\begin{equation}\label{rStarStarFirst}
	\frac{1}{r}\bigg(\fint_{\bb_r}\vert (\phi^{(1)},\sigma^{(1)})\vert^2\bigg)^{\frac{1}{2}}\leq \Big(\frac{r_{\star\star}}{r}\Big)^{1-\varepsilon} \quad \text{for $r\geq \rss$}.
\end{equation}
In the case $ \beta\wedge d\le 2$, the above definition of $\rss$ does not make sense, instead it collapses to $r_\star$ (defined in Proposition \ref{LSregSoutionItself}).  We would like to comment that $r_\star$ is the characteristic scale that captures the sublinear growth condition of  $(\phi^{(1)},\sigma^{(1)})$, which in turns yields large-scale elliptic and parabolic $C^{0,1}$ estimates, see \cite{gloria2014regularity}. Meanwhile, $\rss$ captures the growth of $(\phi^{(2)},\sigma^{(2)})$, which controls the two-scale expansion error up to second-order (as opposed to first-order for $r_\star$).

\medskip

{\sc Step 2.  Effective multipoles using whole-space correctors. }This step is completely deterministic, only relying on \cite[Theorem 2]{bella2020effective} and \cite[Corollary 2.2]{lu2021optimal}, which are applicable since their assumptions follow directly from the quantitative ergodicity of the coefficient field \eqref{msLSI}. The result can be stated as follows : Let $\hat{u}$ be defined via
\begin{equation*}
	\left\{
	\begin{array}{ll}
		-\nabla\cdot a\nabla \hat{u}=\nabla\cdot h & \text{in $\QQ_L$,} \\
		\hat{u}=u_{\mathrm{bc}} & \text{on $\partial\QQ_L$,}
	\end{array}
	\right.
\end{equation*} with $u_{\text{bc}}$ defined by \eqref{EquationBCIntro}, then for any $R\in[\rss,L]$, \begin{equation*}
	\Bigl(\fint_{\bb_R} |\nabla(\hat{u}-u)|^2\Bigr)^{\frac{1}{2}} \lesssim \Big(\frac{\ell}{L}\Big)^d\Big(\frac{\rss}{L}\Big)^{\gamma_{\beta}(1-\varepsilon)}.
\end{equation*} 

\medskip

{\sc Step 3. Approximation of correctors. }Fix $L\gg 1$, set $\sqrt{M}= L^{1-\frac{\varepsilon}{2}}$ and define $(\phi_M^{(1)},\sigma_M^{(1)}, \phi_M^{(2)})$ as in \eqref{eq:1stOrderMassCor}, \eqref{eq:1stOrderMassFlux} and \eqref{eq:2ndOrderMassCor} respectively. Our estimates \eqref{eq:2ndCorApproxBd} and \eqref{eq:1stCorApprox} guarantee that the following approximation bound holds with high probability,
	\begin{equation}\label{eq:algevent1}
	\Bigl(\fint_{\bb_R} \bigl|\bigl(\sqrt{M}\nabla(\phi_M^{(1)}-\phi^{(1)}),\phi_M^{(1)}-\phi^{(1)},\nabla(\phi_{M}^{(2)}-\phi^{(2)})\bigr)\bigr|^2\Bigr)^\frac{1}{2}  \le  \sqrt{M}\Big(\dfrac{\rss}{L}\Big)^{\gamma_{\beta}(1-\varepsilon)}\ \mbox{ for }  R\in \bigl\{\ell, \dfrac{5}{4}L\bigr\}.
\end{equation}
Further, the growth condition \eqref{BoundPhi2} also ensures
\begin{equation}\label{eqn:algevent2}
\Big(\fint_{\bb_L} |\phi^{(2)}|^2\Big)^\frac{1}{2} \le  L^2\Big(\dfrac{\rss}{L}\Big)^{\gamma_{\beta}(1-\varepsilon)}.
\end{equation}
As for $\phi_M^{(2)}$, the same bound holds with high probability due to \eqref{eq:phi2ptws}. These massive correctors can be further replaced by their Dirichlet approximations $(\phi_{M,L}^{(1)},\sigma_{M,L}^{(1)}, \phi_{M,L}^{(2)})$, defined in \eqref{eqn:phiML}, \eqref{eqn:algsigma}, \eqref{eqn:2ndcorapprox}, with a deterministic error smaller than any algebraic power of $L$, thanks to \cite[Proposition 2.8]{lu2021optimal}. Finally the approximate homogenized coefficient
$a_{\text{hom},L}$, defined as in \eqref{eqn:algahL}, approximates well $a_{\text{hom}}$ in the sense that the following holds with high probability :
\begin{equation}\label{eqn:algevent3}
		|a_{\text{hom},L}-a_{\text{hom}}| \le  \Big(\dfrac{\rss}{L}\Big)^{\gamma_{\beta}(1-\varepsilon)},
\end{equation}
which only uses the fact that $\nabla \phi_{M,L}^{(1)}$ is a good approximation of $\nabla \phi^{(1)}$, as well as the fluctuation estimate of $a_{\text{hom}}$ developed in \cite[Corollary 2]{clozeau2021optimal}. The proof follows the same lines of \cite[Lemma 2.6]{lu2021optimal}.

\medskip

{\sc Step 4. Conclusion. }Step $2$ verifies that all conditions of \cite[Proposition 2.3]{lu2021optimal} are satisfied with high probability, and repeating the same arguments as \cite[Proposition 2.3]{lu2021optimal} is sufficient to show that \eqref{eq:AlgApprox} is typically true. The failure probability of the events \eqref{eq:phi2ptws}, \eqref{eq:algevent1}, \eqref{eqn:algevent2}, \eqref{eqn:algevent3} can be estimated using Chebyshev's inequality and the stochastic moments of the random constants in the estimates \eqref{eq:phi2ptws}, \eqref{eq:2ndCorApproxBd} and \eqref{eq:1stCorApprox}. Finally the moment bound of $\rss$ can be estimated using \eqref{BoundPhi2}, using the strategy developed in the proof of \cite[Theorem 1 (ii)]{fischer2017sublinear}, which captures the typical growth of $(\phi^{(2)},\sigma^{(2)})$ and also characterizes its failure probability.

\subsection{Proof of Proposition \ref{OptimalScaling}: Optimality of the algorithm for weakly correlated and small contrasted media}\label{OptimalitySection}
We split the proof into three steps. % In the first step, we appeal to the small contrast approximation and we estimate the error in terms of the contrast. At this stage, we reduce the proof to estimating the variance for the approximate solution that we express, 
In the first step, we express the covariance and the conditional covariance in terms of the Green function of the discrete Laplace operator. In the second step, we study the behaviour of the conditional covariance according to the correlation range $L^{\frac{d}{\beta}}$. In the third step, we study a double sum that we identify as the leading order term. We conclude in the fourth step by combining all estimates together. In the following, we use $\nabla$ to denote both the finite difference gradient and the continuum gradient operator.

\medskip
{\sc Step 1. Representation formula. }We now express the solution $\bar{u}$ of \eqref{SmallContrastEquation} in terms of the discrete Green function $G_D$ of $-\Delta$ on $\mathbb{Z}^d$, that is\footnote{For notational convenience we write $G_D(n)$ for $G_D(n,0)$}
$$\nabla\bar{u}(0)=\sum_{n\in\mathbb{Z}^d}\mathcal{G}_n\nabla^2 G_D(n)\nabla v(n),$$
so that 
\begin{equation}\label{ExpressVariance}
\begin{aligned}
&\mathbb{E}\big[\vert\nabla \bar{u}(0)-\mathbb{E}[\nabla\bar{u}(0)\vert \mathcal{F}_L]\vert^2\big]\\
&=\sum_{n\in\mathbb{Z}^d}\sum_{n'\in\mathbb{Z}^d}\nabla^2 G_D(n)\nabla v(n)\cdot \nabla^2 G_D(n')\nabla v(n')\mathbb{E}\big[(\mathcal{G}_n-\mathbb{E}[\mathcal{G}_n\vert\mathcal{F}_L])(\mathcal{G}_{n'}-\mathbb{E}[\mathcal{G}_{n'}\vert\mathcal{F}_L])\big].
\end{aligned}
\end{equation}
Now note that by definition of $\mathcal{F}_L$ we have $\mathbb{E}[\mathcal{G}_n\vert\mathcal{F}_L]=\mathcal{G}_n$ as long as $\vert n\vert_\infty\leq L$. In addition, by orthogonality, it holds
$$\mathbb{E}\big[\mathcal{G}_n\mathbb{E}[\mathcal{G}_{n'}\vert\mathcal{F}_L]\big]=\mathbb{E}\big[\mathcal{G}_{n'}\mathbb{E}[\mathcal{G}_{n}\vert\mathcal{F}_L]\big]=\mathbb{E}\big[\mathbb{E}[\mathcal{G}_{n}\vert\mathcal{F}_L]\mathbb{E}[\mathcal{G}_{n'}\vert\mathcal{F}_L]\big].$$
Thus, \eqref{ExpressVariance} turns into
\begin{equation}\label{ExpressionVarianceAlmostFinal}
\begin{aligned}
&\mathbb{E}\big[\vert\nabla \bar{u}(0)-\mathbb{E}[\nabla\bar{u}(0)\vert \mathcal{F}_L]\vert^2\big]\\
&=\sum_{\vert n\vert_\infty\geq L}\sum_{\vert n'\vert_\infty\geq L}\nabla^2 G_D(n)\nabla v(n)\cdot \nabla^2 G_D(n')\nabla v(n')(c(n-n')-\hat{c}(n,n')),
\end{aligned}
\end{equation}
where 
\begin{equation}\label{Formulachat}
\hat{c}(n,n'):=\mathbb{E}\big[\mathcal{G}_n\mathbb{E}[\mathcal{G}_{n'}\vert\mathcal{F}_L]\big].
\end{equation}
Finally, we further expand $\nabla^2 G_D(n)$ and $\nabla v(n)$ in the regime $\vert n\vert_\infty\uparrow\infty$. Denoting by $G$ the Green function of $-\Delta$ on $\mathbb{R}^d$, we appeal to the asymptotic of the discrete Green function in \cite[Theorem $4.5$]{mangad1967asymptotic} which states that there exists $C_d>0$ such that for any $n,n'\in\mathbb{Z}^d$
\begin{equation}\label{GreenExp}
\nabla^k G_D(n,n')=C_d\nabla^k G(n-n')+O(\vert n-n'\vert^{-d-k}_{\infty})\quad\text{for any $k\geq 0$}.
\end{equation}
Expressing $\nabla v$ in terms of the Green function $G_D$, using that $h$ is compactly supported in $\bb_\ell$ with $L\gg \ell$, the expansion \eqref{GreenExp} and recalling that $\sum_{\vert n\vert_\infty\leq \ell}h(n)=\ell^d\,\nu$, we find for all $n\in\mathbb{Z}^d$
\begin{equation}\label{Expansionv}
\nabla v(n)=\sum_{\vert n'\vert_\infty\leq \ell}\nabla^2 G_D(n,n')h(n')=C_d\ell^d\nabla^2G(n)\nu+O(\vert n\vert^{-d-1}_\infty).
\end{equation}
Combining \eqref{ExpressionVarianceAlmostFinal}, \eqref{GreenExp} and \eqref{Expansionv} together, we arrive at 
\begin{equation}\label{ExpressionVarianceFinal}
\begin{aligned}
&\mathbb{E}\big[\vert\nabla \bar{u}(0)-\mathbb{E}[\nabla\bar{u}(0)\vert \mathcal{F}_L]\vert^2\big]\\
=&C^2_d\,\ell^{2d}\underbrace{\sum_{\vert n\vert_\infty\geq L}\sum_{\vert n'\vert_\infty\geq L}\nabla^2 G(n)\nabla^2 G(n)\nu\cdot \nabla^2 G(n')\nabla^2 G(n')\nu\,(c(n-n')-\hat{c}(n,n'))}_{:=\mathcal{S}_L}\\
&+\sum_{\vert n\vert_\infty\geq L}\sum_{\vert n'\vert_\infty\geq L}O(\vert n\vert^{-2(d+1)}_\infty)O(\vert n'\vert^{-2(d+1)}_\infty)\,(c(n-n')-\hat{c}(n,n')).
\end{aligned}
\end{equation}
In the following we focus on $\mathcal{S}_L$, the second term can be treated the same way and can be identified as a higher order term.

\medskip

{\sc Step 2. Gaussian conditioning. }We show that there exists $\beta_d>d$ such that for all $\beta\geq \beta_d$ the two following assertions hold : 

\medskip

\begin{itemize}
\item[(i)]For any $R>0$ and $n'\in\mathbb{Z}^d$,
\begin{equation}\label{Boundchat1}
0\le \hat{c}(n,n')\lesssim_\beta R^{-\beta}L^{-d}\quad\text{provided $\vert n\vert_\infty\geq L+R L^{\frac{d}{\beta}}$},
\end{equation}
which also holds symmetrically by exchanging the role of $n$ and $n'$.

\medskip

\item[(ii)]For any $R>0$
\begin{equation}\label{Boundchat2}
0\le \hat{c}(n,n') \lesssim_\beta R^{-\beta}L^{-d}\quad\text{provided $\vert n-n'\vert_\infty\geq R L^{\frac{d}{\beta}}$.}
\end{equation}
\end{itemize}
To begin with, by Gaussian conditioning, for any $n\in\mathbb{Z}^d$ there exists $\Gamma^L_n:=\{\gamma^L_{k,n}\}_{\vert k\vert_\infty\leq L}\subset \mathbb{R}$ such that 
$$\mathbb{E}\big[\mathcal{G}_n\vert \mathcal{F}_L\big]=\sum_{\vert k\vert_\infty\leq L}\gamma^L_{k,n}\mathcal{G}_k.$$
In particular, $\hat{c}$ defined in \eqref{Formulachat} can be expressed as
\begin{equation}\label{ExpressCondiCovariance}
\hat{c}(n,n')=\sum_{\vert k\vert_\infty\leq L}\gamma^L_{k,n}\,c(n'-k).
\end{equation}
Furthermore, by orthogonality $\mathbb{E}\big[(\mathcal{G}_n-\mathbb{E}[\mathcal{G}_n\vert\mathcal{F}_L])\mathcal{G}_\ell\big]=0$ for any $\vert\ell\vert_\infty\leq L$, $\Gamma^L_n$ solves the linear system
\begin{equation}\label{SystemCorrelation}
\mathcal{C}\Gamma^L_n=\mathcal{C}_n,
\end{equation}
where $\mathcal{C}:=\{c(\ell-k)\}_{\vert \ell\vert_\infty,\vert k\vert_\infty\leq L}$ and $\mathcal{C}_n:=\{c(n-k)\}_{\vert k\vert_\infty\leq L}$. We now show that, for $\beta$ large enough, the matrix $\mathcal{C}$ is a perturbation of the identity in the following quantitative way : there exists $\beta_d>d$ such that for any $\beta\geq \beta_d$ it holds
\begin{equation}\label{PerturbationOfID}
\|\text{Id}-\mathcal{C}\|_{\infty}\ll 1.
\end{equation}
To do so, we express
\begin{align*}
\|\text{Id}-\mathcal{C}\|_{\infty}&=\sup_{\vert k\vert_{\infty}\leq L}\sum_{\vert\ell\vert_\infty\leq L}c(\ell-k)\mathds{1}_{\ell\neq k}\stackrel{\eqref{correlationExample}}{\leq} \sum_{n\in\mathbb{Z}^d\backslash \{0\}}(1+\vert n\vert)^{-\beta}.
\end{align*}
Since $\sum_{n\in\mathbb{Z}^d\backslash \{0\}}(1+\vert n\vert)^{-\beta}$ vanishes as $\beta\uparrow\infty$, we obtain \eqref{PerturbationOfID}. In particular, we deduce that the components of $\Gamma^L_n$ are all non-negative and it follows from \eqref{SystemCorrelation} that 
$$0\leq \gamma^L_{k,n}=c(n-k)-\sum_{\ell\neq k}\gamma^L_{\ell,n}c(\ell-k)\leq c(n-k)\quad\text{for any $\vert k\vert_\infty\leq L$}.$$
This together with \eqref{ExpressCondiCovariance} immediately gives \eqref{Boundchat1}, as
\begin{equation*}
	0\le \hat{c}(n,n') \le \sum_{|k|_\infty \le L} c(n-k)c(n'-k) \le R^{-\beta}L^{-d}  \sum_{|k|_\infty \le L} c(n'-k) \lesssim_\beta R^{-\beta}L^{-d}.
\end{equation*}
For \eqref{Boundchat2}, we split the sum in \eqref{ExpressCondiCovariance} into two parts and make use of the regime $\vert n-n'\vert_\infty\geq RL^{\frac{d}{\beta}}$ to get
\begin{align*}
\hat{c}(n,n')= \sum_{\vert k-n\vert_\infty\leq \frac{R}{2}L^{\frac{d}{\beta}}} c(n-k)c(n'-k)+\sum_{\vert k-n\vert_\infty\geq \frac{R}{2}L^{\frac{d}{\beta}}} c(n-k)c(n'-k)\lesssim_\beta R^{-\beta}L^{-d}.
\end{align*}
{\sc Step 3. Lower bound of the main term.} The goal of this step is to show that the main contribution of $\mathcal{S}_L$ (defined in \eqref{ExpressionVarianceFinal}) can be bounded from below by
\begin{equation}\label{AsymptoticDoubleSum}
	\begin{aligned}
		\Sigma_L:=\sum_{\vert n\vert_\infty\geq L}\sum_{\vert n'\vert_\infty\geq L}\nabla^2 G(n)\nabla^2 G(n)\nu\cdot \nabla^2 G(n')\nabla^2 G(n')\nu\,c(n-n')\gtrsim_{d,\beta} L^{-d}L^{-2d}.
	\end{aligned}
\end{equation}
Without loss of generality we set $\nu=e_1$, so that $\nu^\top n = n_1$ is the first coordinate of the vector $n$. Let $R\ge 1$ be a parameter to be determined later but independent of $L$. We divide the double sum into three regimes :
\begin{itemize}
	\item[(i)] Off-diagonal regime $D_1:=\big\{|n'-n|_\infty \ge \frac{L}{R}\big\}$;
	\item[(ii)] Diagonal regime perpendicular to $e_1$: $D_2:=\big\{|n'-n|_\infty\leq \frac{L}{R}$ with $\max\{|n_1|_\infty, |n_1'|_\infty \} \le \frac{2L}{R}\big\}$,
	\item[(iii)] Diagonal regime aligned with $e_1$: $D_3:=\big\{|n'-n|_\infty\leq \frac{L}{R}$ with $\max\{|n_1|_\infty, |n_1'|_\infty \} \ge \frac{2L}{R}\big\}$.
\end{itemize}
We carry out the estimates in the following steps, for which we use the notations $\Sigma^{D_i}_L$ for $i\in\{1,2,3\}$ which split accordingly the double sum in \eqref{AsymptoticDoubleSum}. 

\medskip

{\sc Substep 3.1. Off-diagonal regime. }We treat the contribution from $D_1$ as an error term, which can be estimated as follows:
\begin{align*}
	|\Sigma_L^{D_1}| & = \Big|\sum_{|n|_\infty \ge L} \sum_{|n'|_\infty\ge L,|n'-n| \ge \frac{L}{R}}  \nabla^2 G(n)\nabla^2 G(n)\nu\cdot \nabla^2 G(n')\nabla^2 G(n')\nu\,c(n-n')\Big| \\ & \lesssim (1+\frac{L}{R})^{-\beta} \sum_{|n'|_\infty\ge L}  \sum_{|n|_\infty \ge L} |n|^{-2d}|n'|^{-2d}  \lesssim R^\beta L^{-2d-\beta},
\end{align*} 
which is $\ll L^{-3d}$ for any choice of $R$ independent of $L$. 

\medskip

{\sc Substep 3.2. Diagonal regimes. } In this regime we first simplify the double sum using the precise expression of the Green function $G$ and then divide into two substeps, where we treat the contributions from $D_2$ and $D_3$ separately. Notice that when $|n'-n|_\infty \le \frac{L}{R}$, we have $|n|_\infty\sim |n'|_\infty$.

\smallskip

Since $G(n) = c_d|n|^{2-d}$ (here we assume $d>2$, but the estimates for $d=2$ are identical) is the Green function for the whole space Laplacian, a direct computation yields 
\begin{equation}\label{ActualG}
	c^{-1}_d\nabla^2 G(n) = \frac{d(d-2)}{|n|^{d+2}} nn^\top - \frac{d-2}{|n|^d}\mathrm{Id}.
\end{equation}
Substituting \eqref{ActualG} into \eqref{AsymptoticDoubleSum} and using the expression \eqref{correlationExample} of $c$, we obtain (up to a constant that only depends on $d$) :
\begin{equation}
	\label{App:Expd2d3}
	\begin{aligned}
		\Sigma_L^{D_2}+\Sigma_L^{D_3}&=\sum_{|n|_\infty \ge L}  \sum_{|n'|_\infty\ge L,|n'-n| \le \frac{L}{R}} \nabla^2 G(n)\nabla^2 G(n)\nu\cdot \nabla^2 G(n')\nabla^2 G(n')\nu\,c(n-n')\\ & \sim \sum_{|n|_\infty \ge L} \sum_{|n'|_\infty\ge L,|n'-n| \le \frac{L}{R}} (1+|n'-n|)^{-\beta} |n|^{-2d-2}|n'|^{-2d-2}\Big((d^2-2d)^2 n_1n_1'n^\top n' \\ & \qquad + (d^2-2d)\big(|n|^2(n_1')^2 + |n'|^2n_1^2\big) + |n|^2|n'|^2  \Big).
	\end{aligned}
\end{equation}
We notice immediately that the contributions from the last three terms above are positive. Furthermore, one can directly compute the contribution from the last term alone, which yields
\begin{equation}\label{App:Contripos}
	\begin{aligned}
		\sum_{|n|_\infty \ge L}  \sum_{|n'|_\infty\ge L,|n'-n| \le \frac{L}{R}} (1+|n'-n|)^{-\beta} |n|^{-4d}   \sim  \sum_{|n|_\infty \ge L}  |n|^{-4d}  \sim L^{-3d}.
	\end{aligned}
\end{equation}
Thus the remaining steps of the proof involve establishing lower bounds on the first term of \eqref{App:Expd2d3}, for which we divide into the following two substeps.

\medskip

{\sc Subsubstep 3.2.1. Diagonal regime perpendicular to $e_1$. }In this regime, one can perform direct estimation
\begin{equation}
	\label{App:lowerbdD2}
	\begin{aligned}
		\sum_{|n|_\infty \ge L} & \sum_{|n'|_\infty\ge L,|n'-n| \le \frac{L}{R}}\mathds{1}_{D_2}(n,n') (1+|n'-n|)^{-\beta} |n|^{-2d-2}|n'|^{-2d-2}n_1n_1'n^\top n' \\ & \ge -\sum_{|n|_\infty \ge L} \sum_{|n'|_\infty\ge L,|n'-n| \le \frac{L}{R}} \mathds{1}_{D_2}(n,n') (1+|n'-n|)^{-\beta} |n|^{-4d-2} |n_1||n_1'|\\ & \gtrsim -L^2R^{-2}\sum_{|n|_\infty \ge L}  \sum_{|n'|_\infty\ge L,|n'-n| \le \frac{L}{R}} (1+|n'-n|)^{-\beta} |n|^{-4d-2} \\ & \sim -L^2R^{-2}\sum_{|n|_\infty \ge L}  |n|^{-4d-2} \sim -R^{-2}L^{-3d}.
	\end{aligned}
\end{equation}
{\sc Subsubstep 3.2.2. Diagonal regime aligned with $e_1$. }If $n_1\ge \frac{2L}{R}$ then $n_1'\ge \frac{L}{R}$, since $|n-n'|\le \frac{L}{R}$; if $n_1\le -\frac{2L}{R}$ then $n_1'\le -\frac{L}{R}$. In any case we have $n_1n_1'\ge 0$. Meanwhile $n^\top n' = \frac{1}{2}(|n|^2+|n'|^2 - |n-n'|^2)>0$. Hence in $D_3$, 
\begin{equation}\label{App:D3T1pos}
	\sum_{|n|_\infty \ge L}  \sum_{|n'|_\infty\ge L,|n'-n| \le \frac{L}{R}}\mathds{1}_{D_3}(n,n') (1+|n'-n|)^{-\beta} |n|^{-2d-2}|n'|^{-2d-2}n_1n_1'n^\top n' \ge 0.
\end{equation}
Substituting \eqref{App:Contripos}, \eqref{App:lowerbdD2} and \eqref{App:D3T1pos} into \eqref{App:Expd2d3}, we get the bound from below
\begin{align*}
	\Sigma_L^{D_2}+ \Sigma_L^{D_3} \gtrsim (1-R^{-2})L^{-3d}.
\end{align*}
To conclude, we have shown
\begin{align*}
	\Sigma_L \gtrsim (1-R^{-2} - R^\beta L^{-(\beta-d)})L^{-3d},
\end{align*}
and we finish the proof by choosing $R\gg 1$ but independent of $L$.

\medskip

{\sc Step 4. Conclusion. }
Substituting \eqref{AsymptoticDoubleSum} into \eqref{ExpressionVarianceFinal} yields
\begin{equation}\label{ConclusionPartialEsti}
\mathcal{S}_L\geq C_d L^{-d}L^{-2d}-\sum_{\vert n\vert_\infty\geq L}\sum_{\vert n'\vert_\infty\geq L}\nabla^2 G(n)\nabla^2 G(n)\nu\cdot \nabla^2 G(n')\nabla^2 G(n')\nu\,\hat{c}(n,n').
\end{equation}
Similar to what we have done in Step 3, we choose a parameter $R>0$ independent of $L$ (that will be fixed later), then treat the double sum in \eqref{ConclusionPartialEsti} in different ways according to the three regimes: 
\begin{itemize}
\item[(i)] Far-field regime $D_1:=\big\{\max\{\vert n\vert_\infty,\vert n'\vert_\infty\}\geq L+RL^{\frac{d}{\beta}}\big\}$, 
\item[(ii)] Short-range dependence regime $D_2:=\big\{\max\{\vert n\vert_\infty,\vert n'\vert_\infty\}\leq L+RL^{\frac{d}{\beta}}$ with $\vert n-n'\vert_\infty\leq RL^{\frac{d}{\beta}}\big\}$,
\item[(iii)] Long-range dependence regime $D_3:=\big\{\max\{\vert n\vert_\infty,\vert n'\vert_\infty\}\leq L+RL^{\frac{d}{\beta}}$ with $\vert n-n'\vert_\infty\geq RL^{\frac{d}{\beta}}\big\}$.
\end{itemize}

\medskip

We carry out the estimates in the three following substeps, introducing the notations $\mathcal{S}^{D_i}_L$ for $i\in\{1,2,3\}$ which split accordingly the double sum in \eqref{ConclusionPartialEsti}.

\medskip

{\sc Substep 4.1. Far-field regime. }We combine \eqref{Boundchat1} with the bound on the Green function $\vert \nabla^2 G(n)\vert\lesssim \vert n\vert^{-d}$, to get
\begin{equation}
\vert\mathcal{S}^{D_1}_L\vert\lesssim_\beta R^{-\beta}L^{-d}\Big(\sum_{\vert n\vert_\infty\geq L}\vert\nabla^2 G(n)\vert^2\Big)^{2}\lesssim_\beta R^{-\beta}L^{-d}L^{-2d},
\label{AsympSum1}
\end{equation}
where $C_d,C_{d,\beta}>0$ depend on $d$ and $\beta$ respectively.

\medskip

{\sc Substep 4.2. Short-range dependence regime. }We use $\vert\hat{c}(n,n')\vert\le 1$ and $\vert\nabla^2 G(n)\vert\lesssim \vert n\vert^{-d}$ to obtain
\begin{equation}\label{EstiEnd2}
\begin{aligned}
\vert\mathcal{S}^{D_2}_L\vert &\lesssim \sum_{L\leq \vert n\vert_\infty\leq L+RL^{\frac{d}{\beta}}}\sum_{L\leq \vert n'\vert_\infty\leq L+RL^{\frac{d}{\beta}}}\vert\nabla^2 G(n)\vert^2\vert\nabla^2 G(n')\vert^2\mathds{1}_{D_2}(n,n')\\
&\lesssim R^d L^{-4d}\Big((L+RL^{\frac{d}{\beta}})^d-L^d\Big)L^{\frac{d^2}{\beta}}\lesssim R^{d+1}L^{-d-1+\frac{d^2}{\beta}+\frac{d}{\beta}}L^{-2d}.
\end{aligned}
\end{equation}
{\sc Substep 4.3. Long-range dependence regime. }We use \eqref{Boundchat2} in form of 
\begin{align}
\vert\mathcal{S}^{D_3}_L\vert\lesssim_\beta R^{-\beta}L^{-d}\Big(\sum_{\vert n\vert_\infty\geq L}\vert\nabla^2 G(n)\vert^2\Big)^2\lesssim_\beta R^{-\beta}L^{-d}L^{-2d}.
\label{EstiEnd1}
\end{align}
\medskip

The combination of \eqref{ConclusionPartialEsti}, \eqref{AsympSum1}, \eqref{EstiEnd2} and \eqref{EstiEnd1} provides
$$\mathcal{S}_L\geq C_dL^{-d}L^{-2d}-\mathcal{S}^{D_1}_L+\mathcal{S}^{D_2}_L+\mathcal{S}^{D_3}_L\geq (C_d-C_{d,\beta}R^{-\beta}-R^{d+1}L^{-1+\frac{d^2}{\beta}+\frac{d}{\beta}})L^{-d}L^{-2d}.$$
We conclude by taking $R$ large enough (but independent of $L$) and $\beta$ so large that $-1+\frac{d^2}{\beta}+\frac{d}{\beta}<0$.

%
%\subsection{Proof of Proposition \ref{LSregSoutionItself}: Large-scale $\cc^{0,1}$-estimates}

\subsection{Proof of Proposition \ref{RegResultParabolic}: Large-scale parabolic regularity}
%
%To begin with, the second item of \eqref{LSreg:PointwiseParabolic} is already contained in \cite[Lemma 4]{clozeau2021optimal}, we therefore focus on the proof of the first item. 
The proofs of all estimates largely follow the argument for \cite[Lemma 4]{clozeau2021optimal}. The main difference is on the assumptions we make : compared to \cite[Lemma 4]{clozeau2021optimal}, the \rhs of \eqref{EquationParaLargeScale} is not controlled deterministically but only in moments in probability, cf. \eqref{LSreg:Eq1}, \eqref{LSreg:Eq1p} and \eqref{LSreg:Eq1pq}. Therefore, a careful care on the estimates is done and we repeat the arguments while pointing out where it differs and where a change of argument is applied. Furthermore, our result is more general and treats as well non-divergence form \rhs that requires new arguments. The divergence form case and the non-divergence form case are treated separately in Case $1$ and Case $2$ below.

\medskip

We briefly recall the two main steps (that apply for each cases) that follows the strategy in \cite[Lemma 4]{clozeau2021optimal}. In the first step, we treat the particular case where $F_r$ (respectively $f_r$ in the non-divergence form case) is compactly supported in a ball $\bb_r$ for some $r\geq 1$. We show that  the estimates \eqref{LSreg:PointwiseParabolicnondiv}, \eqref{LSreg:PointwiseParabolic}, \eqref{LSreg:PointwiseParabolicalt} hold without the logarithmic corrections. In the second step, we treat the general case by decomposing $\mathbb{R}^d$ into dyadic annuli $(\mathcal{B}_k)_{k\in\mathbb{N}}$ defined by
$\mathcal{B}_k := \bb_{2^k r}\backslash\bb_{2^{k-1}r}$ for $k \geq 1$ and $\mathcal{B}_0:= \bb_r$, and writing $F_r=\sum_{k\geq 1}F_r \chi_k$ (respectively $f_r=\sum_{k\geq 1}f_r \chi_k$ for the non-divergence form case), where $(\chi_k)_{k\in\mathbb{N}}$ is a partition of unity according to the decomposition $(\mathcal{B}_k)_{k\in\mathbb{N}}$. We then apply the result of the compactly supported case for each $k \in\mathbb{N}$. For this step, which is unchanged, we simply refer to Step 2 of the proof of \cite[Lemma 4]{clozeau2021optimal}. %For notational convenience, we will drop in the notation the random constant $\C(r,x)$ that comes from the contribution of the \rhs $f_r$ and satisfies \eqref{MomentBoundGeneric}, and absorb that into $\lesssim$.

\medskip

{\sc \large Case 1. Divergence-form \rhs }We only prove \eqref{LSreg:PointwiseParabolic} with $f_r = \nabla\,\cdot\, aF_r$ under \eqref{LSreg:Eq1p}, the cases assuming $f_r = \nabla\cdot F_r$ or \eqref{LSreg:Eq1pq} are identical. %In y prove estimates for $v_r$; the bounds for $\nabla v_r$ are proved in \cite[Lemma 4]{clozeau2021optimal}. 

\medskip

{\sc Step 1. }We prove that under the stronger condition
\begin{equation}\label{LSreg:PointwisePara5}
	\text{supp } F_r\subset \bb_r\quad \text{and}\quad r^d \sup \vert F_r\vert+r^{d+1}\sup \vert \nabla F_r\vert\le \C(r,x),
\end{equation}
\eqref{LSreg:PointwiseParabolic} holds without logarithmic corrections, namely
\begin{equation}\label{LSreg:PointwisePara6}
	\bigl\vert \bigl(\nabla v_r(t,x), \tfrac{1}{\sqrt{|t|}}v_r(t,x)\bigr)\bigr\vert\le \frac{\C(r,x)}{(\vert x\vert+r)^d} \quad \text{for any $x\in\mathbb{R}^d$ and $\sqrt{-t}\geq 2r_\star(x)$}.
\end{equation}%\lw{I think we can improve here by removing the 1 in the max. See proofs below. This is useful since it makes all expressions simpler. Also, if we consider the solution of massive elliptic equation $\frac{1}{M} u -\nabla\cdot a \nabla u = \nabla \cdot af_r$, then by integrating against $e^{-t/M}$ in time, we get the consistent pointwise estimate \[|(\frac{1}{\sqrt{M}} u(x), \nabla u(x))| \lesssim (|x|+r)^{-d}.\]}
%
%The estimate \eqref{LSreg:PointwisePara6} follows from a combination of the following three estimates:
%
	%
	We claim that it suffices to prove 
	\begin{equation}\label{LSreg:PointwisePara6avg}
			\fint_{\bb_{r_\star(x)}(x)} 	\bigl\vert \bigl(\nabla v_r(t,\cdot), \tfrac{1}{\sqrt{|t|}}v_r(t,\cdot)\bigr)\bigr\vert^2 \le \frac{ \C(r,x)}{(|x|+r)^{2d}},
	\end{equation}
where we recall that $r_\star$ is defined in Proposition \ref{LSregSoutionItself}. The pointwise version \eqref{LSreg:PointwisePara6} is then obtained by the combination of Schauder's estimate recalled in Proposition \ref{SchauderTheoryStandard} and the moment bounds on $r_\star$ given in \eqref{Momentrstar}. We divide the proof into four substeps. In the first step, we prove a plain energy estimate on $(\nabla v_r,\tfrac{1}{\sqrt{\vert t\vert}}v_r)$. The three last sub-steps treat separately the regimes $\vert x\vert\geq 8(r_\star(x) \vee r)$, $r\geq (r_\star(x)\vee \tfrac{\vert x\vert}{8})$ and $r_\star(x)\geq (r\vee \frac{\vert x\vert}{8})$.
% In Substeps 1.1, 1.2 and 1.4 the proof for $\nabla v_r$ is done in \cite[Lemma 4]{clozeau2021optimal} so we only prove the estimates for $v_r$; the proof in Substep 1.3 is different from the previous work due to the involvement of random constants, and we therefore focus the proof on $\nabla v_r$ since the one for $v_r$ is simpler.
	
	\medskip
	
{\sc Substep 1.1. The plain energy estimate.} We claim that
	\begin{equation}\label{LSreg:PointwisePara7}
	\int_{\mathbb{R}^d}	\bigl\vert \bigl(\nabla v_r(t,\cdot), \tfrac{1}{\sqrt{|t|}}v_r(t,\cdot)\bigr)\bigr\vert^2 \le \C(r) r^{-d}\quad\text{for any $t\in (-\infty,0)$}.
\end{equation}We disintegrate $v_r$ in time using the auxiliary function $w_r$ satisfying in the weak sense
\begin{equation}\label{LSreg:PointwisePara11}
	\left\{
	\begin{array}{ll}
		\partial_\tau w_r-\nabla\cdot a\nabla w_r=0 & \text{in $(0,\infty)\times \mathbb{R}^d$,} \\
		w_r(0)=-\nabla\cdot aF_r,& 
	\end{array}
	\right.
\end{equation}
so that
\begin{equation}\label{LSreg:PointwisePara17}
	v_r(t,\cdot)=\int_{0}^{-t} \dd s\,w_r(s,\cdot)\quad \text{for any $t\in (-\infty,0)$,}
\end{equation}
see for instance \cite[Equation (4.4)]{clozeau2021optimal}.
Then applying the energy estimate \eqref{GOlm1} gives
\begin{equation} \label{LSreg:PointwisePara14}\begin{aligned}
		\int_{\R^d} \eta_{\sqrt{-t}} (x-\cdot) \big|\big(\nabla v_r(t,\cdot),\tfrac{1}{\sqrt{\vert t\vert}} v_r(t,\cdot)\big)\big|^2 & = \int_{\R^d} \eta_{\sqrt{-t}} (x-\cdot) \Big|\int_0^{-t} \dd s\,\big(\nabla w_r(s,\cdot),\tfrac{1}{\sqrt{\vert t\vert}}w_r(s, \cdot)\big) \Big|^2  \\ & \stackrel{\eqref{GOlm1}}{\lesssim} \int_{\R^d}\eta_{\sqrt{-t}} (x-\cdot) |F_r|^2  \\ & \stackrel{\eqref{LSreg:PointwisePara5}}{\leq} r^{-d}\fint_{\bb_r} \C(r,\cdot)\eta_{\sqrt{-t}}(x-\cdot).
\end{aligned}\end{equation}
We finally get \eqref{LSreg:PointwisePara7} by integrating over $x\in \R^d$ and noticing that $\fint_{\bb_r}\mathcal{C}(r,\cdot)=\mathcal{C}(r)$.

\medskip

{\sc Substep 1.2. Proof of \eqref{LSreg:PointwisePara6avg} when $\vert x\vert\geq 8(r_\star(x) \vee r)$.} We first prove by a duality argument that 
\begin{equation}\label{LSreg:PointwisePara3}
	\int_{\mathbb{R}^d\backslash \bb_{R}}	\bigl\vert \bigl(\nabla v_r(t,\cdot), \tfrac{1}{\sqrt{|t|}}v_r(t,\cdot)\bigr)\bigr\vert^2 \le %\max\Big\{\frac{r_\star(0)}{r}, 1\Big\}^{\frac{d}{2}}
\C(r)R^{-d}\quad\text{for all $R\geq 4(r_\star(0)\vee r)$ and $t\in (-\infty,0)$}.
\end{equation}
We only present the proof for $v_r$, the argument for $\nabla v_r$ are contained in \cite[(4.7)]{clozeau2021optimal} and in the arguments below. Let $\zeta\in \cc^{\infty}_c(\mathbb{R}^d)$ supported in $\mathbb{R}^d\backslash \bb_R$. Let $t\in (-\infty,0)$, $s\in [0,-t]$ and $k^s$ be the weak solution of the backward parabolic equation 
$$\left\{
\begin{array}{ll}
	\partial_{\tau} k^s+\nabla\cdot a\nabla k^s=0 & \text{in $(-\infty,s)\times\mathbb{R}^d$,} \\
	k^s(s,\cdot)= \zeta. & 
\end{array}
\right.
$$
Following the same computations leading to \cite[Estimate (4.13)]{clozeau2021optimal}, we have 
\begin{equation}\label{LSreg:PointwisePara4}
	\bigg\vert\int_{\R^d} \zeta\,v_r(t,\cdot)\bigg\vert\leq \C(r)\bigg(\fint_{\bb_r}\Big\vert\int_{0}^{-t}\nabla k^{s}(0,\cdot)\dd s\Big\vert^2\bigg)^{\frac{1}{2}}.
\end{equation}
The \rhs of \eqref{LSreg:PointwisePara4} is then estimated by introducing two auxiliary functions $\tilde{v}$ and $\bar{v}$, which are respectively the weak solutions of 
\[\left\{
\begin{array}{ll}
	\partial_{\tau} \tilde{v}+\nabla\cdot a\nabla \tilde{v}= \zeta & \text{in $(-\infty,0)\times\mathbb{R}^d$,} \\
	\tilde{v}(0,\cdot)=0. & 
\end{array}
\right.\quad\text{and}\quad \left\{
\begin{array}{ll}
	\partial_{\tau} \bar v+\nabla\cdot a\nabla \bar v=\mathds{1}_{(-\infty,0)}  \zeta & \text{in $\mathbb{R}^{d+1}$,} \\
	\bar v=0 &\text{ in $\mathbb{R}^+\times\mathbb{R}^d$.} 
\end{array}
\right.\]
Notice that $\bar v$ can be obtained by a zero extension of $\tilde{v}$ in $\R_+ \times \R^d$, since $\bar v(s,\cdot) = \tilde{v}(s,\cdot)$ whenever $s\le 0$. We perform the same estimates leading to \cite[Estimates (4.15)]{clozeau2021optimal}: Using the fact that $\zeta\equiv 0$ in $\bb_R$ and the relation $\tilde{v}(t,\cdot)=\int_{0}^{-t}\dd s\, k^{s}(0,\cdot)$, we obtain
\begin{equation}\label{LSreg:ConcluDuality3}
	\begin{aligned}
		\bigg(\fint_{\bb_r}\Big\vert\int_{0}^{-t}\nabla k^{s}(0,\cdot)\dd s\Big\vert^2\bigg)^{\frac{1}{2}} \lesssim \Big(\frac{r_\star(0)}{r} \vee 1\Big)^\frac{d}{2}\Big(\fint_t^{t+R^2} \dd s \fint_{\bb_{R}}  |\nabla \bar{v} (s,\cdot)|^2\Big)^\frac{1}{2}.
	\end{aligned}
\end{equation}
We then claim that 
\begin{equation}\label{ConcluDualityClaimvBar}
\fint_t^{t+R^2} \dd s \fint_{\bb_{R}}  |\nabla \bar{v} (s,\cdot)|^2\lesssim \vert t\vert R^{-d}\int_{\mathbb{R}^d}\vert \zeta\vert^2, 
\end{equation}
and we argue differently depending on the regime in $t$ and $R$ :
	\begin{itemize}
		\item For $\sqrt{-t}\le R$, we use \eqref{EnergyEstiAppendix:Eq1} with $g=h=0$ in form of
		\begin{align*}
			\fint_t^{t+R^2} \dd s \fint_{\bb_{R}}  |\nabla \bar{v} (s,\cdot)|^2&=\frac{|t|}{R^2}\fint_t^{0} \dd s \fint_{\bb_{R}}|\nabla \tilde{v} (s,\cdot)|^2  \stackrel{\eqref{EnergyEstiAppendix:Eq1}}{\lesssim} \frac{|t|^2}{R^2} \int_{\mathbb{R}^d} \eta_R |\zeta|^2 \lesssim |t|R^{-d} \int_{\mathbb{R}^d} |\zeta|^2. \stepcounter{equation} \tag{\theequation} \label{LSreg:rttleqR}
		\end{align*}
	\item For $\sqrt{-t}\ge R$, we claim that it suffices to prove 
	\begin{equation}\label{eq:extrpoGOlm1}
		\sup_{s \in [t,0]} \int_{\R^d} \eta_R |\bar{v}(s,\cdot)|^2 \lesssim |t|R^{2-d}\int_{\R^d} |\zeta|^2.
	\end{equation}
 Indeed, with \eqref{eq:extrpoGOlm1} we can close the estimate by using \eqref{EnergyEstiAppendix:Eq1} with initial condition $h=\bar{v}(t+R^2,\cdot)$ in form of
\begin{align*}
	\fint_t^{t+R^2} \dd s \fint_{\bb_{R}}  |\nabla \bar{v} (s,\cdot)|^2 \stackrel{\eqref{EnergyEstiAppendix:Eq1}}{\lesssim} R^{-2} \int_{\R^d} \eta_R |\bar{v}(t+R^2,\cdot)|^2 + R^2 \int_{\R^d} \eta_R |\zeta|^2 \stackrel{\eqref{eq:extrpoGOlm1}}{\lesssim} |t|R^{-d}\int_{\R^d} |\zeta|^2.
\end{align*}
We prove \eqref{eq:extrpoGOlm1} in a more general from, namely there exists a constant $C$ depending on $d$ and $\lambda$ such that for any $k\geq 2$
\begin{equation}\label{eq:extrpoGOlm1Bis}
\sup_{s \in [-kR^2,0]} \int_{\R^d} \eta_R |\bar{v}(s,\cdot)|^2 \leq C k R^{4-d}\int_{\R^d} |\zeta|^2,
\end{equation}
we then obtain \eqref{eq:extrpoGOlm1} by choosing the integer $k$ such that $(k-1)R^2\le |t| <kR^2$. We show \eqref{eq:extrpoGOlm1Bis} by induction. The proof for $k=2$ is a direct application of \eqref{EnergyEstiAppendix:Eq1}, more precisely
\begin{align*}
	 	\sup_{s \in [-2R^2,0]} \int_{\R^d} \eta_R |\bar{v}(s,\cdot)|^2 \stackrel{\eqref{EnergyEstiAppendix:Eq1}}{\le} C R^4 \int_{\R^d} \eta_R |\zeta|^2 \le C R^{4-d}\int_{\R^d} |\zeta|^2.
\end{align*}
Now suppose that \eqref{eq:extrpoGOlm1Bis} holds for $k$, again appealing to \eqref{EnergyEstiAppendix:Eq1} with initial condition $h=\bar{v}(-kR^2,\cdot)$ on the time interval $[-(k+1)R^2,-kR^2]$, we obtain
\begin{align*}
	\sup_{s \in [-(k+1)R^2,-kR^2]} \int_{\R^d} \eta_R |\bar{v}(s,\cdot)|^2 &  \stackrel{\eqref{EnergyEstiAppendix:Eq1}}{\leq} C\int_{\mathbb{R}^d} \eta_R |\bar{v}(-kR^2,\cdot)|^2 + CR^2\int_{-(k+1)R^2}^{-kR^2} \int_{\mathbb{R}^d} \eta_R|\zeta|^2 \\ & \leq CkR^{4-d} \int_{\mathbb{R}^d} |\zeta|^2 + CR^{4-d} \int_{\mathbb{R}^d} |\zeta|^2 = C(k+1)R^{4-d} \int_{\mathbb{R}^d} |\zeta|^2,
\end{align*}
which concludes the induction.
	\end{itemize}

The combination of \eqref{LSreg:PointwisePara4}, \eqref{LSreg:ConcluDuality3}, \eqref{ConcluDualityClaimvBar} yields \eqref{LSreg:PointwisePara3} by arbitrariness of $\zeta$. 
\medskip

We now prove \eqref{LSreg:PointwisePara6avg}. Let $R:=\frac{1}{2}\vert x\vert\geq 4 (r_\star(x)\vee r)$. We note that this further implies $R\geq 2r_\star(0)$, since due to the $\frac{1}{8}$-Lipschitz property of $r_\star$, we have $r_\star (0)\leq r_\star(x) + \frac{|x|}{8}  \le \frac{R}{4}+\frac{R}{4} = \frac{R}{2}$.
%and $\frac{3}{2}r_\star(0)\geq r_\star(x)$, so that using the plain energy estimate \eqref{LSreg:PointwisePara7} we obtain
%
%$$\fintstar \dd y\, \vert v_r(t,y)\vert^2\lesssim r^{-d}_\star(x)\int \dd y\,\vert v_r(t,y)\vert^2\lesssim r^d_\star(0)\frac{t^2}{(\vert x\vert+r)^{2d}}.$$
%
We now proceed as for \cite[Estimate (4.19)]{clozeau2021optimal}: Let $\bar{v}_r$ be the weak solution of
$$\left\{
\begin{array}{ll}
	\partial_\tau \bar{v}_r+\nabla\cdot a\nabla \bar{v}_r=\nabla\cdot \mathds{1}_{(-\infty,0)}aF_r & \text{in $\mathbb{R}^{d+1}$},\\
	\bar{v}_r=0 & \text{in $\mathbb{R}^+\times \mathbb{R}^d$.}
\end{array}
\right.$$
By construction $\bar{v}_r(s,\cdot)=v_r(s,\cdot)$ as long as $s\leq 0$ and $\bar{v}_r(s,\cdot)\equiv 0$ for $s\geq 0$. Note that $\bb_R(x) \subseteq \R^d \backslash \bb_R$ so that, thanks to \eqref{LSreg:PointwisePara5}, $F_r\equiv 0$ in $\bb_R(x)$. Therefore, a combination of Lemma \ref{PointwiseTimeEsti} and the $\cc^{0,1}$-large scale estimate \eqref{LargeScaleC01Paranondiv}  yields
\begin{equation}\label{LSreg:PointwisePara12}
	\fintstar  \vert v_r(t,\cdot)\vert^2\lesssim \fint_{t}^{t+R^{2}}\dd s\fint_{\bb_R(x)}\vert\bar{v}_r(s,\cdot)\vert^2.
\end{equation}
We thus deduce from \eqref{LSreg:PointwisePara3} that
\begin{itemize}
	\item for $\sqrt{-t}\geq R$,
	\begin{align*}
		\fint_{t}^{t+R^{2}}\dd s\fint_{\bb_R(x)}\vert\bar{v}_r(s,\cdot)\vert^2=\fint_{t}^{t+R^{2}}\dd s\fint_{\bb_R(x)}\vert v_r(s,\cdot)\vert^2\stackrel{\eqref{LSreg:PointwisePara3}}{\le} \C(r)|t|R^{-2d};
	\end{align*}
	\item for $R\geq \sqrt{-t}$,
	$$\fint_{t}^{t+R^{2}}\dd s\fint_{\bb_R(x)}\vert\bar{v}_r(s,\cdot)\vert^2=\Big(\frac{\sqrt{-t}}{R}\Big)^{2}\fint_{t}^{0}\dd s\fint_{\bb_R(x)}\vert v_r(s,\cdot)\vert^2\stackrel{\eqref{LSreg:PointwisePara3}}{\le} \C(r) |t| R^{-2d}.$$
\end{itemize}
The combination of \eqref{LSreg:PointwisePara12} with the two previous estimates concludes the argument for \eqref{LSreg:PointwisePara6avg}.

\medskip

{\sc Substep 1.3. Proof of \eqref{LSreg:PointwisePara6avg} when $r\ge (r_\star(x) \vee\tfrac{|x|}{8})$. } 
In this regime, it is sufficient to show that 
\begin{equation}\label{SmallScalePointwiseEsti}
\fint_{\bb_{r_\star(x)}(x)}	\bigl\vert \bigl(\nabla v_r(t,\cdot), \tfrac{1}{\sqrt{|t|}}v_r(t,\cdot)\bigr)\bigr\vert^2 \leq \mathcal{C}(r,x)r^{-2d}.
\end{equation}
We treat only the regime $\sqrt{-t}\geq r$, the proof for $\sqrt{-t}\leq r$ runs through the same argument by dividing the domain dyadically up to the scale $\sqrt{-t}$. Furthermore, we only expose the argument for $\nabla v_r$ as it shows the requirement of the regularity of $F_r$ in \eqref{LSreg:PointwisePara5}. The proof for $v_r$ is treated the same way and is simpler as it only requires the first item of \eqref{LSreg:PointwisePara5} and does not need to introduce the average \eqref{E:defxirk}. Note that our proof is different from the one from \cite[(4.9)]{clozeau2021optimal} since the \rhs of \eqref{EquationParaLargeScale} is only controlled in probabilistic norms (cf.  \eqref{LSreg:PointwisePara5}) instead of in $\LL^\infty$ sense. This prevents us to directly apply the large-scale regularity result \eqref{LargeScaleC01Para} as it would involve quantities like $\sup_{r_\star(x) \le \rho\le r} \fint_{\bb_{\rho}(x)} \C(r,\cdot)$ which has suboptimal moment bounds (in terms on the scaling in $r$)\footnote{For the purpose of this work, we may avoid the issue by replacing the pointwise assumptions in \eqref{LSreg:Eq1} and \eqref{LSreg:Eq1p} with averaged $\LL^2$ norm in $\bb_{r_\star(x)}(x)$, which allows us to obtain bounds with space-independent random constants, but this would make the presentation of the whole paper even more complicated and thus we avoid this treatment.}. Instead, we argue by dyadic decomposition, that involves dyadic sums that behave better under probabilistic norms.

\medskip

We fix $t\leq 0$ and $x\in\mathbb{R}^d$ and we introduce for any $k\ge 0$ a family of dyadic cylinders 
\[\mathcal{A}_{r,k}:=[t,t+2^{-2k}r^2]\times \bb_{2^{-k}r}(x),%\big\backslash [t,t+2^{-2k}r^2]\times \bb_{2^{-k}r}(x),
\]
and define (recalling that $F_r$ is time independent)
\begin{equation}\label{E:defxirk}
	\xi_{r,k} := \fint_{\mathcal{A}_{r,k}} F_r = \fint_{\bb_{2^{-k}r}(x)} F_r.
\end{equation}
We then decompose $\nabla v_r-\xi_{r,0}=\sum_{k\geq 0} \nabla v_{r,k}$ dyadically (that holds by uniqueness in \eqref{EquationParaLargeScale}), where for all $k\geq 1$
\begin{equation}\label{LSreg:PointwisePara13}
\left\{
    \begin{array}{ll}
        \partial_\tau v_{r,k}+\nabla\cdot a\nabla v_{r,k}=\nabla\cdot  aF_{r,k} &\text{in $(-\infty,0)\times \mathbb{R}^d$,} \\
        v_{r,k}(0,\cdot)\equiv 0, &
    \end{array}
\right.
\end{equation}
with
\begin{equation}\label{E:defFrk}
	F_{r,k} = \mathds{1}_{\mathcal{A}_{r,k-1}}(F_r-\xi_{r,k-1})-\mathds{1}_{\mathcal{A}_{r,k}}(F_r-\xi_{r,k}).
\end{equation}
and
\begin{equation}\label{LSreg:PointwisePara20}
	\left\{
	\begin{array}{ll}
		\partial_\tau v_{r,0}+\nabla\cdot a\nabla v_{r,0}=\nabla\cdot \mathds{1}_{[t+r^2,0]\times \bb_r(x)} a(F_{r}-\xi_{r,0}) &\text{in $(-\infty,0)\times \mathbb{R}^d$,} \\
		v_{r,0}(0,\cdot) = -x\cdot \xi_{r,0}. &
	\end{array}
	\right.
\end{equation}
%The motivation is that $w_r=\sum_{k\geq 0} v_{r,k}$ satisfies
We first claim that for any $k\ge 1$,
\begin{equation}\label{LSreg:PointwisePara16}
\fint_{t}^{t+\rho^2}\fint_{\bb_{\rho}(x)}\vert  \nabla v_{r,k}\vert^2\leq\mathcal{C}(k,r,x) r^{-2d}\left\{
    \begin{array}{ll}
           \Big(\frac{2^{-k}r}{\rho}\Big)^{d+2} & \text{if $\rho\geq 2^{-k}r$,}\\
      2^{-2k} & \text{if $r_\star(x)\leq \rho\leq 2^{-k}r$.}
    \end{array}
\right.
\quad
\end{equation}
In the regime $\rho\ge 2^{-k}r$ we obtain \eqref{LSreg:PointwisePara16} by applying the whole space energy estimate \eqref{FullSpaceEnergyEstimate} (combined with a time reflection $t\mapsto -t$) in form of
\begin{align*}
\fint_{t}^{t+\rho^2}\fint_{\bb_{\rho}(x)}\vert  \nabla v_{r,k}\vert^2\lesssim \rho^{-2-d}\int_{t}^0\int_{\mathbb{R}^d}\vert\nabla v_{r,k}\vert^2\stackrel{\eqref{FullSpaceEnergyEstimate}}{\lesssim}&\rho^{-2-d}\int_{t}^0\int_{\mathbb{R}^d} \vert F_{r,k}\vert^2\\
= & 2^{-2k}r^2\rho^{-2-d}\int_{\bb_{2^{-k}r}(x)}\vert F_{r,k}\vert^2\\
\stackrel{\eqref{LSreg:PointwisePara5}}{\leq}&\C(k,r,x)r^{-2d}\left(\frac{2^{-k}r}{\rho}\right)^{d+2} .
\end{align*}
In the regime $r_\star(x) \le \rho\leq 2^{-k}r$, we note that $\partial_{\tau}v_{r,k} + \nabla \cdot a(\nabla v_{r,k} + \xi_{r,k-1} - \xi_{r,k})=0$ in $\mathcal{A}_{r,k}$, hence applying the large-scale regularity estimate \eqref{LargeScaleC01Para} we obtain
\begin{equation}\label{LSreg:PointwisePara15}
\fint_{t}^{t+\rho^2}\fint_{\bb_{\rho}(x)}\vert \nabla v_{r,k}\vert^2\lesssim \fint_{\mathcal{A}_{r,k}}\vert \nabla v_{r,k}\vert^2 + |\xi_{r,k-1}-\xi_{r,k}|^2.
\end{equation}
For the first \rhs term of \eqref{LSreg:PointwisePara15}, we use the energy estimate \eqref{FullSpaceEnergyEstimate} to \eqref{LSreg:PointwisePara13} and the Poincar\'e inequality in $\bb_{2^{-(k-1)}r}(x)$ (recalling \eqref{E:defxirk}) in form of
\begin{align}
\fint_{\mathcal{A}_{r,k}}\vert \nabla v_{r,k}\vert^2  \lesssim (2^{-k}r)^{-2-d} \int_{t}^0 \int_{\mathbb{R}^d}\vert \nabla v_{r,k}\vert^2  
& \stackrel{\eqref{FullSpaceEnergyEstimate}}{\lesssim} (2^{-k}r)^{-2-d}\int_{t}^0 \int_{\mathbb{R}^d}\vert F_{r,k}\vert^2\nonumber\\
& \lesssim \fint_{\bb_{2^{-(k-1)}r}}\vert F_{r} - \xi_{r,k-1}\vert^2 +\fint_{\bb_{2^{-k}r}}\vert F_{r} - \xi_{r,k}\vert^2 \nonumber \\ & \lesssim (2^{-(k-1)}r)^2 \fint_{\bb_{2^{-(k-1)}r}}\vert \nabla F_{r} \vert^2 \stackrel{\eqref{LSreg:PointwisePara5}}{\le} \C(k,r,x)r^{-2d} 2^{-2k}. \label{E:EstigradvMVP}
\end{align}
The second \rhs term of \eqref{LSreg:PointwisePara15} can be treated directly using Poincar\'e's inequality and \eqref{LSreg:PointwisePara5} in form of
\begin{align*}
	|\xi_{r,k-1}-\xi_{r,k}|^2 = \Big| \fint_{\bb_{2^{-k}r}}(F_r - \xi_{r,k-1})\Big|^2 \lesssim \fint_{\bb_{2^{-(k-1)}r}}\vert F_r - \xi_{r,k-1}\vert^2 &\lesssim 2^{-2k}r^2 \fint_{\bb_{2^{-(k-1)}r}} |\nabla F_r|^2\\
	& \stackrel{\eqref{LSreg:PointwisePara15}}{\le} \C(k,r,x)r^{-2d} 2^{-2k}.
\end{align*}
This establishes \eqref{LSreg:PointwisePara16}. We now prove \eqref{SmallScalePointwiseEsti} and we first treat the term for $k=0$. Since $v_{r,0}$ is $a$-caloric in $\mathcal{A}_{r,0}$, we deduce from the large-scale regularity estimate \eqref{LargeScaleC01Para} and \eqref{LSreg:PointwisePara16} applied with $\rho=r$
\begin{align}
\fint_{t}^{t+4r^2_\star(x)}\fint_{\bb_{2r_\star(x)}}\vert \nabla v_{r,0}\vert^2 & \stackrel{\eqref{LargeScaleC01Para}}{\lesssim} \fint_{t}^{t+r^2}\fint_{\bb_r(x)}\vert \nabla v_{r,0}\vert^2\nonumber\\
&\lesssim \fint_{t}^{t+r^2}\fint_{\bb_r(x)}\vert \nabla v_r\vert^2+\vert\xi_{r,0}\vert^2+\Big(\sum_{k\geq 1}\Big(\fint_{t}^{t+r^2}\fint_{\bb_r(x)}\vert \nabla v_{r,k}\vert^2\Big)^{\frac{1}{2}}\Big)^2.\label{vr0}
\end{align}
The first and second \rhs term of \eqref{vr0} are controlled using \eqref{LSreg:PointwisePara7} and \eqref{LSreg:PointwisePara5} respectively in form of
\begin{equation}\label{vr0Bis}
\fint_{t}^{t+r^2}\fint_{\bb_r(x)}\vert \nabla v_r\vert^2+\vert\xi_{r,0}\vert^2\leq \Big(\mathcal{C}(r)+\Big(\fint_{\bb_r(x)}\mathcal{C}(r,\cdot)\Big)^2\Big)r^{-2d}=\mathcal{C}(r)r^{-2d}.
\end{equation}
For the last \rhs term of \eqref{vr0}, we appeal to \eqref{LSreg:PointwisePara16} with $\rho=r$ in form of
\begin{equation}\label{vr0Ters}
\begin{aligned}
\sum_{k\geq 1}\Big(\fint_{t}^{t+r^2}\fint_{\bb_r(x)}\vert \nabla v_{r,k}\vert^2\Big)^{\frac{1}{2}}\leq r^{-d}\sum_{r\geq 2^{-k}r} \mathcal{C}(k,r,x)2^{-(\frac{d}{2}+1)k}=\mathcal{C}(r,x)r^{-d}.
\end{aligned}
\end{equation}
The combination of \eqref{vr0}, \eqref{vr0Bis} and \eqref{vr0Ters} yields 
\begin{equation}\label{vr0Quatro}
\fint_{t}^{t+4r^2_\star(x)}\fint_{\bb_{2r_\star(x)}}\vert \nabla v_{r,0}\vert^2 \leq \mathcal{C}(r,x)r^{-2d}.
\end{equation}
We finally conclude using Lemma \ref{PointwiseTimeEsti}, the triangle inequality, applying \eqref{LSreg:PointwisePara16} with $\rho=r_\star(x)$ the same way as before, and \eqref{vr0Quatro} in form of
\begin{align*}
\fintstar\vert \nabla v_r(t,\cdot)\vert^2\stackrel{\eqref{Blowgrad}}{\lesssim} &\fint_{t}^{t+4r^2_\star(x)}\fint_{\bb_{2r_\star(x)}(x)}\vert \nabla v_r\vert^2 + \fint_{\bb_{2r_\star(x)}(x)}\vert F_r\vert^2\\
\stackrel{\eqref{LSreg:PointwisePara5}}{\lesssim} & \bigg(\sum_{k\geq 0} \bigg(\fint_{t}^{t+4r^2_\star(x)}\fint_{\bb_{2r_\star(x)}(x)}\vert \nabla v_{r,k}\vert^2\bigg)^{\frac{1}{2}}\bigg)^2+\mathcal{C}(r,x) r^{-2d}
\stackrel{\eqref{LSreg:PointwisePara16}, \eqref{vr0Quatro}}{\leq}\mathcal{C}(r,x) r^{-2d},
\end{align*}
which gives \eqref{SmallScalePointwiseEsti}. %The proofs for $v_r$ or \rhs given by $\nabla\cdot F_r$ are simpler: one can consider $v_{r,k}$ which is the solution of \eqref{LSreg:PointwisePara13} but with $F_{r,k} = (\mathds{1}_{\mathcal{A}_{r,k-1}}-\mathds{1}_{\mathcal{A}_{r,k}})F_r$ instead of \eqref{E:defFrk} (and $\xi_{r,k}$ are no longer necessary), and directly set $v_{r,0} = v_r - \sum_{k\geq 1} v_{r,k}$. For the proof with \rhs given by $\nabla\cdot F_r$, when using \eqref{EnergyEstiAppendix:Eq1} in \eqref{E:EstigradvMVP}, we may treat \rhs as if in non-divergence form and directly apply the bounds for $\nabla F_r$. When estimating $v_r$ instead of $\nabla v_r$, the bounds for $\nabla F_r$ are not used. 

\medskip

{\sc Substep 1.4. Proof of \eqref{LSreg:PointwisePara6avg} when $r_\star(x) \ge (r\vee \tfrac{|x|}{8})$. }We only prove the estimate for $v_r$, the estimate for $\nabla v_r$ is treated the same wayIn this regime, since $r_\star$ is $\tfrac{1}{8}$-Lipschitz, we have
\begin{equation}\label{LSreg:Rstarrelations}
	r_\star(0) \lesssim r_\star(x) \quad\text{and}\quad r_\star(x) \gtrsim r+|x|.
\end{equation}
Hence, using the plain energy estimate \eqref{LSreg:PointwisePara7} yields
\begin{align*}
	\fint_{\bb_{r_\star(x)}(x)} |(\nabla v_r(t,\cdot),\tfrac{1}{\sqrt{\vert t\vert}}v_r(t,\cdot))|^2 & \lesssim r_\star(x)^{-d} \int_{\R^d}  |(\nabla v_r(t,\cdot),\tfrac{1}{\sqrt{\vert t\vert}}v_r(t,\cdot))|^2\\
	&\stackrel{\eqref{LSreg:PointwisePara7}}{\lesssim}  r_\star^{-d}(0) \C(r)r^{-d} \stackrel{\eqref{LSreg:Rstarrelations}}{\lesssim} \frac{\C(r)r_\star^d(0)}{r^d(|x|+r)^{2d}}.
\end{align*}
\medskip

{\sc Step 2. The general case. }
We recover the general case by decomposing $F_r$ according to a family of dyadic annuli $(\mathcal{B}_k)_{k\in\mathbb{N}}$ defined by
$\mathcal{B}_k := \bb_{2^k r}\backslash\bb_{2^{k-1}r}$ for $k \geq 1$ and $\mathcal{B}_0:= \bb_r$. Setting $f_r=\sum_{k\geq 1}f_r \chi_k$, where $(\chi_k)_{k\in\mathbb{N}}$ is a partition of unity according to the decomposition $(\mathcal{B}_k)_{k\in\mathbb{N}}$, we write by linearity and uniqueness $v_r=\sum_{k\geq 0}v_{r,k}$. We then follow line by line the argument of Substep 1.2 of the proof of \cite[Lemma 4]{clozeau2021optimal}.

\medskip

{\sc \large Case 2. Proof of \eqref{LSreg:PointwiseParabolicnondiv}, non-divergence form \rhs }We again prove that under the stronger condition
\begin{equation}\label{LSreg:PointwisePara5pr}
	\text{supp } f_r\subset \bb_r\quad \text{and}\quad r^d \sup \vert f_r\vert \le \C(r,x),
\end{equation}
our result \eqref{LSreg:PointwiseParabolic} holds without logarithmic correction, that is
\begin{equation}\label{LSreg:PointwisePara6pr}
	|\nabla v_r(t,x)| \lesssim \frac{\C(r,x)}{(|x|+r)^{d-1}}.
\end{equation}
The proof largely follows from the roadmap of Case 1, and we only emphasize the differences.
	%

	%
%	\textbf{The large-scale regularity estimate:} for any $t\in (-\infty,0)$ and $\vert x\vert\geq 8(r_\star(x)\vee  r)$,
%	%
%	\begin{equation}\label{LSreg:PointwisePara8pr}
%		%
%		\fint_{\bb_{r_\star(x)}(x)} \dd y\,\vert \nabla v_r(t,y)\vert^2\le  \Big(\frac{r_\star(0)}{r}\vee 1\Big)^d\frac{1}{(\vert x\vert+r)^{2(d-1)}}.
%		%
%	\end{equation}
%	%
%	\textbf{The large-scale $\cc^{0,1}$-estimates:} for all $x\in\mathbb{R}^d$, $\sqrt{-t}\geq 2r_\star(x)$ and $r\geq r_\star(x)$,
%	
%	%
%	\medskip
%	%
%	
%	\begin{itemize}
%		%
%		\item For $\sqrt{-t}\leq r$, 
%		%
%		\begin{equation}\label{LSreg:PointwisePara9pr}
%			%
%			\fintstar \dd y\,\vert \nabla v_r(t,y)\vert^2\lesssim \fint_{t}^0\dd s\fint_{\bb_{\sqrt{-t}}(x)}\dd y\,\vert \nabla v_r(s,y)\vert^2+r^{-2(d-1)}.
%			%
%		\end{equation}
%		%
%		\item For $\sqrt{-t}\geq r$,
%		%
%		\begin{equation}\label{LSreg:PointwisePara10pr}
%			%
%			\fintstar \dd y\,\vert \nabla v_r(t,y)\vert^2\lesssim \fint_{t}^{t+r^2}\dd s\fint_{\bb_r(x)}\dd y\,\vert \nabla v_r(s,y)\vert^2+r^{-2(d-1)}.
%			%
%		\end{equation}
%		%
%	\end{itemize}
\medskip

{\sc Substep 1.1. The plain energy estimate.} We claim that	\begin{equation}\label{LSreg:PointwisePara7pr}
	\int_{\R^d} \vert \nabla v_r(t,\cdot)\vert^2\le \C(r) r^{2-d}\quad\text{for any $t\in (-\infty,0)$}.
\end{equation}
Let $w \in \text{H}^1(\mathbb{R}^d)/ \mathbb{R}$ be the energy solution of the elliptic equation
\begin{equation}\label{eq:Deltaweqf}
	-\Delta w = f_r,
\end{equation} then by Green's function representation $\nabla w (x) = \int_{\R^d} \dd y \nabla G(x-y) f_r(y)$ and \eqref{LSreg:PointwisePara5pr}, we have the following pointwise bound
\begin{equation}\label{BoundAuxilaryPbW}
	|\nabla w(x)| \le  \frac{\C(r,x)}{(|x|+r)^{d-1}}.
\end{equation}
Hence, applying plain energy estimate on \eqref{EquationParaLargeScale} by treating $f_r=\nabla\cdot \nabla w$ as a divergence-form term, we obtain as for \eqref{LSreg:PointwisePara7}
\begin{equation}\label{eq:engyestw}
	\int_{\R^d}  |\nabla v_r(t,\cdot)|^2 \le\int_{\R^d}  |\nabla w(t,\cdot)|^2 \le  \C(r)r^{2-d}.
\end{equation}
{\sc Substep 1.2. Proof of \eqref{LSreg:PointwisePara6pr} when $|x|\ge 8(r_\star(x) \vee r)$. }The argument requires a slight change from Case 1. Let $R\ge 4(r_\star(0)\vee r)$ and $\zeta\in (C_c^\infty(\R^d))^d$ be supported in $\R^d \backslash B_R$. For any $t\in (-\infty,0)$, let $\bar{v}$ be the weak solution of \[\left\{
\begin{array}{ll}
	\partial_{\tau} \bar{v} +\nabla\cdot a\nabla \bar{v}=\nabla \cdot \zeta & \text{in $(-\infty,0)\times\mathbb{R}^d$,} \\
	\bar{v}(0,\cdot)= 0,& 
\end{array}
\right.\]
then an identical duality argument to \cite[(4.13)]{clozeau2021optimal} yields 
\begin{align*}
	\Big|\int_{\R^d}\zeta\cdot \nabla v_r(t,\cdot)\Big| = \Big|\int_{\R^d} f_r\,\bar{v}(t,\cdot) \Big|.
\end{align*}
In the case $\sqrt{-t}\le R$, the proof is identical to the previous case \eqref{LSreg:rttleqR}
\begin{align*}
	\Big|\int_{\R^d}  f_r\,\bar{v}(t,\cdot) \Big| &  \stackrel{\eqref{LSreg:PointwisePara5pr}}{\le} \C(r)\Big(\fint_{\bb_{2r}}|\bar{v}(t,\cdot)|^2\Big)^\frac{1}{2} \stackrel{\eqref{LSreg:ConcluDuality3},\eqref{LSreg:rttleqR}}{\le} \C(r)\Big(\frac{|t|}{R^2}\fint_t^0 \fint_{\bb_{R}} |\bar{v}(s,\cdot)|^2\Big)^\frac{1}{2} \\ & \stackrel{\eqref{EnergyEstiAppendix:Eq1}}{\le} \C(r)\sqrt{-t}\Big( \int_{\R^d} \eta_{R} |\zeta|^2\Big)^\frac{1}{2}  \le \C(r) R^{1-\frac{d}{2}}\Big(\int_{\R^d}|\zeta|^2\Big)^\frac{1}{2}.
\end{align*}
The case $\sqrt{-t}\ge R$ is a bit more involved, since directly appealing to energy estimate \eqref{EnergyEstiAppendix:Eq1} on $\bar{v}$ would provide time-dependent bounds. The idea is to use the function $w$ defined in \eqref{eq:Deltaweqf}, so that
\begin{equation*}
	\Big|\int_{\R^d} f_r\bar{v}(t,\cdot) \Big| =  \Big|\int_{\R^d} \nabla w\cdot \nabla \bar{v}(t,\cdot) \Big|.
\end{equation*} We next divide the domain into domains $|x| \le R$ and $|x|\ge R$, and use the whole-space plain energy estimate $\int_{\R^d}|\nabla \bar{v}(t,\cdot)|^2\lesssim \int_{\R^d }|\zeta|^2$. Without loss of generality we assume $R,r$ are dyadic scales. When $x\in \bb_R$, we can estimate
\begin{align*}
	\Big|\int_{\bb_R}   \nabla w\cdot \nabla \bar{v}(t,\cdot) \Big| & \le  \Big|\int_{\bb_r} \nabla w(y)\cdot \nabla \bar{v}(t,\cdot) \Big| + \sum_{2^n r\le R} \Big|\int_{\bb_{2^n r} \backslash \bb_{2^{n-1}r}} \nabla w\cdot \nabla \bar{v}(t,\cdot) \Big|\\  & \stackrel{\eqref{BoundAuxilaryPbW}}{\le}  \int_{\bb_r} r^{1-d}\C(r,\cdot)|\nabla \bar{v}(t,\cdot) |+ \sum_{2^n r\le R} (2^nr)^{1-d}\int_{\bb_{2^n r} \backslash \bb_{2^{n-1}r}} \C(r,\cdot) |\nabla \bar{v}(t,\cdot) | \\ & \le r \C(r) \Big( \fint_{\bb_r} |\nabla \bar{v}(t,\cdot)|^2  \Big)^\frac{1}{2} + \sum_{2^n r\le R}  2^n r \C(2^nr)\Big( \fint_{\bb_{2^n r} }  |\nabla \bar{v}(t,\cdot)|^2  \Big)^\frac{1}{2}  \\ & \stackrel{\eqref{Blowgrad},\eqref{LargeScaleC01Para}}{\lesssim } r\C(r)  \Big( \fint_t^{t+R^2}\dd s\fint_{\bb_R} |\nabla \bar{v}(s,\cdot)|^2  \Big)^\frac{1}{2}  + \sum_{2^n r\le R} 2^n r\C(2^nr) \Big( \fint_t^{t+R^2}\dd s \fint_{\bb_R}  |\nabla \bar{v}(t,\cdot)|^2  \Big)^\frac{1}{2} \\ & \le \mathcal{C}(R)R^{1-\frac{d}{2}}\Big( \fint_t^{t+R^2}\dd s\int_{\R^d}  |\nabla \bar{v}(t,\cdot)|^2  \Big)^\frac{1}{2} \le \mathcal{C}(R)R^{1-\frac{d}{2}}\Big(\int_{\R^d} |\zeta|^2\Big)^\frac{1}{2}.
\end{align*}
The estimate for $\mathbb{R}^d\backslash\bb_R$ is more straightforward:
\begin{align*}
	\Big|\int_{\mathbb{R}^d\backslash \bb_R}    \nabla w\cdot \nabla \bar{v}(t,\cdot) \Big| & \le \sum_{2^n r\ge R} \Big|\int_{\bb_{2^n r} \backslash \bb_{2^{n-1}r}} \nabla w\cdot \nabla \bar{v}(t,\cdot) \Big| \\ & \stackrel{\eqref{BoundAuxilaryPbW}}{\le}\sum_{2^n r\ge R}  2^n r \C(2^nr)\Big( \fint_{\bb_{2^n r} }  |\nabla \bar{v}(t,\cdot)|^2  \Big)^\frac{1}{2} \\ & \le \sum_{2^n r\ge R} (2^n r)^{1-\frac{d}{2}}\C(2^nr) \Big( \int_{\R^d }|\nabla \bar{v}(t,\cdot)|^2  \Big)^\frac{1}{2} \le \C(R)R^{1-\frac{d}{2}} \Big(\int_{\R^d} |\zeta|^2\Big)^\frac{1}{2}.
\end{align*}
To combine, we have by duality argument that for any $R\ge 4(r_*(x)\vee r)$, 
\begin{equation*}
	\Big(\int_{\mathbb{R}^d\backslash \bb_R} |\nabla v_r(t,\cdot)|^2\Big)^\frac{1}{2} \le  \C(R)R^{1-\frac{d}{2}}.
\end{equation*}
The rest of the whole proof, including Substeps 1.2, 1.3, 1.4 and Step 2, proceed identically to the divergence-form case, with the only exception being that \eqref{LargeScaleC01Para} is used instead of \eqref{LargeScaleC01Paranondiv}.

\subsection{Proof of Corollary \ref{RegResultElliptic}: Large-scale elliptic regularity}
The proof relies on the following observation: Define $w$ as the weak solution of the backward parabolic equation
\begin{equation*}
	\left\{
	\begin{array}{ll}
		\partial_{\tau}w+\nabla\cdot a\nabla w=-f_r & \text{ in $(-\infty,0)\times \mathbb{R}^d$,} \\
		w(0,\cdot)=0,& 
	\end{array}
	\right.
\end{equation*}
a direct computation shows that the solution $u_r$ of \eqref{LSreg:Eq3} is given by 
\begin{equation*}
	u_r = M\int_{-\infty}^0 \dd t\, e^{\frac{t}{M}} w(t,\cdot).
\end{equation*}
Hence, by time integration, the bounds of Proposition \ref{RegResultParabolic} translate into these in Corollary \ref{RegResultElliptic}.

\subsection{Proof of Lemma \ref{lem:funcderiv}: Functional derivative of the $2^{\text{nd}}$-order time-dependent flux}
%
%\nc{I think we miss an approximating argument here. Indeed, \eqref{EnergyPassToTheLimit:Eq1} is not directly true since we need to estimate the $\LL^2$-norm of the initial deta in \eqref{eq:deltahu2} which requires regularity of the test function $\delta a$ that we don't have. What I propose: we smooth a bit the initial data replacing $a$ by $a_\varepsilon:=a\star \chi_\varepsilon$. Then it is fine, we can show the Lemma but the fourth right-hand side term of \eqref{eqn:funcderiv} will change and an additional convolution with the kernel $\chi_\varepsilon$ will appear. Then we can proceed and show estimates uniform in this regularization parameter. Finally, thanks to \eqref{SmoothnessCoef}, we can show that the regularized $q^{(2)}_{r,\varepsilon}$ converges to $q^{(2)}_r$ in every $\LL^{p}(\Omega)$  and we conclude.} \lw{Is it possible to place the $\int_0^1 \dd t$ inside the absolute value?}
For notational convenience, we drop the dependences on the indices $ij$. Using the definition of the functional derivative \eqref{DefDerivative}, we proceed by taking finite differences and showing that it passes to the limit as $h\downarrow 0$. Let us fix $x\in\mathbb{R}^d$, $h\in (0,1)$, $\ell\geq 1$ and $\delta a\in \LL^{\infty}(\mathbb{R}^d)$ such that $\text{supp } \delta a\subset \bb_{\ell}(x)$. Introducing the following notation for any random field $F$ of the coefficient field $a$ : 
\[\delta^h F(a,\cdot):=\frac{F(a+h\delta a,\cdot)-F(a,\cdot)}{h},\]
we take finite differences on \eqref{2NDOrder:Eq8} which leads to, for any unit vector $\upsilon$,
\begin{equation}\label{eq:1stfd}
\begin{aligned}
\delta^h \big(q^{(2)}_{M,\varepsilon}(T,\cdot)+\sigma^{(1)}_{M}e_j\big)_r(0)\cdot \upsilon & = \int_{\R^d} g_r\,\upsilon\cdot \delta a\,\phi_M^{(1)} (a+h\del a,\cdot)+\int_{\R^d} \int_0^T \dd s\,g_r\, \upsilon\cdot \delta a \nabla u^{(2)}_{M,\varepsilon}(a+h\delta a,s,\cdot) \\ & \qquad +\int_{\R^d} \int_0^T \dd s\,g_r\,\upsilon\cdot a \nabla \delta^h u^{(2)}_{M,\varepsilon}(s,\cdot)+\int_{\R^d} g_r\,\upsilon\cdot a\, \del^h \phi_M^{(1)}e_j,
\end{aligned} \end{equation}
where, using \eqref{2NDOrder:Eq1} and \eqref{eq:1stOrderMassCor}, $\delta^h u^{(2)}_{M,\varepsilon}$ and $\delta^{h}\phi^{(1)}_M$ satisfy 
\begin{equation}\label{eq:deltahu2}
	\left\{ \begin{aligned}
		& \partial_\tau \delta^h u^{(2)}_{M,\varepsilon}-\nabla \cdot a \nabla \delta^h u^{(2)}_{M,\varepsilon}= \nabla \cdot \delta a \nabla u^{(2)}_{M,\varepsilon}(a+h\delta a,\cdot,\cdot) \\ & \delta^h u^{(2)}_{M,\varepsilon}(0) = \nabla \cdot \big((\delta a\star \chi_\varepsilon)\, \phi^{(1)}_M(a+h\delta a,\cdot) + a_\varepsilon\, \delta^h \phi^{(1)}_M-\delta^h \sigma^{(1)}_M)e_j\big),
	\end{aligned}   \right.
\end{equation}
and 
\begin{equation}\label{Eq:DeltahPhiM}
\tfrac{1}{M}\delta^h\phi^{(1)}_M-\nabla \cdot a \nabla \delta^h \phi^{(1)}_M = \nabla \cdot \delta a (e_i+\nabla \phi^{(1)}_M(a+h\delta a,\cdot)).
\end{equation}
We now make use of the localized energy estimates of Proposition \ref{lem:GOlm1} to show the continuity of the maps $h\mapsto \big(\phi^{(1)}_M(a+h\delta a,\cdot),\nabla \phi^{(1)}_M(a+h\delta a,\cdot)\big)$, $h\mapsto \big(\sigma^{(1)}_M(a+h\delta a,\cdot),\nabla \sigma^{(1)}_M(a+h\delta a,\cdot)\big)$ and $h\mapsto \nabla u^{(2)}_{M,\varepsilon}(a+h\delta a,\cdot)$ in $\LL^{2}_{\text{loc}}(\mathbb{R}^d)$. First, applying \eqref{LocalizedMassiveTerm} to \eqref{eq:1stOrderMassCor} and to its difference with \eqref{Eq:DeltahPhiM} yield for any $R>0$
\begin{equation}\label{EnergyPassToTheLimit:Massive}
\begin{aligned}
&\fint_{\bb_{R}(x)} \big\vert \big(\tfrac{1}{\sqrt{M}}(\phi^{(1)}_M(a+h\delta a,\cdot)-\phi^{(1)}_{M}(a,\cdot)),\nabla\phi^{(1)}_M(a+h\delta a,\cdot)-\nabla\phi^{(1)}_{M}(a,\cdot)\big)\big\vert^2\\
&\lesssim \Big(\frac{\sqrt{M}}{R}\vee 1\Big)^d h^2 \int_{\R^d}\eta_{\sqrt{M}}(\cdot-x)\vert\delta a\vert^2\vert e_i+\nabla\phi^{(1)}_M(a+h\delta a,\cdot)\vert^2\lesssim_{R,M} h^2.
\end{aligned}
\end{equation}
Likewise, using \eqref{eq:1stOrderMassFlux} and \eqref{EnergyPassToTheLimit:Massive},
\begin{equation}\label{EnergyPassToTheLimit:Massive2}
\begin{aligned}
&\fint_{\bb_{R}(x)} \big\vert \big(\tfrac{1}{\sqrt{M}}(\sigma^{(1)}_M(a+h\delta a,\cdot)-\sigma^{(1)}_{M}(a,\cdot)),\nabla\sigma^{(1)}_M(a+h\delta a,\cdot)-\nabla\sigma^{(1)}_{M}(a,\cdot)\big)\big\vert^2\lesssim_{R,M} h^2.
\end{aligned}
\end{equation}
Second, we apply \eqref{EnergyEstiAppendix:Eq1} to \eqref{eq:deltahu2} and to its difference with \eqref{2NDOrder:Eq1Regularized} that we combine with \eqref{EnergyPassToTheLimit:Massive} to obtain for any $R>0$
\begin{equation}\label{EnergyPassToTheLimit:Eq1}
\begin{aligned}
&\fint_{\bb_{R}(x)}\int_{0}^{T}\dd s\,\big\vert\nabla u^{(2)}_{M,\varepsilon}(a+h\delta a,s,\cdot)-\nabla u^{(2)}_{M,\varepsilon}(a,s,\cdot) \big\vert^2 \\
\lesssim\, &h^2 \Big(\frac{\sqrt{T}}{R}\vee 1\Big)^d\int_{\R^d} \eta_{\sqrt{T}}(\cdot-x)\vert \nabla \cdot \big((\delta a\star \chi_\varepsilon)\, \phi^{(1)}_M(a+h\delta a,\cdot) + a_\varepsilon\, \delta^h \phi^{(1)}_M-\delta^h \sigma^{(1)}_M)e_j\big)\vert^2\\
&+\int_{0}^T\int_{\R^d}\eta_{\sqrt{T}}(\cdot-x)\vert\delta a\vert^2\vert\nabla u^{(2)}_{M,\varepsilon}(a+h\delta a,\cdot)\vert^2\stackrel{\eqref{EnergyPassToTheLimit:Massive},\eqref{EnergyPassToTheLimit:Massive2}}{\lesssim_{R,T,M,\varepsilon}} h^2.
\end{aligned}
\end{equation}
The combination of \eqref{EnergyPassToTheLimit:Massive} and \eqref{EnergyPassToTheLimit:Eq1} allows us to pass to the limit as $h\downarrow 0$ in the first two terms of \eqref{eq:1stfd} and give the first and second terms of \eqref{eqn:funcderiv} respectively.

\medskip

It remains to treat the two last terms of \eqref{eq:1stfd}. We first rewrite the third term using the dual equation \eqref{eq:vT} satisfied by $v^T$: Testing $\delta^h u^{(2)}_{M,\varepsilon}$ in \eqref{eq:vT} and $v^T$ in \eqref{eq:deltahu2} together with an integration by parts gives 
\begin{equation}\label{eq:der4terms}
\begin{aligned}
\int_{\R^d} \int_0^T \dd s\, g_r \upsilon\cdot a \nabla \delta^h u^{(2)}_{M,\varepsilon}(s,\cdot)  \stackrel{\eqref{eq:vT}}{=}&-\int_{\R^d} \int_0^T\dd s\, \partial_\tau v^T \delta^h u^{(2)}_{M,\varepsilon}(s,\cdot)+\int_{\R^d} \int_0^T\nabla\delta^h u^{(2)}_{M,\varepsilon}(s,\cdot)\cdot a\nabla v^T \\ 
\stackrel{\eqref{eq:deltahu2}}{=}& -\int_{\R^d} \nabla v^T(0,\cdot)\cdot (\delta a\star \chi_\varepsilon) \phi^{(1)}_M(a+h\delta a,\cdot)+ \int_{\R^d} v^{T}(0,\cdot) \nabla\cdot (a_\varepsilon\, \delta^h \phi^{(1)}_M-\delta^h \sigma^{(1)}_M)e_j\\
&+ \int_{\R^d} \int_0^T\dd s\, v^T(s,\cdot) \nabla \cdot \delta a \nabla u^{(2)}_{M,\varepsilon}(s,\cdot)(a+h\delta a,s,\cdot). 
\end{aligned}
\end{equation}
Using \eqref{EnergyPassToTheLimit:Massive} and \eqref{EnergyPassToTheLimit:Eq1}, we can pass to the limit as $h\downarrow 0$ in the first and third \rhs terms of \eqref{eq:der4terms} yielding the third and fourth terms of \eqref{eqn:funcderiv} respectively. We now deal with the second \rhs term of \eqref{eq:der4terms}. Expending the divergence and using an integration by parts yield
\begin{equation}\label{MassiveTermLast:Eq1}
\begin{aligned}
\int_{\R^d} v^{T}(0,\cdot)\nabla\cdot \big(a_\varepsilon\,\delta^h\phi^{(1)}_M-\delta^{h}\sigma^{(1)}_M\big)e_j=-\int_{\R^d} \nabla v^{T}(0,\cdot)\cdot a_\varepsilon\,\delta^h\phi^{(1)}_M e_j-\int_{\R^d} v^{T}(0,\cdot)\nabla\cdot \delta^h\sigma^{(1)}_Me_j.
\end{aligned}
\end{equation}
The second term of \eqref{MassiveTermLast:Eq1} can be simplified using the equations on $\sigma^{(1)}_M$. Indeed, introducing the auxiliary vector field
\begin{equation}\label{AxuliaryVectorField}
(\tfrac{1}{M}-\Delta)\zeta_M=\tfrac{1}{M}\big(q^{(1)}_M-\mathbb{E}[q^{(1)}_M]-\nabla\phi^{(1)}_M\big),
\end{equation}
it holds
\begin{equation}\label{DivergenceSigmaMassive}
\nabla\cdot \sigma^{(1)}_M=q^{(1)}_M-\mathbb{E}[q^{(1)}_M]-\zeta_M,
\end{equation}
where $(\nabla\cdot \sigma^{(1)}_{M})_{ij}=\sum_k \partial_k \sigma^{(1)}_{ijk,M}$, i.e the divergence is taken with respect to the last coordinate. The identity \eqref{DivergenceSigmaMassive} can be shown by applying the divergence in \eqref{eq:1stOrderMassFlux} to the effect of 
$$\tfrac{1}{M}(\nabla\cdot \sigma^{(1)}_M)_{ij}-\Delta (\nabla\cdot \sigma^{(1)}_M)_{ij}=\partial_j\nabla\cdot q^{(1)}_{i,M}-\Delta q^{(1)}_{ij,M}.$$
Inserting \eqref{eq:1stOrderMassCor} together with the definition $q^{(1)}_{i,M}=a(\nabla\phi^{(1)}_{i,M}+e_i)$ yields
$$\tfrac{1}{M}(\nabla\cdot \sigma^{(1)}_M)_{i}-\Delta (\nabla\cdot \sigma^{(1)}_M)_{i}=\frac{1}{M}\nabla \phi^{(1)}_{i,M}-\Delta q^{(1)}_{ij,M}.$$
Finally, adding \eqref{AxuliaryVectorField} leads to 
$$(\tfrac{1}{M}-\Delta)((\nabla\cdot \sigma)_i+\zeta_{i,M})=\tfrac{1}{M}(q^{(1)}_{i,M}-\mathbb{E}[q^{(1)}_{i,M}])-\Delta q^{(1)}_{i,M},$$
which yields \eqref{DivergenceSigmaMassive} by invertibility of $\tfrac{1}{M}-\Delta$ on bounded fields. 

\medskip

Now, taking finite differences in \eqref{DivergenceSigmaMassive} and \eqref{AxuliaryVectorField} provides
\begin{equation}\label{FiniteDiffSigmaMassive}
\nabla\cdot \delta^h\sigma^{(1)}_M=\delta^h q^{(1)}_M-\delta^{h}\zeta_M,
\end{equation}
and
\begin{equation}\label{DiffZeta}
(\tfrac{1}{M}-\Delta)\delta^h \zeta_M=\tfrac{1}{M}\big(\delta^h q^{(1)}_M-\nabla\delta^h\phi^{(1)}_M\big).
\end{equation}
Applying \eqref{FiniteDiffSigmaMassive}, the second term of \eqref{MassiveTermLast:Eq1} turns into
$$\int_{\R^d} v^{T}(0,\cdot)\nabla\cdot \delta^h\sigma^{(1)}_Me_j=\int_{\R^d} v^{T}(0,\cdot)\delta^{h}q^{(1)}_M e_j-\int_{\R^d} v^{T}(0,\cdot)\delta^h\zeta_Me_j.$$
Testing $\delta^h \zeta_M$ in \eqref{eqn:phi1dual} and $w_1$ in \eqref{DiffZeta} together with $\delta^h q^{(1)}_M=\delta a\big(e+\nabla \phi^{(1)}_M(a+h\delta a,\cdot)\big)+a \nabla\delta^h\phi^{(1)}_M$ yields
\begin{equation} \label{eq:deltahqM}
\begin{aligned}
\int_{\R^d} v^{T}(0,\cdot)e_j\cdot \delta^{h}q^{(1)}_M-\int_{\R^d} v^{T}(0,\cdot)e_j\cdot\delta^h\zeta_M=&\int_{\R^d}  (v^{T}(0,\cdot)e_j-\tfrac{1}{M}w_1)\cdot\delta^h q^{(1)}_M+\int_{\R^d} \tfrac{1}{M} w_1\cdot \nabla \delta^h \phi^{(1)}_M\\
=&\int_{\R^d} (v^{T}(0,\cdot)e_j-\tfrac{1}{M}w_1)\cdot \delta a(e_i+\nabla\phi^{(1)}_M(a+h\delta a,\cdot))\\
&+\int_{\R^d} (v^{T}(0,\cdot)a e_j-\tfrac{1}{M}(a-\mathrm{Id})w_1)\cdot \nabla\delta^h\phi^{(1)}_M. 
\end{aligned}
\end{equation}
Applying \eqref{EnergyPassToTheLimit:Massive}, the first integral in \eqref{eq:deltahqM} gives the sixth term of \eqref{eqn:funcderiv} as $h\downarrow 0$. Finally, collecting all remaining terms, which include the last \rhs term of \eqref{eq:1stfd}, the first \rhs term of \eqref{MassiveTermLast:Eq1} and the second \rhs term of \eqref{eq:deltahqM}, we derive testing $\delta^h \phi^{(1)}_M$ in \eqref{eqn:2ndphi1dual} and $w_2$ in \eqref{Eq:DeltahPhiM}
\begin{align*}
	\int_{\R^d} g_r\,\upsilon\cdot a & \del^h  \phi_M^{(1)}e_j-\int_{\R^d} \nabla v^{T}(0,\cdot)\cdot a_\varepsilon \,\delta^h\phi^{(1)}_M e_j -\int_{\R^d} (v^{T}(0,\cdot)a e_j+\tfrac{1}{M}(a-\mathrm{Id})w_1)\cdot \nabla\delta^h\phi^{(1)}_M \\ & = \int_{\R^d} \del^h \phi_M^{(1)} \Big(g_re_j\cdot a \upsilon - \nabla v^T(0,\cdot)\cdot a_\varepsilon\,e_j\, +\nabla\cdot a v^T(0,\cdot)e_j-\tfrac{1}{M}\nabla\cdot \big((a-\mathrm{Id})w_1\big)\Big) \\ & \stackrel{\eqref{eqn:2ndphi1dual},\eqref{Eq:DeltahPhiM}}{=} -\int_{\R^d}\nabla w_2\cdot \delta a (e_i+\nabla \phi^{(1)}_M(a+h\delta a,\cdot)),
\end{align*}
which finally corresponds to the fifth term of \eqref{eqn:funcderiv} after passing to the limit as $h\downarrow 0$ using \eqref{EnergyPassToTheLimit:Massive}.

\subsection{Proof of Lemma \ref{lem:ptwisew1}: Pointwise bounds for the massive Laplacian operator}
	%We start the proof with recalling the bound on $v^T$:
%\begin{align*}
	%|v^T(0,x)| \le \mathcal{C}(r,x) \frac{\log(1+\frac{|x|}{r})}{(|x|+r)^{d-1}}\max\Big\{\frac{\sqrt{T}}{|x|+r},1 \Big\}.
%\end{align*}
The proof is based on the following explicit formula for $w_1$ using the Green function $G_M$ of the elliptic operator $\tfrac{1}{M}-\Delta$, that is
\begin{equation}\label{eqn:w1green}
w_1= \int_{\R^d} \dd y\,G_M(\cdot-y) v^T(0,y).
\end{equation}
We recall that $G_M$ satisfies the following bounds: there exists a constant $C>0$ such that for any $x\in\mathbb{R}^d$
\begin{equation}\label{eqn:greenGM}
|G_M(x)| \lesssim |x|^{2-d} \exp\Big(-\frac{1}{C}\frac{|x|}{\sqrt{M}}\Big)\quad\mbox{and}\quad|\nabla G_M(x)| \lesssim |x|^{1-d} \exp\Big(-\frac{1}{C}\frac{|x|}{\sqrt{M}}\Big).
\end{equation}
As a consequence of the large-scale regularity estimate \eqref{LSreg:PointwiseParabolic}, the localized energy estimate \eqref{EnergyEstiAppendix:Eq1} and Schauder's estimates recalled in Proposition \ref{SchauderTheoryStandard}, we have
\begin{equation}\label{BoundvTForMassive}
|v^T(0,y)| \le \mathcal{C}(r,y) \frac{\sqrt{T}\log(2+\tfrac{|y|}{r})}{(|y|+r)^d}\quad\text{for any $y\in\mathbb{R}^d$,}
\end{equation}
where the proof of its more general form is postponed in \eqref{eq:11}. We are now in position to prove \eqref{eqn:bdw1}. The combination of \eqref{eqn:w1green}, \eqref{eqn:greenGM}, and \eqref{BoundvTForMassive} and the definition \eqref{BoundGenericAlgebraic} of $\mathcal{C}$ as well as a polar change of coordinates yields for any $x\in\mathbb{R}^d$
\begin{equation}\label{IntegralBoundW1}
	\mathbb{E}\big[|w_1(x)|^p\big]^{\frac{1}{p}} \lesssim p^{\frac{1}{\gamma}}\sqrt{T}\int_{\R^d}\dd y\, |x-y|^{2-d} \exp\Big(-\frac{|x-y|}{C\sqrt{M}}\Big) \frac{\log(2+\frac{|y|}{r})}{(|y|+r)^d}\quad\text{for any $p<\infty$},
\end{equation}
for some $\gamma>0$. We then decompose the integral into the far-field regime $\vert y\vert\geq 2\vert x\vert$ and the near field regime $\vert y\vert\leq 2\vert x\vert$. 

\medskip

In the far-field regime $|y|\ge 2|x|$ we have $ \frac{|y|}{2} \le |x-y| \le \frac{3}{2} |y|$, so that, after a polar change of coordinates,
\begin{equation}\label{IntegralBoundW1:2}
\int_{\mathbb{R}^d\backslash \bb_{2\vert x\vert}} \dd y\, |x-y|^{2-d} \exp\Big(-\frac{|x-y|}{C\sqrt{M}}\Big) \frac{\log(2+\frac{|y|}{r})}{(|y|+r)^d}\lesssim \int_{2\vert x\vert}^{\infty}\dd \rho\, \rho \exp\Big(-\frac{\rho}{C\sqrt{M}}\Big) \frac{\log(2+\frac{\rho}{r})}{(\rho+r)^d}.
\end{equation}
We further distinguish between the cases $\vert x\vert\leq \sqrt{M}$ and $\vert x\vert\geq \sqrt{M}$. For $\vert x\vert\leq \sqrt{M}$, we neglect the exponential factor in \eqref{IntegralBoundW1:2} and integrate directly into
\begin{align*}
\int_{2\vert x\vert}^{\infty}\dd \rho\, \rho \exp\Big(-\frac{\rho}{C\sqrt{M}}\Big) \frac{\log(2+\frac{\rho}{r})}{(\rho+r)^d}&\leq r^{-d}\vert x\vert^2\int_{2}^{\infty}\dd\rho\, \rho \frac{\log(2+\tfrac{\vert x\vert}{r}\rho)}{(\tfrac{\vert x\vert}{r}\rho+1)^d}\\
&\lesssim 
\left\{
    \begin{array}{ll}
        r^{2-d} & \text{if $\vert x\vert\leq r$}\\
        \vert x\vert^{2-d}\log(2+\tfrac{\vert x\vert}{r})\int_1^{\infty}\frac{\log(2+\rho)}{\rho^{d-1}}& \text{if $\vert x\vert\geq r$}
    \end{array}
\right.\\
&\lesssim M\frac{\log(2+\tfrac{\vert x\vert}{r})}{(\vert x\vert+r)^d}.
\end{align*}
%
%\begin{align*}
%
	%&\int_{\mathbb{R}^d\backslash \bb_{2\vert x\vert}} \dd y \, \mathcal{C}(r,y) |x-y|^{2-d} \exp\Big(-\frac{|x-y|}{\sqrt{M}}\Big)  \frac{\log(2+\frac{|y|}{r})}{(|y|+r)^d}\\ 
	%& \lesssim \int_{\mathbb{R}^d\backslash \bb_{2\vert x\vert}} \dd y\,  \mathcal{C}(r,y)|x-y|^{2-d} \exp\Big(-\frac{|x-y|}{\sqrt{M}}\Big) \frac{\log(2+\frac{|x-y|}{r})}{(|x-y|+r)^d} \\ & \leq \mathcal{C}(r,T,M,x)\int_{|x|}^\infty \dd \rho\, \rho\exp\Big(-\frac{\rho}{\sqrt{M}}\Big) \frac{\log(2+\frac{\rho}{r})}{(\rho+r)^d}\\ &  \leq \mathcal{C}(r,T,M,x) \left\{\begin{aligned}
		%&  \frac{\log(2+\frac{\sqrt{M}}{r})}{(|x|+r)^{d-2}}& \text{for $|x| \le \sqrt{M}$,} \\ 
		%& \sqrt{M}\frac{\log(2+\frac{|x|}{r})}{(|x|+r)^{d-1}} \exp(-\frac{|x|}{\sqrt{M}})  & \text{for $|x|\ge \sqrt{M}$}.
	%\end{aligned} \right. 
%
%\end{align*}
%
For $\vert x\vert\geq \sqrt{M}$, we integrate the exponential factor in form of 
\begin{align*}
\int_{\vert x\vert}^{\infty}\dd \rho\, \rho \exp\Big(-\frac{\rho}{C\sqrt{M}}\Big) \frac{\log(2+\frac{\rho}{r})}{(\rho+r)^d}\leq M(\vert x\vert+r)^{-d}\int_{0}^{\infty} \dd\rho\, \rho\exp(-\frac{\rho}{C})\log(2+\tfrac{\sqrt{M}}{r}\rho)
\lesssim M\frac{\log(2+\tfrac{\vert x\vert}{r})}{(\vert x\vert+r)^d}.
\end{align*}
For the near-field regime $|y|\le 2|x|$, we further distinguish between the cases $\vert x\vert\leq r$ and $\vert x\vert\geq r$. For $\vert x\vert\leq r$, we have $|y| \le 2r$ and $|x-y| \le 3|x|$ thus, neglecting the exponential, the integral in \eqref{IntegralBoundW1} is dominated by
\begin{align*}
\int_{\bb_{2\vert x\vert}} \dd y\, |x-y|^{2-d} \exp\Big(-\frac{|x-y|}{C\sqrt{M}}\Big) \frac{\log(2+\frac{|y|}{r})}{(|y|+r)^d}\lesssim r^{-d}  \int_{\bb_{3\vert x\vert}(x)} \dd y\,|x-y|^{2-d}\lesssim M(\vert x\vert+r)^{-d}.
\end{align*}
For $|x| \ge r$, we further sub-divide into $|y-x|\le \frac{|x|}{2}$ and $|y-x|\ge \frac{|x|}{2}$. When $|y-x|\le \frac{|x|}{2}$ holds, we also have $\frac{|x|}{2}\leq |y| \leq \frac{3}{2}\vert x\vert$ and thus using a polar change of coordinates
\begin{align*}
\int_{\bb_{2\vert x\vert}\cap \bb_{\frac{\vert x\vert}{2}}(x)} \dd y\, |x-y|^{2-d} \exp\Big(-\frac{|x-y|}{C\sqrt{M}}\Big) \frac{\log(2+\frac{|y|}{r})}{(|y|+r)^d}&\lesssim\frac{\log(2+\frac{|x|}{r})}{(|x|+r)^d} \int_{1}^{\frac{\vert x\vert}{2}} \dd \rho\, \rho \exp\Big(-\frac{\rho}{C\sqrt{M}}\Big) \\ 
& \lesssim M\frac{\log(2+\frac{|x|}{r})}{(|x|+r)^d}.
\end{align*}
Finally, when $|y-x|\ge \frac{|x|}{2}$, we have $|y-x|\leq 3|x|$ so that
\begin{align*}
\int_{\bb_{2\vert x\vert}\cap \bb_{\frac{\vert x\vert}{2}(x)}} \dd y\, |x-y|^{2-d} \exp\Big(-\frac{|x-y|}{C\sqrt{M}}\Big) \frac{\sqrt{T}\log(2+\frac{|y|}{r})}{(|y|+r)^d}&\lesssim \vert x\vert^2\exp(-\tfrac{\vert x\vert}{C\sqrt{M}})\frac{\log(2+\tfrac{\vert x\vert}{r})}{(\vert x\vert+r)^{d}}\\
&\lesssim M \frac{\log(2+\tfrac{\vert x\vert}{r})}{(\vert x\vert+r)^{d}}.
\end{align*}
The bound for $\nabla w_1$ can be obtained similarly by replacing $G_M$ with $\nabla G_M$ in \eqref{eqn:w1green}.
\subsection{Proof of Theorem \ref{SecondOrderDecayFlux}: Fluctuations of the $2^{\text{nd}}$-order time-dependent flux}
%%We obtain the stochastic estimates \eqref{2NDOrder:Eq2} in three steps. Introducing the notation $g_r:=r^{-d}e^{-\frac{\vert \cdot\vert^2}{r^2}}$ for the Gaussian kernel, we prove in a first step stochastic estimates on the averaged time dependent flux defined as
%
%\begin{equation}\label{2NDOrder:Eq9}
%(q^{(2)}_M)_r:=\int g_r(a\phi_M^{(1)}-\sigma_M^{(1)})e'+\int\int_{0}^T g_r a\nabla u^{(2)}_M\text{ for any $T<\infty$ and $r<\sqrt{T}$,}
%\end{equation}
%
%in form of\footnote{Throughout the paper we hide the polynomial $p$-dependence inside $\lesssim$, since as long as the $p$-dependence is polynomial (which easily follow throughout our proofs), the random variable on \lhs has exponential moments, c.f. \cite[Lemma 9]{clozeau2021optimal}, \cite[Lemma 4.8]{lu2021optimal}.}
%
%\begin{equation}\label{2NDOrder:Eq6}
%\big\langle \big\vert \big( q^{(2)}_M(T)-\langle q^{(2)}_M(T) \rangle\big)_r\big\vert^p \big\rangle ^\frac{1}{p}\lesssim r^{-\frac{d}{2}}\sqrt{T} \log(T)\log^2 (\frac{\sqrt{T}}{r})\mu_{\beta}(T).
%\end{equation}
%
%Here, our main tools are the spectral gap inequality \eqref{msPIp}, the energy estimates in Lemma \ref{lem:GOlm1}, as well as the large-scale regularity results of Proposition \ref{RegResultElliptic} and Proposition \ref{RegResultParabolic}. In the second step, we obtain the stochastic estimates \eqref{2NDOrder:Eq2} by combining \eqref{2NDOrder:Eq6} with \cite[Lemma 6]{gloria2015corrector} \eqref{} \nc{lemma which converts weak to strong, put in the appendix}. We start with preliminary estimates. 
Since the structure of $(\nabla \phi^{(2)}_M(T,\cdot))_r$ is similar to that of $(q_M^{(2)}(T,\cdot))_r$, we only prove the fluctuations of the later. We split the proof into four steps. In the first step, we provide preliminary estimates which will serve as a toolbox for the rest of the proof. It summarizes the large-scale regularity estimates established in Proposition \ref{RegResultParabolic} and Corollary \ref{RegResultElliptic} while applied to the equations \eqref{2NDOrder:Eq1Regularized}, \eqref{eq:vT}, \eqref{eqn:phi1dual} and \eqref{eqn:2ndphi1dual}. In the second step, we treat the $(\sigma^{(1)}_{i,M}e_j)_r$ contribution in \eqref{2NDOrder:Eq8} in a direct way using the fluctuation estimates of its gradient in \cite[Corollary 2]{clozeau2021optimal}. In the third step, we control the functional derivative of $(q^{(2)}_{ij,M,\varepsilon}(T,\cdot)+\sigma^{(1)}_M e_j)_r$, for which a representation formula has been established in \eqref{eqn:funcderiv}, uniformly in the regularizion parameter $\varepsilon$. We treat each terms in \eqref{eqn:funcderiv} separately and argue differently for the two regimes $\ell\leq \sqrt{T}$ and $\ell\geq \sqrt{T}$, which splits naturally this step into two sub-steps. In the fourth step, we combine the bound on the functional derivative with the spectral gap inequality \eqref{msPIp}, yielding the control of moments \eqref{SecondOrderDecayFlux:Eq} for the regularized quantity, for which we can pass to the limit as $\varepsilon\downarrow 0$ thanks to the regularity property \eqref{SmoothnessCoef} of the coefficient field.

\medskip

For notational convenience, we drop the dependence on the indices $ij$. We also recall the notation $\mathcal{C}$ given in \eqref{MomentBoundGeneric} which will be used all along the proof.

\medskip

{\sc Step 1. Preliminary estimates. }In this step we state the regularity estimates we have on the solutions of \eqref{2NDOrder:Eq1Regularized}, \eqref{eq:vT}, \eqref{eqn:phi1dual} and \eqref{eqn:2ndphi1dual} that we use in the next steps. 

\medskip

\textbf{Energy estimates. }Applying the localized energy estimate \eqref{EnergyEstiAppendix:Eq1} with $T=R^2$ to the equation \eqref{2NDOrder:Eq1Regularized}, together with \eqref{SmoothnessCoef} and \eqref{eq:phiptws}, yields for any $R>0$ and $x\in\mathbb{R}^d$
\begin{equation}\label{eq:9} 
\begin{aligned}
&\int_{0}^{R^2}\dd s\int_{\R^d}\eta_R(\tfrac{\cdot-x}{C})\vert\nabla u_{M,\varepsilon}^{(2)}(s,\cdot)\vert^2\stackrel{\eqref{EnergyEstiAppendix:Eq1}}{\lesssim}\int_{\R^d}\eta_R(\tfrac{\cdot-x}{C})\vert\nabla\cdot [(a_\varepsilon\,\phi_M^{(1)}-\sigma_M^{(1)})e_j]\vert^2 \stackrel{\eqref{SmoothnessCoef},\eqref{eq:phiptws}}{\leq} \mathcal{C}(R,x),
\end{aligned}
\end{equation}
where we recall the notation $\eta_R:=R^{-d}\exp(-\frac{\vert \cdot\vert}{R})$. In addition, the localized energy estimates \eqref{GOlm1} combined with \eqref{RegA} and \eqref{eq:phiptws} gives for any $R>0$
\begin{equation}\label{EnergyIntInside}
\int_{\R^d}\eta_{R}(\tfrac{\cdot-x}{C})\bigg\vert\int_{0}^{R^2}\dd s\,\nabla u^{(2)}_{M,\varepsilon}(s,\cdot)\bigg\vert^2\stackrel{\eqref{GOlm1}}{\lesssim} \int_{\R^d}\eta_R(\tfrac{\cdot-x}{C})\vert(a_\varepsilon\,\phi_M^{(1)}-\sigma_M^{(1)})e_j\vert^2\stackrel{\eqref{RegA},\eqref{eq:phiptws}}{\leq} \mathcal{C}(R,x),
\end{equation}
and for any $R\geq T$
\begin{equation}\label{EnergyToUpgrade}
	\int_{\R^d} \eta_{R}(\tfrac{\cdot-x}{C})\vert \nabla u_{M,\varepsilon}^{(2)}(T,\cdot)\vert^2 \lesssim T^{-1}\int_{\R^d}\eta_{R}(\tfrac{\cdot-x}{C})\vert(a_\varepsilon\,\phi_M^{(1)}-\sigma_M^{(1)})e_j\vert^2\leq \mathcal{C}(T,x)\, T^{-1}.
\end{equation}
Furthermore, we upgrade \eqref{EnergyToUpgrade} using Schauder's estimates  and the large-scale regularity \eqref{LargeScaleC01Para}  into
\begin{equation}\label{eq:10}
\vert \nabla u_{M,\varepsilon}^{(2)}(T,x)\vert\leq \mathcal{C}(T,x)\, T^{-1}\quad\text{for any $T>0$ and $x\in\mathbb{R}^d$.}
\end{equation}
Indeed, while \eqref{eq:10} immediately follows from \eqref{SchauderParaEsti} and \eqref{EnergyToUpgrade} in the regime $\sqrt{T}\leq r_\star (x)$, in the regime $\sqrt{T}> r_\star(x)$ we appeal to \eqref{SchauderParaEsti} together with the large-scale regularity \eqref{LargeScaleC01Para} (with an additional time reflection $t\mapsto -t$ and change of kernel $\fint_{\bb_{R}} \rightarrow \int_{\R^d}\eta_R$, which can be done thanks to \cite[(94)]{lu2021optimal}) and \eqref{EnergyToUpgrade} in form of 
\begin{align*}
\vert \nabla u_{M,\varepsilon}^{(2)}(T,x)\vert&\stackrel{\eqref{SchauderParaEsti}}{\lesssim}\mathcal{C}(x)r_\star^d(x)\fint_{T-r^2_\star(x)}^T\dd\tau \int_{\R^d}\eta_{r_\star(x)}(\cdot-x)\vert \nabla u_{M,\varepsilon}^{(2)}(\tau,\cdot)\vert^2\\
&\stackrel{\eqref{LargeScaleC01Para}}{\leq} \mathcal{C}(x)\fint_{0}^T\dd\tau\int_{\R^d}\eta_{\sqrt{T}}(\cdot-x)\vert \nabla u_{M,\varepsilon}^{(2)}(\tau,\cdot)\vert^2 \stackrel{\eqref{EnergyToUpgrade}}{\leq}\mathcal{C}(T,x)T^{-1}.
\end{align*}
\textbf{Large-scale regularity. }As a consequence of the large-scale regularity estimate \eqref{LSreg:PointwiseParabolic}, the localized energy estimate \eqref{EnergyEstiAppendix:Eq1} and Schauder's estimates in Proposition \ref{SchauderTheoryStandard}, we have for any $(t,x)\in (-\infty,T]\times \mathbb{R}^d$
\begin{equation}\label{eq:11}
	 \vert v^T(t,x)\vert\leq \mathcal{C}(r,x)\frac{\sqrt{T-t}\log(2+\tfrac{\vert x\vert}{r})}{(\vert x\vert+r)^d}\quad\text{and}\quad\vert\nabla v^T(t,x)\vert\leq \mathcal{C}(r,x)\frac{\log(2+\tfrac{\vert x\vert}{r})}{(\vert x\vert+r)^{d}}.
\end{equation}
Indeed, we first note that $v^T=(v_k(\cdot-T,\cdot))_{k\in\{1,\cdots,d\}}$ where $v_k$ solves \eqref{EquationParaLargeScale} with $\bar{e}=e_k$ and $f_r=g_r$. In the regime $\sqrt{T-t}\geq 2r_\star(x)$, \eqref{eq:11} follows directly from \eqref{LSreg:PointwiseParabolic} applied to $v_k$. In the regime $\sqrt{T-t}\leq 2r_\star(x)$ we further notice that, as for \eqref{LSreg:PointwisePara17},
\begin{equation}\label{IntegralVersionvT}
v^{T}(t,\cdot)=\int_{0}^{T-t}\dd \tau\,w(\tau,\cdot)\quad\text{for any $t<T$,}
\end{equation}
where $w:=(w_k)_{k\in\{1,\cdots,d\}}$ solves 
\begin{equation}\label{AuxiliaryEquationFluctu}
\left\{
    \begin{array}{ll}
        \partial_\tau w_k-\nabla\cdot a\nabla w_k=0 & \text{$(0,\infty)\times \mathbb{R}^d$,}\\
        \omega_k(0)=-\nabla\cdot a\,g_r e_k.&
    \end{array}
\right.
\end{equation}
Therefore, applying the energy estimate \eqref{GOlm1} on \eqref{AuxiliaryEquationFluctu} gives
\begin{align*}
	\fint_{\bb_{2r_\star(x)}(x)} \dd y\big|(\tfrac{1}{\sqrt{T-t}} v^T(t,y), \nabla v^T(t,y))\big|^2&\stackrel{\eqref{IntegralVersionvT}}{=}\fint_{\bb_{2r_\star(x)}(x)}\dd y\bigg|\int_{0}^{T-t}\dd \tau\,(\tfrac{1}{\sqrt{T-t}} w(\tau,\cdot), \nabla w(\tau,\cdot))\bigg|^2\\
& \lesssim \int_{\R^d}\eta_{2r_\star(x)}(x-\cdot) |g_r(y)|^2  \lesssim \frac{1}{(|x|+r)^d},
\end{align*}
which, combined with Schauder's estimates in Proposition \ref{SchauderTheoryStandard}, yields \eqref{eq:11} in the regime $\sqrt{T-t}\leq 2r_\star(x)$.

\medskip

Finally, recalling the equation \eqref{eqn:2ndphi1dual} for $w_2$ together with the bounds \eqref{eq:11} for $v^T(0,x)$ as well as \eqref{eqn:bdw1} for $w_1$, $\frac{w_2}{\sqrt{T}}$ satisfies the assumptions of (ii) in Corollary \ref{RegResultElliptic} from which we deduce 
\begin{equation}\label{eqn:bdnablaw2}
	|\nabla w_2(x)| \le \mathcal{C}(r,x)\sqrt{T}\frac{\log^2(2+\tfrac{\vert x\vert}{r})}{(\vert x\vert+r)^{d}}\left(1 \vee\frac{\vert x\vert+r}{\sqrt{T}}\right).
\end{equation}
{\sc Step 2. Fluctuations of $(\sigma_M^{(1)}e_j)_r$. }Proof that 
\begin{equation}\label{FluctuationsSigmaM}
\vert (\sigma^{(1)}_M e_j)_r(x)\vert\leq \mathcal{C}(x,r)r^{-\frac{d}{2}}\mu_\beta(T)\quad\text{for any $x\in\mathbb{R}^d$},
\end{equation}
where we recall that $\mu_\beta$ is defined in \eqref{eq:defmubeta}. Since $\sigma^{(1)}_M$ is stationary, we fix $x=0$ without loss of generality and simply write $(\sigma^{(1)}_M e_j)_r$ for its evaluation at $0$.

\medskip

The bound \eqref{FluctuationsSigmaM} is a direct consequence of the fluctuation bounds on the gradient given in \cite[Corollary 2]{clozeau2021optimal} which reads
\begin{equation}\label{FluctuationsBoundGradientM}
	\vert\nabla( \sigma^{(1)}_M e_j)_{\sqrt{\tau}}\vert \leq \mathcal{C}(x,\tau)\tilde{\mu}_\beta(\tau)\quad\text{for any $\tau\geq 1$ and $x\in\mathbb{R}^d$,}
\end{equation}
with 
\begin{equation}
\tilde{\mu}_\beta(\tau):=\left\{
    \begin{array}{ll}
        \tau^{-\frac{\beta}{4}} & \text{for $\beta<d$,} \\
        \tau^{-\frac{d}{4}}\log^{\frac{1}{2}}(\tau) & \text{for $\beta=d$,}\\
        \tau^{-\frac{d}{4}} & \text{for $\beta>d$.}
    \end{array}
\right.
\end{equation}
Indeed, since $\mathbb{E}[\sigma^{(1)}_M]=0$ and $\sigma^{(1)}_M$ is stationary, one has by ergodicity $(\sigma^{(1)}_{M})_{\sqrt{\tau}}\,e_j\underset{\tau\uparrow\infty}{\rightarrow}0$. Thus, using in addition that $\partial_\tau g_{\tau}=\Delta g_{\tau}$ for any $\tau>0$ and $\|\sigma^{(1)}_M\|_{\LL^{\infty}(\mathbb{R}^d)}\lesssim_M 1$ (as a consequence of the localized energy estimate \eqref{LocalizedMassiveTerm} applied to \eqref{eq:1stOrderMassFlux} and \eqref{eq:1stOrderMassCor} successively together with Schauder's estimate), we can perform an integration by integration by parts to obtain
\begin{equation*}
	(\sigma^{(1)}_{M})_{r}e_j=-\int_{r^2}^{\infty}\dd\tau\int_{\R^d}\partial_{\tau}g_{\sqrt{\tau}}\,\sigma_M^{(1)}e_j=-\int_{r^2}^{\infty}\dd\tau\int_{\R^d}\Delta g_{\sqrt{\tau}}\sigma_M^{(1)}e_j=\int_{r^2}^{\infty}\dd\tau\int_{\R^d} \nabla g_{\sqrt{\tau}}\cdot \nabla (\sigma_M^{(1)}e_j).
\end{equation*}
Now, using the semigroup property $g_{\sqrt{\tau}}=g_{\sqrt{\tfrac{\tau}{2}}}\star g_{\sqrt{\tfrac{\tau}{2}}}$, we arrive at 
$$(\sigma^{(1)}_{M})_{r}e_j=\int_{r^2}^{+\infty} \nabla( \sigma^{(1)}_M e_j)_{\sqrt{\tfrac{\tau}{2}}}\star \nabla g_{\sqrt{\tfrac{\tau}{2}}}(0)\, \dd\tau,$$
so that using \eqref{FluctuationsBoundGradientM}, we get
$$\vert(\sigma^{(1)}_{M})_{r}e_j\vert\leq \int_{r^2}^{\infty}\mu_\beta\Big(\sqrt{\tfrac{\tau}{2}}\Big)\mathcal{C}(\tau,\cdot)\star \vert\nabla g_{\sqrt{\frac{\tau}{2}}}\vert(0).$$
Using finally that $\|\nabla g_{\sqrt{\tfrac{\tau}{2}}}\|_{\LL^1(\mathbb{R}^d)}\lesssim \tau^{-\frac{1}{2}}$ and  $r\leq \sqrt{T}$, we deduce that 
$$\int_{r^2}^{\infty}\mu_\beta\Big(\sqrt{\tfrac{\tau}{2}}\Big)\mathcal{C}(\tau,\cdot)\star \vert\nabla g_{\sqrt{\frac{\tau}{2}}}\vert(0)\leq \mathcal{C}(r)r^{-\frac{d}{2}}\mu_\beta(T),$$
which yields \eqref{FluctuationsSigmaM}.

\medskip

{\sc Step 3. Control of the functional derivative $\partial_{\cdot,\ell}\, (q^{(2)}_{M,\varepsilon}+\sigma^{(1)}_M e_j)_r$. }We split this step into two sub-steps treating the regimes $\ell\leq \sqrt{T}$ and $\ell\geq \sqrt{T}$ separately.

\medskip

{\sc Sub-step 3.1. Regime $\ell\leq \sqrt{T}$. }Proof that for any $k\in \{1,\cdots, d\}$
\begin{equation}\label{ControlFunctionalEllSmall}
\begin{aligned}
\|\partial^{\text{fct}}_{\cdot,\ell}\,(q^{(2)}_{M,\varepsilon}(T,\cdot)+\sigma^{(1)}_{M}e_j)_r\cdot e_k\|^2_{\LL^2(\mathbb{R}^d)}\leq \mathcal{C}(\ell,M,r,T,\varepsilon)\, T\ell^{2d} r^{-d}(\log^2(\tfrac{\sqrt{T}}{r})+\log^2(\tfrac{\sqrt{T}}{\ell})).
\end{aligned}
\end{equation}
The bound \eqref{ControlFunctionalEllSmall} is a consequence of the combination of the following five estimates \eqref{FirstRHSDerivative}, \eqref{SecondRHSDerivative}, \eqref{ThirdRHSDerivative}, \eqref{FourthRHSDerivative} and \eqref{FifithSixthRHSDerivative}, treating separately each terms of \eqref{eqn:funcderiv} with $\nu=e_k$.

\medskip

\textbf{First \rhs term of \eqref{eqn:funcderiv}. }Using Cauchy-Schwarz' inequality, \eqref{eq:phiptws} and $\|g_r\|_{\LL^{2}(\mathbb{R}^d)}^2\lesssim r^{-d}$, we directly obtain
\begin{equation}\label{FirstRHSDerivative}
\begin{aligned}
	\int_{\R^d}\dd x \bigg(\int_{\bb_\ell(x)} \vert\phi^{(1)}_M\, g_r\,e_j\otimes e_k\vert \bigg)^2  \lesssim \ell^{2d}\int_{\R^d} g_r^2 |\phi_M^{(1)}|^2 &\stackrel{\eqref{eq:phiptws}}{\lesssim} \ell^{2d} 	\int_{\R^d}\C g_r^2\\
	&\leq \mathcal{C}(r)\,\ell^{2d}r^{-d}. 
\end{aligned}
\end{equation}

\textbf{Second \rhs term of \eqref{eqn:funcderiv}. }We claim that 
\begin{equation}\label{SecondRHSDerivative}
\int_{\R^d}\dd x  \bigg(\int_{\bb_{\ell}(x)}g_r\bigg\vert\int_{0}^T\dd s\,\nabla u_{M,\varepsilon}^{(2)}(s,\cdot)\bigg\vert\bigg)^2\leq \mathcal{C}(\ell,r,T,\varepsilon)\ell^{2d} r^{-d}(1+\log^2(\tfrac{\sqrt{T}}{r})).
\end{equation}
We split the time integral into two parts $\int_{0}^{T}=\int_{0}^{r^2}+\int_{r^2}^T$ and we use the triangle inequality to the effect of 
\begin{equation}\label{SplitFirstDeriveSmall}
\begin{aligned}
\int_{\R^d}\dd x & \biggl(\int_{\bb_{\ell}(x)}g_r\bigg\vert\int_{0}^T\dd s\,\nabla u_{M,\varepsilon}^{(2)}(s,\cdot)\bigg\vert\bigg)^2  \\ & \lesssim 	\int_{\R^d}\dd x \bigg(\int_{\bb_{\ell}(x)} g_r\bigg\vert\int_{0}^{r^2}\dd s\,\nabla u_{M,\varepsilon}^{(2)}(s,\cdot)\bigg\vert\bigg)^2 +	\int_{\R^d}\dd x\bigg(\int_{\bb_{\ell}(x)}g_r\bigg\vert\int_{r^2}^T\dd s\nabla u_{M,\varepsilon}^{(2)}(s,\cdot)\bigg\vert\bigg)^2.
\end{aligned}
\end{equation}
For the first \rhs term of \eqref{SplitFirstDeriveSmall}, we use Cauchy-Schwarz' inequality, the fact that $g^2_r\lesssim r^{-d}\eta_r$ and \eqref{EnergyIntInside} to obtain
\begin{equation}\label{SecondTermCompute1}
\begin{aligned}
\int_{\R^d}\dd x \Bigl(\int_{\bb_{\ell}(x)}g_r\Big\vert\int_{0}^{r^2}\dd s\nabla u_{M,\varepsilon}^{(2)}(s,\cdot)\Big\vert\Bigr)^2\lesssim& \ell^{2d}\int_{\R^d}g^2_r\bigg\vert\int_{0}^{r^2}\dd s\nabla u_{M,\varepsilon}^{(2)}(s,\cdot)\bigg\vert^2\\
\lesssim &  \ell^{2d}r^{-d}\int_{\R^d} \eta_r\bigg\vert\int_{0}^{r^2}\dd s\nabla u_{M,\varepsilon}^{(2)}(s,\cdot)\bigg\vert^2\\
\stackrel{\eqref{EnergyIntInside}}{\leq} & \mathcal{C}(r)\ell^{2d}r^{-d}.
\end{aligned}
\end{equation}
For the second \rhs term \eqref{SplitFirstDeriveSmall}, we proceed the same way using \eqref{eq:10} instead
\begin{equation}\label{SecondTermCompute2}
\begin{aligned}
\int_{\R^d}\dd x\bigg(\int_{\bb_{\ell}(x)}g_r\bigg\vert\int_{r^2}^T\dd s\,\nabla u_{M,\varepsilon}^{(2)}(s,\cdot)\bigg\vert\bigg)^2\stackrel{\eqref{eq:10}}{\lesssim} \ell^{2d}\int_{\R^d}g^2_r\bigg(\int_{r^2}^{T}\dd s\, s^{-1}\mathcal{C}(s,\cdot)\bigg)^2\leq \mathcal{C}(r,T) \ell^{2d}r^{-d}\log^2(\tfrac{\sqrt{T}}{r}).
\end{aligned}
\end{equation}

%We expect the scaling of the first four terms of \eqref{eqn:funcderiv} to be smaller, since they also appear in \cite{clozeau2021optimal} with $u^{(1)}$ replaced by $u_M^{(2)}$, which, in view of \eqref{eq:9}, \eqref{eq:10}, has the same deterministic scaling as $u^{(1)}$, and the proof strategies are identical. The latter two terms of \eqref{eqn:funcderiv}, however, have the leading order behavior and therefore require new techniques.
%
\textbf{Third \rhs term of \eqref{eqn:funcderiv}. }We claim that 
\begin{equation}\label{ThirdRHSDerivative}
\int_{\R^d}\dd x\bigg(\int_{\bb_\ell(x)}\vert (\nabla v^{T}(0,\cdot)\otimes \phi^{(1)}_M e_j)\star\chi_\varepsilon\vert\bigg)^2\leq \mathcal{C}(\ell,r)\ell^{2d}r^{-d}.
\end{equation}
Using Cauchy-Schwarz' inequality, the second item of \eqref{eq:11} and \eqref{eq:phiptws}, we have
\begin{align*}
\int_{\R^d}\dd x\bigg(\int_{\bb_\ell(x)}\vert(\nabla v^{T}(0,\cdot)\otimes \phi^{(1)}_M e_j)\star\chi_\varepsilon\vert\bigg)^2\leq& \int_{\R^d} \dd x \int_{\bb_\ell(x)}\vert \nabla v^T(0,\cdot)\vert^2\star \chi_\varepsilon\int_{\bb_{\ell}(x)}\vert\phi^{(1)}_M\vert^2\star\chi_\varepsilon\\
\stackrel{\eqref{eq:11},\eqref{eq:phiptws}}{\leq} &\ell^{2d}\int_{\R^d} \dd x \fint_{\bb_\ell(x)}\mathcal{C}^2(r,\cdot)\frac{\log^2(2+\tfrac{\vert\cdot\vert}{r})}{(\vert\cdot\vert+r)^{2d}}\star\chi_\varepsilon\fint_{\bb_\ell(x)}\mathcal{C}^2\star\chi_\varepsilon.
\end{align*}
We then check that
\begin{equation}\label{FluctuationCheckSto}
\mathcal{I}_{\ell,r}:=\int_{\R^d} \dd x \fint_{\bb_\ell(x)}\mathcal{C}^2(r,\cdot)\frac{\log^2(2+\tfrac{\vert\cdot\vert}{r})}{(\vert\cdot\vert+r)^{2d}}\star\chi_\varepsilon\fint_{\bb_\ell(x)}\mathcal{C}^2\star\chi_\varepsilon=r^{-d}\mathcal{C}(\ell,r,\varepsilon).
\end{equation}
We obtain \eqref{FluctuationCheckSto} using the equivalent definition \eqref{BoundGenericAlgebraic}. For any $p<\infty$, we have by Minkowski's inequality and Cauchy-Schwarz' inequality 
$$\mathbb{E}[\mathcal{I}^p_{\ell,r}]^{\frac{1}{p}}\leq \int_{\R^d} \dd x\fint_{\bb_\ell(x)}\mathbb{E}[\mathcal{C}^{4p}(r,\cdot)]^{\frac{1}{2p}}\frac{\log^2(2+\tfrac{\vert\cdot\vert}{r})}{(\vert\cdot\vert+r)^{2d}}\star\chi_\varepsilon\fint_{\bb_{\ell}(x)}\mathbb{E}[\mathcal{C}^{4p}]^{\frac{1}{2p}}\star\chi_\varepsilon.$$
Using then \eqref{BoundGenericAlgebraic}, there exist $\gamma>0$ and $C>0$ such that 
$$\mathbb{E}[\mathcal{I}^p_{\ell,r}]^{\frac{1}{p}}\leq C^4 p^{\frac{4}{\gamma}}\int_{\R^d}\dd x\fint_{\bb_{\ell}(x)}\frac{\log^2(2+\tfrac{\vert\cdot\vert}{r})}{(\vert\cdot\vert+r)^{2d}}\star\chi_\varepsilon=C^4 p^{\frac{4}{\gamma}}\int_{\R^d}\frac{\log^2(2+\tfrac{\vert\cdot\vert}{r})}{(\vert\cdot\vert+r)^{2d}}\lesssim C^4 p^{\frac{4}{\gamma}} r^{-d},$$
which shows \eqref{FluctuationCheckSto}.

\medskip

\textbf{Fourth \rhs term of \eqref{eqn:funcderiv}. }We claim that 
\begin{equation}\label{FourthRHSDerivative}
\int\dd x\bigg(\int_{\bb_\ell(x)}\bigg\vert\int_{0}^T\dd s\,\nabla v^{T}(s,\cdot)\otimes \nabla u^{(2)}_{M,\varepsilon}(s,\cdot)\bigg\vert\bigg)^2\leq \mathcal{C}(\ell,r,\varepsilon)\ell^{2d}r^{-d}(\ell^2+\log^2(\tfrac{\sqrt{T}}{\ell})).
\end{equation}
We split the time integral into two contributions using the triangle inequality
\begin{equation}\label{SplitThirdTermSmall}
\begin{aligned}
&\int_{\R^d}\dd x\bigg(\int_{\bb_\ell(x)}\bigg\vert\int_{0}^T\dd s\,\nabla v^{T}(s,\cdot)\otimes \nabla u^{(2)}_{M,\varepsilon}(s,\cdot)\bigg\vert\bigg)^2\\
&\lesssim \int_{\R^d}\dd x\bigg(\int_{\bb_\ell(x)}\bigg\vert\int_{0}^{\ell^2}\dd s\,\nabla v^{T}(s,\cdot)\otimes \nabla u^{(2)}_{M,\varepsilon}(s,\cdot)\bigg\vert\bigg)^2+\int_{\R^d}\dd x\bigg(\int_{\bb_\ell(x)}\bigg\vert\int_{\ell^2}^T\dd s\,\nabla v^{T}(s,\cdot)\otimes \nabla u^{(2)}_{M,\varepsilon}(s,\cdot)\bigg\vert\bigg)^2.
\end{aligned}
\end{equation}
For the first \rhs term of \eqref{SplitThirdTermSmall}, we use Cauchy-Schwarz' inequality combined with \eqref{eq:9} and the second-item of \eqref{eq:11} to obtain
\begin{align*}
\int_{\R^d}\dd x\bigg(\int_{\bb_\ell(x)}\bigg\vert\int_{0}^{\ell^2}\dd s\,\nabla v^{T}(s,\cdot)\otimes \nabla u^{(2)}_{M,\varepsilon}(s,\cdot)\bigg\vert\bigg)^2&\lesssim \ell^{2d}\int_{\R^d}\dd x\fint_{\bb_{\ell}(x)}\int_0^{\ell^2}\vert\nabla v^T(s,\cdot)\vert^2\fint_{\bb_{\ell}(x)}\int_0^{\ell^2}\vert\nabla u^{(2)}_{M,\varepsilon}(s,\cdot)\vert^2\\
&\stackrel{\eqref{eq:9},\eqref{eq:11}}{\leq} \ell^{2(d+1)}\int_{\R^d}\dd x\,\mathcal{C}(\ell,\varepsilon,x)\fint_{\bb_{\ell}(x)}\frac{\log^2(2+\tfrac{\vert\cdot\vert}{r})}{(\vert\cdot\vert+r)^{2d}}\\
&\leq \mathcal{C}(\ell,\varepsilon) \ell^{2(d+1)}r^{-d},
\end{align*}
where for the last line we proceed as for \eqref{FluctuationCheckSto}. For the second \rhs term of \eqref{SplitThirdTermSmall}, we instead appeal to \eqref{eq:10} and the second-item of \eqref{eq:11} to obtain
\begin{align*}
\int_{\R^d}\dd x\bigg(\int_{\bb_\ell(x)}\bigg\vert\int_{\ell^2}^T\dd s\,\nabla v^{T}(s,\cdot)\otimes \nabla u^{(2)}_{M,\varepsilon}(s,\cdot)\bigg\vert\bigg)^2 &\stackrel{\eqref{eq:10},\eqref{eq:11}}{\leq} \ell^{2d}\int_{\R^d}\bigg(\int_{\ell^2}^{T}\dd s\,s^{-1}\mathcal{C}(s,\cdot)\bigg)^2\mathcal{C}^2(r,\cdot)\frac{\log^2(2+\tfrac{\vert\cdot\vert}{r})}{(\vert\cdot\vert+r)^{2d}}\\
&\leq \mathcal{C}(\ell,T,r)\ell^{2d}r^{-d}\log^2(\tfrac{\sqrt{T}}{\ell}).
\end{align*}
\textbf{Fifth and Sixth \rhs terms of \eqref{eqn:funcderiv}. }We claim that
\begin{equation}\label{FifithSixthRHSDerivative}
\begin{aligned}
\int_{\R^d}\dd x\bigg(\int_{\bb_{\ell}(x)}\vert\nabla w_2\otimes (e_i+\nabla \phi^{(1)}_M)\vert\bigg)^2&+\int_{\R^d}\dd x\bigg(\int_{\bb_{\ell}(x)}\vert (v^T(0,\cdot)e_j-\tfrac{1}{M}w_1)\otimes (e_i+\nabla\phi^{(1)}_M)\vert\bigg)^2\\
&\leq \mathcal{C}(\ell,r)T\ell^{2d}r^{-d}.
\end{aligned}
\end{equation}
Since the bounds of $v^T$, $\frac{1}{M}w_1$ are smaller than that of $\nabla w_2$, see \eqref{eq:11}, \eqref{eqn:bdw1} and \eqref{eqn:bdnablaw2}, we only treat the contribution from the first \lhs term of \eqref{FifithSixthRHSDerivative}. Combining Cauchy-Schwarz' inequality with \eqref{eqn:bdnablaw2} and \eqref{eq:phiptws}, we have
\begin{align*}
\int_{\R^d}\dd x\bigg(\int_{\bb_{\ell}(x)} \vert\nabla w_2\otimes (e_i+\nabla \phi^{(1)}_M)\vert\bigg)^2 &\leq T\ell^{2d}\int_{\R^d}\dd x\, \mathcal{C}(r,\cdot)\frac{\log^4(2+\tfrac{\vert\cdot\vert}{r})}{(\vert\cdot\vert+r)^{2d}}\Big(1+\frac{(\vert\cdot\vert+r)^2}{T}\Big)\\
&\leq \mathcal{C}(\ell,r)T\ell^{2d}r^{-d},
\end{align*}
where, for the last line, we recall that $T\geq r^{2}$ and  proceed as for \eqref{FluctuationCheckSto}.

\medskip

{\sc Substep 3.2. Regime $\ell\geq \sqrt{T}$. } We claim that for any $k\in\{1,\cdots,d\}$
\begin{equation}\label{ControlFunctionalEllBig}
\|\partial^{\text{fct}}_{\cdot,\ell}\,(q^{(2)}_{M,\varepsilon}(T,\cdot)+\sigma^{(1)}_{M}e_j)_r\cdot e_k\|^2_{\LL^2(\mathbb{R}^d)}\leq \mathcal{C}(\ell,M,r,T,\varepsilon)\, \ell^{d}\big(\ell^2+\log^6(\tfrac{\ell}{r})+\log^2(\tfrac{\sqrt{T}}{r})\big).
\end{equation}
The bound \eqref{ControlFunctionalEllBig} is a consequence of the combination of the following five estimates \eqref{FirstRHSDerivativeLarge}, \eqref{SecondRHSDerivativeLarge}, \eqref{ThirdRHSDerivativeLarge}, \eqref{FourthRHSDerivativeLarge} and \eqref{FifithSixthRHSDerivativeLarge}, treating separately each terms of \eqref{eqn:funcderiv} with $\nu=e_k$.

\medskip

\textbf{First \rhs term of \eqref{eqn:funcderiv}. }Using Jensen's inequality for the measure $g_r\dd x$ and \eqref{eq:phiptws} directly gives 
\begin{equation}\label{FirstRHSDerivativeLarge}
\begin{aligned}
\int_{\R^d}\dd x \bigg(\int_{\bb_\ell(x)} \vert\phi_M^{(1)}g_r e_j\otimes e_k\vert \bigg)^2\lesssim \int_{\R^d}\dd x\int_{\bb_{\ell}(x)}g_r\vert\phi^{(1)}_M\vert^2\stackrel{\eqref{eq:phiptws}}{\lesssim} \ell^d\int_{\R^d}\mathcal{C}g_r\leq \mathcal{C}(r)\ell^d.
\end{aligned}
\end{equation}
\textbf{Second \rhs term of \eqref{eqn:funcderiv}. }We claim that 
\begin{equation}\label{SecondRHSDerivativeLarge}
\int_{\R^d}\dd x  \bigg(\int_{\bb_{\ell}(x)}g_r\bigg\vert\int_{0}^T\dd s\,\nabla u_{M,\varepsilon}^{(2)}(s,\cdot)\bigg\vert\bigg)^2\leq \mathcal{C}(r)\ell^d\log^2(\tfrac{\sqrt{T}}{r}).
\end{equation}
We apply Jensen's inequality for the measure $g_r\dd x$ to get 
\begin{align*}
\int_{\R^d}\dd x  \bigg(\int_{\bb_{\ell}(x)}g_r\bigg\vert\int_{0}^T\dd s\,\nabla u_{M,\varepsilon}^{(2)}(s,\cdot)\bigg\vert\bigg)^2&\lesssim \int_{\R^d}\dd x\int_{\bb_{\ell}(x)} g_r \bigg\vert\int_{0}^T\dd s\,\nabla u_{M,\varepsilon}^{(2)}(s,\cdot)\bigg\vert^2 \lesssim \ell^d\int_{\R^d}g_r\bigg\vert\int_{0}^T\dd s\,\nabla u_{M,\varepsilon}^{(2)}(s,\cdot)\bigg\vert^2,
\end{align*}
and then \eqref{SecondRHSDerivativeLarge} follows from the same computations as in \eqref{SecondTermCompute1} and \eqref{SecondTermCompute2}.

\medskip

\textbf{Third \rhs term of \eqref{eqn:funcderiv}. }We claim that 
\begin{equation}\label{ThirdRHSDerivativeLarge}
\int_{\R^d}\dd x\bigg(\int_{\bb_{\ell}(x)}\vert (\nabla v^{T}(0,\cdot)\otimes \phi^{(1)}_M e_j)\star \chi_\varepsilon\vert\bigg)^2\leq \mathcal{C}(\ell,r)\ell^d \log^{4}(\tfrac{\ell}{r}).
\end{equation}
We only treat the case when $\ell\geq r_\star(0)$, in the other regime we simply use \eqref{ThirdRHSDerivative} and we bound one $\ell^d$ by $r_\star^d(0)$ together with \eqref{Momentrstar}. Using the second-item of \eqref{eq:11} and \eqref{eq:phiptws}, we have 
\begin{equation}\label{ThirdTermContribution1Bis}
\int_{\R^d}\dd x\bigg(\int_{\bb_{\ell}(x)}\vert (\nabla v^{T}(0,\cdot)\otimes \phi^{(1)}_M e_j)\star \chi_\varepsilon\vert\bigg)^2\leq \int_{\R^d}\dd x\bigg(\int_{\bb_{\ell}(x)}\mathcal{C}\frac{\log(2+\tfrac{\vert\cdot\vert}{r})}{(\vert\cdot\vert+r)^d}\star \chi_\varepsilon\bigg)^2.
\end{equation}
Then, note that from Jensen's inequality with the measure $\chi_\varepsilon\dd x$, we have
\begin{align*}
\int_{\R^d}\dd x\bigg(\int_{\bb_{\ell}(x)}\mathcal{C}\frac{\log(2+\tfrac{\vert\cdot\vert}{r})}{(\vert\cdot\vert+r)^d}\star \chi_\varepsilon\bigg)^2&=\int_{\R^d}\dd x\bigg(\int\dd z\,\chi_{\varepsilon}(z)\int_{\bb_{\ell}(z-x)}\dd y\,\mathcal{C}(y)\frac{\log(2+\tfrac{\vert y\vert}{r})}{(\vert y\vert+r)^d}\bigg)^2\\
&\leq \int_{\R^d} \dd x\bigg(\int_{\bb_{\ell}(x)}\mathcal{C}\frac{\log(2+\tfrac{\vert \cdot\vert}{r})}{(\vert \cdot\vert+r)^d}\bigg)^2,
\end{align*}
so that \eqref{ThirdTermContribution1Bis} turns into
\begin{equation}\label{ThirdTermContribution1}
\int_{\R^d}\dd x\bigg(\int_{\bb_{\ell}(x)}\vert (\nabla v^{T}(0,\cdot)\otimes \phi^{(1)}_M e_j)\star \chi_\varepsilon\vert\bigg)^2\leq \int_{\R^d}\dd x\bigg(\int_{\bb_{\ell}(x)}\mathcal{C}\frac{\log(2+\tfrac{\vert\cdot\vert}{r})}{(\vert\cdot\vert+r)^d}\bigg)^2.
\end{equation}
We then split the \rhs of \eqref{ThirdTermContribution1} into the near-field contribution $\vert x\vert\leq 4\ell$ and the far-field contribution $\vert x\vert\geq 4\ell$. For the near-field contribution $\vert x\vert\leq 4\ell$, we use Minkowski's inequality in form of
\begin{align*}
\int_{\bb_{4\ell}} \dd x\bigg(\int_{\bb_{\ell}(x)}\mathcal{C}\frac{\log(2+\tfrac{\vert\cdot\vert}{r})}{(\vert\cdot\vert+r)^d}\bigg)^2 &\leq \bigg(\int_{\R^d}\dd x\, \mathcal{C}\frac{\log(2+\tfrac{\vert x\vert}{r})}{(\vert x\vert+r)^d}\,\vert \bb_{4\ell}\cap \bb_{\ell}(x)\vert^{\frac{1}{2}}\bigg)^2 \lesssim \ell^d \bigg(\int_{\bb_{5\ell}}  \mathcal{C}\frac{\log(2+\tfrac{\vert\cdot\vert}{r})}{(\vert\cdot\vert+r)^d}\bigg)^2.
\end{align*}
We now check that 
\begin{equation}\label{BoundLocalSupremum}
\bigg(\int_{\bb_{5\ell}}  \mathcal{C}\frac{\log(2+\tfrac{\vert\cdot\vert}{r})}{(\vert\cdot\vert+r)^d}\bigg)^2=\mathcal{C}(\ell,r)\log^4(\tfrac{\ell}{r}),
\end{equation}
using the criterion \eqref{BoundGenericAlgebraic}. Let $C,\gamma>0$ such that for any $p<\infty$, $\sup_{x\in\mathbb{R}^d}\mathbb{E}[\mathcal{C}^p(x)]^{\frac{1}{p}}\leq C p^{\frac{1}{\gamma}}$. Using Minkowski's inequality, we have
\begin{align*}
\mathbb{E}\bigg[\bigg(\int_{\bb_{5\ell}}  \mathcal{C}\frac{\log(2+\tfrac{\vert\cdot\vert}{r})}{(\vert\cdot\vert+r)^d}\bigg)^2\bigg)^p\bigg]^{\frac{1}{p}}\leq\bigg( \int_{\bb_{5\ell}}\mathbb{E}[\mathcal{C}^{2p}]^{\frac{1}{2p}}\frac{\log(2+\tfrac{\vert\cdot\vert}{r})}{(\vert\cdot\vert+r)^d}\bigg)^2 \leq C p^{\frac{1}{\gamma}}\bigg(\int_{\bb_{5\ell}}\frac{\log(2+\tfrac{\vert\cdot\vert}{r})}{(\vert\cdot\vert+r)^d}\bigg)^2\lesssim C p^{\frac{1}{\gamma}}\log^{4}(\tfrac{\ell}{r}).
\end{align*}
%
%For the part $\vert x\vert\geq 6\ell$, we note that $\bb_{5\ell}(x)\subset \mathbb{R}^d\backslash \bb_{\ell}$. Therefore, using Minkowski's inequality and Cauchy-Schwarz inequality, it yields
%
%\begin{align*}
%
%\mathbb{E}\bigg[\bigg(\sup_{\vert x\vert\geq 6\ell}\bigg(\int_{\bb_{5\ell}(x)}  \mathcal{C}\frac{\log(2+\tfrac{\vert\cdot\vert}{r})}{(\vert\cdot\vert+r)^d}\bigg)^2\bigg)^p\bigg]^{\frac{1}{p}}&\leq C p^{\frac{1}{\gamma}} \sup_{\vert x\vert\geq 6\ell}\bigg(\int_{\bb_{5\ell}(x)} \frac{\log(2+\tfrac{\vert\cdot\vert}{r})}{(\vert\cdot\vert+r)^d}\bigg)^2\\
%&\lesssim C p^{\frac{1}{\gamma}}\ell^d\int_{\mathbb{R}^d\backslash \bb_{\ell}}\frac{\log^2(1+\tfrac{\vert\cdot\vert}{r})}{(\vert\cdot\vert+r)^{2d}}\\
%&\lesssim C p^{\frac{1}{\gamma}}\log^2(\tfrac{\ell}{r}).
%
%\end{align*}
%
For the far-field contribution $\vert x\vert\geq 4\ell$, we use Cauchy-Schwarz' inequality and Fubini's theorem to obtain
\begin{equation}\label{ArgumentForScalingInte}
\begin{aligned}
\int_{\mathbb{R}^d\backslash\bb_{4\ell}} \dd x\bigg(\int_{\bb_{\ell}(x)}\mathcal{C}\frac{\log(2+\tfrac{\vert\cdot\vert}{r})}{(\vert\cdot\vert+r)^d}\bigg)^2 &\lesssim \ell^d \int_{\mathbb{R}^d\backslash \bb_{4\ell}}\dd x\int_{\bb_{\ell}(x)}\mathcal{C}^2\frac{\log^2(2+\tfrac{\vert\cdot\vert}{r})}{(\vert\cdot\vert+r)^{2d}}\\
&=\int_{\R^d}\dd x\,\mathcal{C}^2\frac{\log^2(2+\tfrac{\vert x\vert}{r})}{(\vert x\vert+r)^{2d}}\vert\mathbb{R}^d\backslash \bb_{4\ell}\cap \bb_{\ell}(x)\vert\\
&\lesssim \ell^{2d}\int_{\mathbb{R}^d\backslash \bb_{3\ell}}\mathcal{C}^2\frac{\log^2(2+\tfrac{\vert\cdot\vert}{r})}{(\vert\cdot\vert+r)^{2d}}\leq \mathcal{C}(\ell,r)\ell^d\log^2(\tfrac{\ell}{r}).
\end{aligned}
\end{equation}
\textbf{Fourth \rhs term of \eqref{eqn:funcderiv}. }We claim that 
\begin{equation}\label{FourthRHSDerivativeLarge}
\int_{\R^d}\dd x\bigg(\int_{\bb_\ell(x)}\bigg\vert\int_{0}^T\dd s\,\nabla v^{T}(s,\cdot)\otimes \nabla u^{(2)}_{M,\varepsilon}\bigg\vert\bigg)^2\leq \mathcal{C}(r,\ell)\ell^{d}(r+\log^6(\tfrac{\ell}{r})).
\end{equation}
We split the time integral into two contributions using the triangle inequality
\begin{equation}\label{SplitFourthLarge}
\begin{aligned}
&\int_{\R^d}\dd x\bigg(\int_{\bb_\ell(x)}\bigg\vert\int_{0}^T\dd s\,\nabla v^{T}(s,\cdot)\otimes \nabla u^{(2)}_{M,\varepsilon}\bigg\vert\bigg)^2\\
&\lesssim \int_{\R^d}\dd x\bigg(\int_{\bb_\ell(x)}\bigg\vert\int_{0}^{r^2}\dd s\,\nabla v^{T}(s,\cdot)\otimes \nabla u^{(2)}_{M,\varepsilon}\bigg\vert\bigg)^2+\int_{\R^d}\dd x\bigg(\int_{\bb_\ell(x)}\bigg\vert\int_{r^2}^T\dd s\,\nabla v^{T}(s,\cdot)\otimes \nabla u^{(2)}_{M,\varepsilon}\bigg\vert\bigg)^2.
\end{aligned}
\end{equation}
For the first \rhs term of \eqref{SplitFourthLarge}, we use Minkowski's inequality, \eqref{eq:9} and \eqref{eq:11} %as well as the combination of \eqref{IntegralVersionvT} with \eqref{EnergyEstiAppendix:Eq1} in form of 
%
%\begin{equation}\label{EnergyConvo}
%
%\fint_{\bb_r(x)}\int_{0}^{r^2}\dd s\,\vert\nabla v^T(s,\cdot)\vert^2\lesssim \int_{\R^d}\eta_r(\tfrac{\cdot-x}{C})g^2_r\quad\text{for any $x\in\mathbb{R}^d$,}
%
%\end{equation}
%
to obtain 
\begin{align*}
\int_{\R^d}\dd x\bigg(\int_{\bb_\ell(x)}\bigg\vert\int_{0}^{r^2}\dd s\,\nabla v^{T}(s,\cdot)\otimes \nabla u^{(2)}_{M,\varepsilon}\bigg\vert\bigg)^2\lesssim& \ell^d\bigg(\int_{\R^d}\bigg\vert\int_{0}^{r^2}\dd s\,\nabla v^{T}(s,\cdot)\otimes \nabla u^{(2)}_{M,\varepsilon}\bigg\vert\bigg)^2\\
=&\ell^d\bigg(\int_{\R^d}\fint_{\bb_r(x)}\int_{0}^{r^2}\dd s\,\vert\nabla v^T(s,\cdot)\vert\vert\nabla u^{(2)}_{M}(s,\cdot)\vert\bigg)^2\\
\leq& \ell^d \bigg(\int_{\R^d}\dd x\bigg(\fint_{\bb_r(x)}\int_{0}^{r^2}\dd s\,\vert\nabla v^T(s,\cdot)\vert^2\bigg)^{\frac{1}{2}}\\
&\times\bigg(\fint_{\bb_r(x)}\int_{0}^{r^2}\dd s\,\vert\nabla u^{(2)}_M(s,\cdot)\vert^2\bigg)^{\frac{1}{2}}\bigg)^{2}\\
\stackrel{\eqref{eq:9},\eqref{eq:11}}{\lesssim}&r\ell^d\bigg(\int_{\R^d}\mathcal{C}(R,x)\bigg(\int_{\R^d}\eta_r(\tfrac{\cdot-x}{C})\mathcal{C}^2\frac{\log^2(2+\tfrac{\vert\cdot\vert}{r})}{(\vert\cdot\vert+r)^{2d}}\bigg)^{\frac{1}{2}}\bigg)^2 \leq\mathcal{C}(R,r)r\ell^d.
\end{align*}
For the second \rhs term of \eqref{SplitFourthLarge}, we use directly \eqref{eq:10} and the second-item of \eqref{eq:11} in form of 
\begin{align*}
\int_{\R^d}\dd x\bigg(\int_{\bb_{\ell}(x)}\bigg\vert\int_{r^2}^{T}\dd s\,\nabla v^T(s,\cdot)\otimes \nabla u^{(2)}(s,\cdot)\bigg\vert\bigg)^{2}&\leq \int_{\R^d}\dd x\bigg(\int_{\bb_{\ell}(x)}\mathcal{C}(r,\cdot)\frac{\log(2+\tfrac{\vert\cdot\vert}{r})}{(\vert\cdot\vert+r)^{d}}\int_{r^2}^T\dd s\,s^{-1}\mathcal{C}(s,\cdot)\bigg)^2 \\ & \leq \mathcal{C}(\ell,r)\ell^d\log^2(\tfrac{\sqrt{T}}{r})\log^4(\tfrac{\ell}{r}),
\end{align*}
where the last line comes from the same argument as for the \rhs of \eqref{ThirdTermContribution1}.

\medskip

\textbf{Fifth and Sixth \rhs terms of \eqref{eqn:funcderiv}. }We claim that
\begin{equation}\label{FifithSixthRHSDerivativeLarge}
\begin{aligned}
\int_{\R^d}\dd x\bigg(\int_{\bb_{\ell}(x)}\vert\nabla w_2\otimes (e+\nabla \phi^{(1)}_M)\vert\bigg)^2&+\int_{\R^d}\dd x\bigg(\int_{\bb_{\ell}(x)}\vert (v^T(0,\cdot)e_j-\tfrac{1}{M}w_1)\otimes (e_i+\nabla\phi^{(1)}_M)\vert\bigg)^2\\
&\leq \mathcal{C}(\ell,r)\ell^{d+2}\log^4(\tfrac{\ell}{r}).
\end{aligned}
\end{equation}
Since the bounds of $v^T$, $\frac{1}{M}w_1$ are smaller than that of $\nabla w_2$, see \eqref{eq:11}, \eqref{eqn:bdw1} and \eqref{eqn:bdnablaw2}, we only treat the first \lhs term of \eqref{FifithSixthRHSDerivativeLarge}. Appealing to \eqref{eqn:bdnablaw2} and \eqref{eq:phiptws} yields
$$\int_{\R^d}\dd x\bigg(\int_{\bb_{\ell}(x)}\vert\nabla w_2\otimes (e+\nabla \phi^{(1)}_M)\vert\bigg)^2\leq T\int_{\R^d}\dd x\bigg(\int_{\bb_{\ell}(x)}\mathcal{C}(r,\cdot)\frac{\log^2(2+\tfrac{\vert\cdot\vert}{r})}{(\vert\cdot\vert+r)^{d}}\Big(1+\frac{\vert\cdot\vert+r}{\sqrt{T}}\Big)\bigg)^2.$$
We then argue as for \eqref{BoundLocalSupremum} and \eqref{ArgumentForScalingInte} to show that 
$$T\int_{\R^d}\dd x\bigg(\int_{\bb_{\ell}(x)}\mathcal{C}(r,\cdot)\frac{\log^2(2+\tfrac{\vert\cdot\vert}{r})}{(\vert\cdot\vert+r)^{d}}\Big(1+\frac{\vert\cdot\vert+r}{\sqrt{T}}\Big)\bigg)^2=\mathcal{C}(\ell,r)\ell^{d}(\ell^2+T)(1+\log^4(\tfrac{\ell}{r})),$$
which is sufficient for \eqref{FifithSixthRHSDerivativeLarge} since in this regime $\ell \ge \sqrt{T}\ge r$.
\medskip

{\sc Step 4. Proof of \eqref{SecondOrderDecayFlux:Eq}. }Let $C,\gamma>0$ such that the random variable $\mathcal{C}(\ell,M,r,T,\varepsilon)$ defined in \eqref{ControlFunctionalEllSmall} and \eqref{ControlFunctionalEllBig} satisfies 
$$\sup_{\ell,M,r,T,\varepsilon}\mathbb{E}[\mathcal{C}^p(\ell,M,r,T,\varepsilon)]^{\frac{1}{p}}\leq C p^{\frac{1}{\gamma}}\quad\text{for any $p<\infty$}.$$
According to \eqref{FluctuationsSigmaM}, it remains to show that there exist $C,\tilde{\gamma}>0$ such that
\begin{equation}\label{ConclusionFluctuations}
\mathbb{E}\big[\vert (q^{(2)}_{ij,M}(T,\cdot)+\sigma^{(1)}_M e_j-\mathbb{E}[q^{(2)}_{ij,M}(T,\cdot)+\sigma^{(1)}_M e_j])_r(0)\vert^p\big]^{\frac{1}{p}}\leq Cp^{\frac{1}{\tilde{\gamma}}}\log^2(\tfrac{\sqrt{T}}{r})r^{-\frac{d}{2}}\mu_\beta(T).
\end{equation}
We first apply the spectral gap inequality \eqref{msPIp} together with Minkowski's inequality and the bounds \eqref{ControlFunctionalEllSmall} and \eqref{ControlFunctionalEllBig} in form of 
\begin{align*}
\mathbb{E}\big[ &\vert (q^{(2)}_{ij,M,\varepsilon}(T,\cdot)+\sigma^{(1)}_M e_j-\mathbb{E}[q^{(2)}_{ij,M,\varepsilon}(T,\cdot)+\sigma^{(1)}_M e_j])_r(0)\vert^p\big]^{\frac{1}{p}}\\
&\stackrel{\eqref{msPIp}}{\lesssim}  \sqrt{p}\,\mathbb{E}\bigg[\Big(\int_1^\infty \dd\ell\,\ell^{-d-\beta-1}\int_{\R^d}\dd x\, \vert\partial^{\mathrm{fct}}_{x,\ell}(q^{(2)}_{ij,M,\varepsilon}(T,\cdot)+\sigma^{(1)}_M e_j)_r(0)\vert^2\Big)^{\frac{p}{2}}\bigg]^{\frac{1}{p}}\\
&\stackrel{\eqref{ControlFunctionalEllSmall},\eqref{ControlFunctionalEllBig}}{\leq} Cp^{\frac{1}{\gamma}+\frac{1}{2}}\bigg(Tr^{-d}\int_{1}^{\sqrt{T}}\dd\ell\,\ell^{d-1-\beta} (\log^2(\tfrac{\sqrt{T}}{r})+\log^2(\tfrac{\sqrt{T}}{\ell})) +\int_{\sqrt{T}}^{\infty}\dd\ell\,\ell^{-1-\beta}(\ell^2+\log^6(\tfrac{\ell}{r}))\bigg)^\frac{1}{2},
\end{align*}
which shows \eqref{ConclusionFluctuations} for the regularized quantity $(q^{(2)}_{ij,M,\varepsilon}(T,\cdot)+\sigma^{(1)}_M e_j)_r(0)$. It remains to show that
\begin{equation}\label{ConvergenceApproxEpsilon}
\mathbb{E}\big[\vert (q^{(2)}_{ij,M,\varepsilon}(T,\cdot))_r(0)-(q^{(2)}_{ij,M}(T,\cdot))_r(0)\vert^p\big]^{\frac{1}{p}}\underset{\varepsilon\downarrow 0}{\rightarrow}0\quad\text{for any $p<\infty$.}
\end{equation}
This comes from the combination of the energy estimate \eqref{EnergyEstiAppendix:Eq1}, the corrector bound \eqref{eq:phiptws} and \eqref{SmoothnessCoef} in form of
\begin{align*}
\mathbb{E}\big[\vert (q^{(2)}_{ij,M,\varepsilon}(T,\cdot))_r(0)-(q^{(2)}_{ij,M}(T,\cdot))_r(0)\vert^p\big]^{\frac{1}{p}}&=\mathbb{E}\bigg[\Big\vert\int_{0}^T\int_{\R^d}g_r\, a(\nabla u^{(2)}_{ij,M,\varepsilon}(s,\cdot)-\nabla u^{(2)}_{ij,M}(s,\cdot)\Big\vert^p\bigg]^{\frac{1}{p}}\\
&\stackrel{\eqref{EnergyEstiAppendix:Eq1}}{\lesssim_{T,r}}\mathbb{E}\bigg[\Big(\int_{\R^d}g_{\sqrt{T}}\vert\nabla\cdot (a_\varepsilon-a)\phi^{(1)}_M\vert^2\Big)^{\frac{p}{2}}\bigg]^{\frac{1}{p}} \stackrel{\eqref{eq:phiptws},\eqref{SmoothnessCoef}}{\lesssim_{T,r,p}}\varepsilon^{\alpha}. 
\end{align*}

\subsection{Proof of Corollary \ref{SecondOrderDecay}: Optimal time decay for the $2^{\text{nd}}$-order semigroup}
For notational convenience, we drop the dependences on the indices $ij$. We also recall that throughout the proof, $\mathcal{C}$ denotes a generic random constant which satisfies \eqref{MomentBoundGeneric} and can change from line to line.

\medskip

We first prove the estimate for $(u^{(2)}_M,\nabla u^{(2)}_M)$ in the average form 
\begin{equation}\label{AverageForm}
\bigg(\int_{\R^d}\eta_{\sqrt{2T}}|\nabla u^{(2)}_M(T,\cdot)|^2  \bigg)^\frac{1}{2} + \frac{1}{\sqrt{T}} \bigg(\int_{\R^d}\eta_{\sqrt{2T}}|u^{(2)}_M(T,\cdot)|^2  \bigg)^\frac{1}{2}\leq \mathcal{C}(T)T^{-1-\frac{d}{4}}\mu_\beta(T),
\end{equation}
and then deduce the analogue for the massive-less counterparts $(u^{(2)},\nabla u^{(2)})$ by passing to the limit as $M\uparrow \infty$. We finally obtain the pointwise formulation in \eqref{2NDOrder:Eq2} by applying the Schauder estimates in Proposition \ref{SchauderTheoryStandard}. The proof crucially uses the $\LL^2-\LL^1$-estimate \cite[Lemma 6]{gloria2015corrector} which reads
\begin{multline*}
	\bigg(\int_{\R^d}\eta_{\sqrt{2T}}|\nabla u^{(2)}_M(T,\cdot)|^2  \bigg)^\frac{1}{2} + \frac{1}{\sqrt{T}} \bigg(\int_{\R^d}\eta_{\sqrt{2T}}|u^{(2)}_M(T,\cdot)|^2  \bigg)^\frac{1}{2}  \\   \lesssim  \frac{1}{T}\fint_{\frac{T}{4}}^{\frac{T}{2}} \dd t \fint_0^{\sqrt{t}} \dd r \Big(\frac{r}{\sqrt{t}}\Big)^{\frac{d}{2}} \int_{\R^d}\eta_{2\sqrt{T}} \big|\big(q^{(2)}_{ij,M}(T,\cdot)-\mathbb{E}[q^{(2)}_{ij,M}(T,\cdot)] \big)_r\big|.
\end{multline*}
We then divide the $r$-integral in two parts $\int_{0}^{\sqrt{t}}=\int_{1}^{\sqrt{t}}+\int_{0}^1$. In the regime $1\leq r\leq \sqrt{t}$, we appeal to Theorem \ref{SecondOrderDecayFlux} in form of 
\begin{equation}\label{DecayRegime1}
\begin{aligned}
	 \frac{1}{T}\fint_{\frac{T}{4}}^{\frac{T}{2}} \dd t & \fint_1^{\sqrt{t}} \dd r \Big(\frac{r}{\sqrt{t}}\Big)^{\frac{d}{2}} \int_{\R^d}\eta_{2\sqrt{T}} \big|\big(q^{(2)}_{ij,M}(T,\cdot)-\mathbb{E}[q^{(2)}_{ij,M}(T,\cdot)] \big)_r\big| \\ & \stackrel{\eqref{SecondOrderDecayFlux:Eq}}{\leq}  T^{-1-\frac{d}{4}}\mu_{\beta}(T) \fint_{\frac{T}{4}}^{\frac{T}{2}} \dd t \fint_1^{\sqrt{t}} \dd r\,\C(r,T,x) \log^2 (\tfrac{\sqrt{T}}{r})\\ 
&\leq T^{-1-\frac{d}{4}}\mu_{\beta}(T) \int_{\frac{1}{2}}^{\sqrt{T}} \dd r\, \C(\tfrac{\sqrt{T}}{r},T,x) \frac{\log^2 (u)}{u^2}  \leq \C(T,x) T^{-1-\frac{d}{4}}\mu_{\beta}(T).
\end{aligned}
\end{equation}
In the regime $r\leq 1$, we combine the definition \eqref{2NDOrder:Eq8} of $q^{(2)}_{M}$ with \eqref{eq:phiptws}, \eqref{eq:9} (with $a_\varepsilon$ replaced by $a$) and \eqref{eq:10} to obtain for any $x\in\mathbb{R}^d$
\begin{align*}
\vert (q^{(2)}_M(T,\cdot))_r(x)\vert&\stackrel{\eqref{2NDOrder:Eq8}}{\leq} \mathcal{C}(x)+\int_{\R^d}g_{r}(x-\cdot)\int_{0}^{r^2}\dd s\,\vert \nabla u^{(2)}_{M}(s,\cdot)\vert+\int_{\R^d}g_{r}(x-\cdot)\int_{r^2}^T\dd s\,\vert \nabla u^{(2)}_{M}(s,\cdot)\vert\\
&\stackrel{\eqref{eq:phiptws},\eqref{eq:9}}{\leq}\mathcal{C}(r,T,x)(1+\log(\tfrac{\sqrt{T}}{r})).
\end{align*}
Using the previous bound, we deduce 
\begin{equation}\label{DecayRegime2}
\begin{aligned}
\frac{1}{T}\fint_{\frac{T}{4}}^{\frac{T}{2}} \dd t  &\frac{1}{\sqrt{t}}\int_0^{1} \dd r \Big(\frac{r}{\sqrt{t}}\Big)^{\frac{d}{2}} \int_{\R^d}\eta_{2\sqrt{T}} \big|\big(q^{(2)}_{ij,M}(T,\cdot)-\mathbb{E}[q^{(2)}_{ij,M}(T,\cdot)] \big)_r\big|\\
&\leq T^{-\frac{3}{2}-\frac{d}{4}}\fint_{\frac{T}{4}}^{\frac{T}{2}} \dd t \int_0^{1} \dd r\, \mathcal{C}(r,T)(1+\log(\tfrac{\sqrt{T}}{r}))\\
&\leq \mathcal{C}(r,T)T^{-\frac{3}{2}-\frac{d}{4}}.
\end{aligned}
\end{equation}
The combination of \eqref{DecayRegime1} and \eqref{DecayRegime2} yields \eqref{AverageForm}.
%As for the regime $r\in [0,1]$, a careful examination of the argument in \cite[Lemma 8]{clozeau2021optimal} sees that the same estimates still hold for the second-order case with the random constant having a worse stochastic integrability. Hence we finish the proof since the contribution from $r\in [0,1]$ has the same scaling as the case of first-order semigroup. This shows 
%\begin{equation*}
		%\Big(\int_{\R^d}\eta_{\sqrt{2T}}|\nabla u^{(2)}_M(T,\cdot)|^2  \Big)^\frac{1}{2} + \frac{1}{\sqrt{T}} \Big(\int_{\R^d}\eta_{\sqrt{2T}}|u^{(2)}_M(T,\cdot)|^2  \Big)^\frac{1}{2} \le \C(T,x) T^{-1-\frac{d}{4}}\mu_{\beta}(T).
%\end{equation*}
%We finally turn this into pointwise-in-space estimates by combining \eqref{SchauderParaEsti}, \eqref{Blowgrad}, \eqref{Blownongrad}, and large-scale $C^{0,1}$-estimates \eqref{LargeScaleC01Para}, \eqref{LargeScaleC01Paranondiv}.

%
\medskip

We now pass the limit as $M\uparrow\infty$ in \eqref{AverageForm}. First, note that from \eqref{strategy:naiveo2} and \eqref{2NDOrder:Eq1}, $e_M:=u^{(2)} - u^{(2)}_M$ satisfies 
 \begin{align*}
	\left\{ \begin{aligned}
		& 	\partial_\tau e_M - \nabla \cdot a \nabla e_M= 0 \quad\text{in $(0,\infty)\times\mathbb{R}^d$}, \\ 
		& e_M(0) = \nabla \cdot \big((\phi^{(1)}-\phi^{(1)}_M)a-(\sigma^{(1)}-\sigma^{(1)}_M) \bigr) e_j.
	\end{aligned}\right.
\end{align*} 
Therefore, for any fixed $R\ge \sqrt{T}\ge 1$, we deduce from the localized energy estimate \eqref{GOlm1}
\begin{equation}\label{ConvergenceM}
\begin{aligned}
\int_{\R^d}\eta_R  \big|(T\nabla e_M(T,\cdot),\sqrt{T}e_M(T,\cdot)) \big|^2 \lesssim \int_{\R^d}\eta_R \big|(\phi^{(1)}-\phi^{(1)}_M, \sigma^{(1)}-\sigma^{(1)}_M)\big|^2. 
\end{aligned} 
\end{equation}
Applying Lemma \ref{lem:1stCorApprox}, we learn that the r.h.s. of \eqref{ConvergenceM} vanishes as $M\uparrow \infty$ almost-surely. This allows us to pass the limits as $M\uparrow\infty$ in \eqref{AverageForm} and obtain the massive-less counterparts in \eqref{2NDOrder:Eq2}.

%%%%%%%%%%%%%%%%%%%%%%%%%%%%%%%%%%%%

%
\subsection{Proof of Corollary \ref{Cor:Growth2ndCor}: Construction of $2^{\text{nd}}$-order sublinear correctors}
The proof follows the lines of the ones for \cite[Corollary $2$\&$3$]{clozeau2021optimal} and splits into three steps. In the first step, we show the existence of the random fields $(\phi^{(2)}_{ij},\sigma^{(2)}_{ij})$ and the moment bound \eqref{BoundPhi2Grad} using a time integration of the second-order semigroup defined in \eqref{2NDOrder:Eq1} together with its optimal time decay estimates established in Corollary \ref{SecondOrderDecay}. In the second and third steps, we prove the fluctuation estimate \eqref{BoundPhi2Fluctu} and the growth estimate \eqref{BoundPhi2}. For notational convenience, we drop the dependences on the indices $ij$. We also recall that throughout the proof, $\mathcal{C}$ denotes a generic random constant which satisfies \eqref{MomentBoundGeneric} and can change from line to line.

\medskip

{\sc Step 1. Construction of $(\phi^{(2)},\sigma^{(2)})$. }We first build $\phi^{(2)}$ and show that it can be built such that
\begin{equation}\label{eqn:timeintphi}
\nabla \phi^{(2)} = \int_0^\infty \dd t\, \nabla u^{(2)}(t,\cdot).
\end{equation}
To begin with, $\int_{0}^\infty \dd t\, \nabla u^{(2)}(t,\cdot)$ is well defined in $\LL^2(\Omega, \mathrm{H}^1(\mathbb{R}^d))$, as an immediate consequence of the combination of \eqref{EnergyIntInside} (applied in the regime $\varepsilon\downarrow 0$),  \eqref{2NDOrder:Eq2} and Schauder's estimates in Proposition \ref{SchauderTheoryStandard} in form of the stronger statement
\begin{equation}\label{MomentBoundWeelDefine}
\bigg\vert\int_{0}^\infty \dd t\, \nabla u^{(2)}(t,x)\bigg\vert\leq \mathcal{C}(x)+\int_{1}^{\infty}\dd t\, \mathcal{C}(t,x)t^{-1-\frac{d}{4}}\mu_{\beta}(t)\leq \mathcal{C}(x)\quad \text{for any $x\in\mathbb{R}^d$}.
\end{equation}
Thus, since it is in addition curl-free, there exists $\phi^{(2)}\in \LL^2(\Omega, \mathrm{H}^1(\mathbb{R}^d))$ satisfying \eqref{eqn:timeintphi}. Note that the integral representation \eqref{eqn:timeintphi} together with \eqref{MomentBoundWeelDefine} implies \eqref{BoundPhi2Grad} for $\phi^{(2)}$. It remains to show that $\phi^{(2)}$ uniquely solves (up to an additive constant) the first equation in \eqref{EqutionPhi2Theorem}. First, since $(u^{(2)}(t,\cdot),\nabla u^{(2)}(t,\cdot))$ are stationary for all time $t>0$, it transfers to $\nabla\phi^{(2)}$ through \eqref{eqn:timeintphi} and incidentally $\mathbb{E}[\nabla\phi^{(2)}]=0$. Second, for any $T\geq 1$ and $\zeta\in\cc^{\infty}_c(\mathbb{R}^d)$, we have by integrating \eqref{strategy:naiveo2} from $0$ to $T$
\begin{equation}\label{eqn:testfnu2}
	\int_{\R^d} u^{(2)}(T,\cdot) \zeta + \int_{\R^d} (a\phi^{(1)}-\sigma^{(1)})e_j\cdot\nabla \zeta + \int_{\R^d} \nabla \zeta \cdot a \int_0^T \dd t\,\nabla u^{(2)}(t,\cdot) =0.
\end{equation}
Using \eqref{2NDOrder:Eq2}, we can pass to the limit as $T\uparrow\infty$ in \eqref{eqn:testfnu2} which shows that $\phi^{(2)}$ solves the first equation in \eqref{EqutionPhi2Theorem} almost-surely in the distributional sense. We now justify the uniqness and to this aim we consider an other solution $\bar{\phi}$ to the first equation of \eqref{EqutionPhi2Theorem} with stationary gradient and satisfying $\mathbb{E}[\nabla\bar{\phi}]=0$ as well as $\mathbb{E}[\vert\nabla\bar{\phi}\vert^2]<\infty$. The error $e:=\bar{\phi}-\phi^{(2)}$ solves
$$-\nabla\cdot a\nabla e=0.$$
Using Caccioppoli's inequality, we have for any $R\geq 1$,
\begin{equation}\label{Cacciopo}
\fint_{\bb_R}\vert\nabla e\vert^2\lesssim \frac{1}{R^2}\fint_{\bb_{2R}}\Big\vert e-\fint_{\bb_{2R}} e\Big\vert^2.
\end{equation}
Since $\nabla\phi^{(2)}$ and $\nabla\bar{\phi}$ are both stationary, of vanishing expectation and have finite second moment, $\phi^{(2)}$ and $\bar{\phi}$ are sub-linear at infinity. Consequently, as $R\uparrow\infty$, the r.h.s of \eqref{Cacciopo} vanishes whereas the l.h.s tends to $\mathbb{E}[\vert\nabla e\vert^2]$ by ergodicity. This shows that $\nabla e=0$ almost-surely so that $\bar{\phi}=\phi^{(2)}$ up to an additive constant. 

\medskip

The construction of $\sigma^{(2)}$ can be done in the same way by defining (which comes from the third equation in \eqref{EqutionPhi2Theorem})
\begin{equation}\label{Sigma2Int}
\nabla\sigma^{(2)}=\int_{0}^{\infty}\dd t\, \nabla (\nabla\times (a\nabla\phi^{(2)}+(a\phi^{(1)}-\sigma^{(1)})e_j))_{\sqrt{t}},
\end{equation}
where we recall that the subscript $\sqrt{t}$ stands for the convolution with the Gaussian kernel $g_{\sqrt{t}}$ defined in \eqref{GaussianKernel}. 
\medskip

{\sc Step 2. Proof of the fluctuations of the second-order correctors gradients. }We prove \eqref{BoundPhi2Fluctu}. We only give the arguments for $\nabla\phi^{(2)}$ as $\nabla\sigma^{(2)}$ can be treated the same way using the representation formula \eqref{Sigma2Int} instead. We split in two parts by triangle inequality
\begin{equation}\label{SplitNablaPhi}
\vert (\nabla\phi^{(2)})_r\vert\leq \vert (\nabla\phi^{(2)})_r-(\nabla\phi^{(2)}(r^2))_r\vert+\vert(\nabla\phi^{(2)}(r^2))_r\vert.
\end{equation}
Recalling the definition of the time dependent corrector $\phi^{(2)}(r^2)$ in \eqref{TimeDepSecondCor}, we deduce using the integral representation \eqref{eqn:timeintphi}, the fluctuation estimate \eqref{SecondOrderDecayFlux:Eq} and of the optimal time decay estimate \eqref{2NDOrder:Eq2} of $u^{(2)}$ that for any $x\in\mathbb{R}^d$
\begin{align*}
\vert ((\nabla\phi^{(2)})_r-(\nabla\phi^{(2)}(r^2))_r)(x)\vert&\stackrel{\eqref{TimeDepSecondCor},\eqref{eqn:timeintphi}}{=}\bigg\vert \int_{r^2}^\infty\dd t\int_{\mathbb{R}^d}g_r(x-\cdot)\nabla u^{(2)}(t,\cdot)\bigg\vert\\
&\stackrel{\eqref{2NDOrder:Eq2}}{\leq}\int_{r^2}^{\infty}\dd t\, \mathcal{C}(t,\cdot)\star g_r(x)\,t^{-1-\frac{d}{4}}\mu_\beta(t)\\
&\leq \mathcal{C}(r,x)r^{-\frac{d}{2}}\mu_\beta(r^2),
\end{align*}
and 
$$\vert(\nabla\phi^{(2)}(r^2))_r\vert\leq\mathcal{C}(r,x)r^{-\frac{d}{2}}\mu_{\beta}(r^2),$$
which, together with \eqref{SplitNablaPhi}, shows \eqref{BoundPhi2Fluctu}.

\medskip
%

%\smallskip
%We prove that 
%\[\langle | q_r^{(2)} - \langle q_r^{(2)} \rangle |^p \rangle ^\frac{1}{p} +\langle | (\nabla \phi^{(2)})_r  |^p \rangle ^\frac{1}{p}\lesssim  r^{-\frac{d}{2}}\mu_{\beta}(r^2).\]

%We start with $q_r^{(2)}$. Using triangle inequality, \eqref{SecondOrderDecayFlux:Eq} with $T=r^2$, and semigroup decay \eqref{2NDOrder:Eq2}, we have (recalling $\mu_\beta$ defined in \eqref{eq:defmubeta})
%\begin{align*}
	%\langle | q_r^{(2)} - \langle q_r^{(2)} \rangle |^p \rangle ^\frac{1}{p} & \le \langle | q_r^{(2)} -(q^{(2)}(r^2))_r|^p \rangle ^\frac{1}{p}+\langle | (q^{(2)}(r^2))_r - \langle  (q^{(2)}(r^2))_r  \rangle |^p \rangle ^\frac{1}{p}+ | \langle q^{(2)}_r -  (q^{(2)}(r^2))_r \rangle | \\ & \lesssim  \Big\langle \Big|\int_{r^2}^\infty \dd t \int_{\R^d}  \dd x g_r(x)a\nabla u^{(2)}(t,x) \Big|^p \Big\rangle ^\frac{1}{p}+ r^{-\frac{d}{2}}\mu_{\beta}(r^2) \\ & \le \int_{r^2}^\infty \dd t \int_{\R^d}\dd x g_r(x)\langle |\nabla u^{(2)}(t,x) |^p \rangle ^\frac{1}{p} + r^{-\frac{d}{2}}\mu_{\beta}(r^2) \sim  r^{-\frac{d}{2}}\mu_{\beta}(r^2).
%\end{align*}
%Here in the second inequality we used $q_r^{(2)} - (q_r^{(2)}(r^2))_r= \int_{r^2}^\infty \dd t \int_{\R^d}  \dd x g_r(x)a\nabla u^{(2)}(t,x) $, and the third term is smaller than the first by Jensen's inequality. As for $(\nabla \phi)_r$, we have similarly

{\sc Step 3. Proof of the growth of the second-order correctors. }We prove \eqref{BoundPhi2}. As before, we only give the arguments for $\phi^{(2)}$ as $\sigma^{(2)}$ can be treated the same way. Using Schauder's estimate \eqref{SchauderEllipticEsti} applied to the first equation of \eqref{EqutionPhi2Theorem} together with \eqref{BoundPhi2Grad}, we have 
$$\vert\phi^{(2)}(x)-\phi^{(2)}_1(x)\vert\leq \mathcal{C}(x)\quad\text{for any $x\in\mathbb{R}^d$}.$$
Thus, it is sufficient to show that 
\begin{equation}\label{GrowthPhi2WTS}
\vert \phi^{(2)}_1(x)-\phi^{(2)}_1(0)\vert\leq \mathcal{C}(x)\nu_\beta(\vert x\vert)\quad\text{for any $x\in\mathbb{R}^d$.}
\end{equation}
We define $R=|x|$ and we use the triangle inequality to split 
\begin{equation}\label{Growth2NDSplit}
|\phi^{(2)}_1(x) - \phi^{(2)}_1(0)| \le |\phi_R^{(2)}(x) - \phi_R^{(2)}(0)|+ |\phi_R^{(2)}(x) - \phi_1^{(2)}(x)|+|\phi_R^{(2)}(0) - \phi_1^{(2)}(0)|. 
\end{equation}
For the first r.h.s term of \eqref{Growth2NDSplit}, we use the fundamental theorem of calculus together with \eqref{BoundPhi2Fluctu} in form of 
 \begin{align*}
 \vert\phi^{(2)}_R(x)-\phi^{(2)}_R(0)\vert&=\vert x\vert\bigg\vert\int_{0}^1\dd\tau\, \frac{x}{\vert x\vert}\cdot (\nabla\phi^{(2)})_R(-\tau x)\bigg\vert\\
 &\stackrel{\eqref{BoundPhi2}}{\leq}R^{1-\frac{d}{2}}\mu_\beta(R^2)\int_{0}^1 \mathcal{C}(R,\tau x)\\
 &\leq \mathcal{C}(x)\nu_\beta(\vert x\vert).
 \end{align*}
The two last terms are bounded the same way. By the fundamental theorem of calculus and the semigroup property $g_\tau = g_{\frac{\tau}{\sqrt{2}}} *g_{\frac{\tau}{\sqrt{2}}} $, we have
\begin{align*}
\phi^{(2)}_R(0)-\phi^{(2)}_1(0) = \int_1^R \dd \tau  \frac{\partial }{\partial \tau} \phi_\tau^{(2)}(0)  &= \int_1^R \dd \tau \int_{\R^d} \dd y\,g_\tau(y)\nabla \phi^{(2)}(y)\cdot\frac{y}{\tau}\\
&= 2 \int_1^R \dd \tau \int_{\R^d} \dd x\,g_{\frac{\tau}{\sqrt{2}}}(y)\nabla \phi^{(2)}_{\frac{\tau}{\sqrt{2}}}(y)\cdot\frac{y}{\tau},
\end{align*}
with which we derive from \eqref{BoundPhi2Fluctu}
$$\vert \phi^{(2)}_R(0)-\phi^{(2)}_1(0)\vert\leq \int_1^R \dd\tau\, \tau^{-\frac{d}{2}}\mu_\beta(\tau^2)\int_{\mathbb{R}^d}\dd y\,\frac{\vert y\vert}{\tau}g_{\frac{\tau}{\sqrt{2}}}(y)\leq \mathcal{C}(x)\nu_\beta(\vert x\vert).$$
The combination of the two previous estimates yields \eqref{GrowthPhi2WTS}.

\subsection{Proof of Corollary \ref{Cor:Approx2ndCor}: Massive approximation of $2^{\text{nd}}$-order correctors}
As for \eqref{eqn:timeintphi}, we can express $\phi_M^{(2)}$ using the semigroup \eqref{2NDOrder:Eq1} :
\begin{equation}\label{IntFormulaMassive}
\phi^{(2)}_M=\int_{0}^\infty \dd t\, e^{-\frac{t}{M}}u^{(2)}_M(t,\cdot).
\end{equation}
We first prove \eqref{eq:phi2ptws}. The idea is to divide into the two regimes $t\le 1$ and $t>1$. In the regime $t\le 1$, an identical argument as for \eqref{eq:9} (for $R=1$ and $\varepsilon\downarrow 0$ together with Schauder's estimate \eqref{SchauderParaEsti}) gives
\begin{equation*}
	\Big| \int_0^1 \dd t \, e^{-\frac{t}{M}} u_M^{(2)}(t,x)\Big|\le \C(x).
\end{equation*}
When $t\ge 1$, we apply \eqref{2NDOrder:Eq2} :
\begin{equation*}
	 \int_1^\infty \dd t \, e^{-\frac{t}{M}} 	|u_M^{(2)}(t,x)| \le  \int_1^\infty \dd t \, e^{-\frac{t}{M}} \C(t,x)t^{-\frac{\beta \wedge 3}{4}} (1+\mathds{1}_{\beta=3}\log^\frac{1}{2}(t))\leq \C(x) M^{1-\frac{\beta \wedge 3}{4}}(1+\mathds{1}_{\beta=3} \log^\frac{1}{2}M).
\end{equation*}

We now prove \eqref{eq:2ndCorApproxBd}. For this purpose, we define the auxiliary  quantity 
\begin{equation}\label{eqn:tildephi2M}
	\nabla \bar{\phi}^{(2)}_M = \int_0^\infty \dd t\, e^{-\frac{t}{M}} \nabla u^{(2)}(t,\cdot).
\end{equation}
The reason why we define $\bar{\phi}^{(2)}_M$ is that it differs from $\phi^{(2)}$ only from the weight $e^{-\frac{t}{M}}$ (see \eqref{eqn:timeintphi}), while it differs from $\phi_M^{(2)}$ due to $u_M^{(2)}$ replaced by $u^{(2)}$, so it serves as the bridge between these two quantities. We treat separately the error between $\nabla \bar\phi^{(2)}_M$ and $\nabla\phi^{(2)}$ and the error between $\nabla \bar\phi^{(2)}_M$ and $\nabla\phi^{(2)}_M$.

\medskip

{\sc Step 1. Error $\nabla \bar\phi^{(2)}_M-\nabla\phi^{(2)}$. }Using the two integral representations \eqref{eqn:tildephi2M} and \eqref{eqn:timeintphi}, we have 
\begin{align*}
\bigg(\int_{\R^d}\eta_R \vert \nabla\bar{\phi}^{(2)}_M-\nabla\phi^{(2)}\vert^2\bigg)^{\frac{1}{2}}=&\bigg(\int_{\R^d} \eta_R \Big\vert \int_{0}^{\infty}\dd t\,(1-e^{-\frac{t}{M}})\nabla u^{(2)}(t,\cdot)\Big\vert^2\bigg)^{\frac{1}{2}}\\
\leq& \bigg(\int_{\R^d} \eta_R \Big\vert \int_{0}^{1}\dd t\,(1-e^{-\frac{t}{M}})\nabla u^{(2)}(t,\cdot)\Big\vert^2\bigg)^{\frac{1}{2}}\\
&+\bigg(\int_{\R^d} \eta_R \Big\vert \int_{1}^{\infty}\dd t\,(1-e^{-\frac{t}{M}})\nabla u^{(2)}(t,\cdot)\Big\vert^2\bigg)^{\frac{1}{2}}.
\end{align*}
In the regime $t\leq 1$, we use \eqref{eq:9} with $\varepsilon\downarrow 0$ in form of 
\begin{align*}
\int_{\R^d}\eta_R \Big\vert \int_{0}^{1}\dd t\,(1-e^{-\frac{t}{M}})\nabla u^{(2)}(t,\cdot)\Big\vert^2 &\leq \frac{1}{M^2}\int_{0}^1\dd t\,\int_{\R^d} \eta_R \vert\nabla u^{(2)}(t,\cdot)\vert^2\\
&\stackrel{\eqref{eq:9}}{\leq} \mathcal{C}(R)M^{-2}.
\end{align*}
In the regime $t\geq 1$, we use Minkowski's inequality and \eqref{2NDOrder:Eq2} in form of 
\begin{align*}
\int_{\R^d}\eta_R \Big\vert \int_{1}^{\infty}\dd t\,(1-e^{-\frac{t}{M}})\nabla u^{(2)}(t,\cdot)\Big\vert^2&\leq \int_1^{\infty}\dd t\,\mathcal{C}(t)(1-e^{-\frac{t}{M}})t^{-\frac{1}{2}-\frac{\beta\wedge 3}{4}}(1+\mathds{1}_{\beta=3}\log^\frac{1}{2}(t))\\
&\leq \mathcal{C}(R,M)M^{\frac{1}{2}-\frac{\beta\wedge 3}{4}}(1+\mathds{1}_{\beta=3}\log^{\frac{1}{2}}(M)).
\end{align*}
{\sc Step 2. Error $\nabla \bar\phi^{(2)}_M-\nabla\phi^{(2)}_M$. } We first appeal to \eqref{eqn:tildephi2M} and \eqref{eqn:timeintphi} together with the triangle inequality to get
\begin{align*}
\bigg(\int_{\R^d}\eta_R \vert \nabla\bar{\phi}^{(2)}_M-\nabla\phi^{(2)}_M\vert^2\bigg)^{\frac{1}{2}}\leq & \bigg(\int_{\R^d}\eta_R\Big\vert\int_{0}^1\dd t\, (\nabla u^{(2)}_M(t,\cdot)-\nabla u^{(2)}(t,\cdot))\Big\vert^2\bigg)^{\frac{1}{2}}\\
&+\bigg(\int_{\R^d}\eta_R\Big\vert \int_{1}^\infty \dd t\, e^{-\frac{t}{M}}(\nabla u^{(2)}_M(t,\cdot)-\nabla u^{(2)}(t,\cdot))\Big\vert^2\bigg)^{\frac{1}{2}}.
\end{align*}
Then, using that from \eqref{strategy:naiveo2} and \eqref{2NDOrder:Eq1} the difference $e_M:=u^{(2)}_M-u^{(2)}$ solves 
$$\left\{
    \begin{array}{ll}
        \partial_\tau e_M-\nabla\cdot a\nabla e_M=0 & \text{in $(0,\infty)\times \mathbb{R}^d$},\\
        e_M(0)=\nabla\cdot (a(\phi^{(1)}_M-\phi^{(1)})-(\sigma^{(1)}_M-\sigma^{(1)}))e_j,&
    \end{array}
\right.$$
we apply \eqref{GOlm1} and Lemma \ref{lem:1stCorApprox}, to obtain on the one hand
\begin{align*}
\bigg(\int_{\R^d}\eta_R\Big\vert\int_{0}^1\dd t\, (\nabla u^{(2)}_M(t,\cdot)-\nabla u^{(2)}(t,\cdot))\Big\vert^2\bigg)^{\frac{1}{2}}&\stackrel{\eqref{GOlm1}}{\lesssim} \bigg(\int_{\R^d}\eta_R\vert(\phi^{(1)}_M-\phi^{(1)},\sigma^{(1)}_M-\sigma^{(1)})\vert^2\bigg)^{\frac{1}{2}}\\
&\stackrel{\eqref{eq:1stCorApprox}}{\leq}\mathcal{C}(R)M^{\frac{1}{2}-\frac{\beta \wedge 3}{4}}(1+\mathds{1}_{\beta=3}\log^{\frac{1}{2}} M),
\end{align*}
and on the other hand, 
\begin{align*}
\bigg(\int_{\R^d}\eta_R\Big\vert \int_{1}^\infty \dd t\, e^{-\frac{t}{M}}(\nabla u^{(2)}_M(t,\cdot)-\nabla u^{(2)}(t,\cdot))\Big\vert^2\bigg)^{\frac{1}{2}}\stackrel{\eqref{GOlm1}}{\lesssim}&\bigg(\int_{\R^d}\eta_R\vert(\phi^{(1)}_M-\phi^{(1)},\sigma^{(1)}_M-\sigma^{(1)})\vert^2\bigg)^{\frac{1}{2}}\int_{1}^\infty \dd t\, t^{-1}e^{-\frac{t}{M}}\\
\stackrel{\eqref{eq:1stCorApprox}}{\leq}&\mathcal{C}(R)M^{\frac{1}{2}-\frac{\beta \wedge 3}{4}}(1+\mathds{1}_{\beta=3}\log^{\frac{1}{2}} M)\log M.
\end{align*}
\section*{Acknowledgement}
We would like to thank our affiliations, Institute of Science and Technology Austria and Max Planck Institute for Mathematics in the Sciences, for supporting the authors to visit each other, which greatly facilitates this work. We would like to thank Marc Josien and Quinn Winters for assistance in numerical implementation.
\appendix
\renewcommand\theequation{\thesection.\arabic{equation}}
\setcounter{equation}{0}
\section{Energy and regularity estimates for elliptic and parabolic systems}
%For notational convenience, we drop the dependences on the two unit vectors $e$ and $e'$, fixed once for all.
%

%Notice that the $v^T$ defined here is the same function as \cite[(3.11)]{clozeau2021optimal}.
%
We state in this section some energy and regularity estimates for elliptic and parabolic systems that are used in the proofs of the paper. The following lemma allows us to obtain pointwise estimates in time, for a proof we refer to \cite[Lemma 8.2]{armstrong2019quantitative}.
\begin{lemma}\label{PointwiseTimeEsti} 
Fix $r>0$, $(s,x)\in\mathbb{R}^{d+1}$ and $h\in\LL^{2}(\mathbb{R}^d)$. Assume that $v$ is a weak solution of
$$\partial_{\tau} v-\nabla\cdot a\nabla v=f+\nabla\cdot h \quad\text{in $(s-4r^2,s)\times \bb_r(x)$.}$$
We have 
\begin{align}\sup_{t\in (s-r^2,s)}\fint_{\bb_r(x)}\vert v(t,\cdot)\vert^2\lesssim \fint_{s-4r^2}^s \dd s'\fint_{\bb_{2r}(x)}\dd y\,\vert v(s',y)\vert^2+r^2\fint_{\bb_{2r}(x)}\vert h\vert^2 + r^4\fint_{\bb_{2r}(x)} f^2.  \label{Blownongrad}\\ \sup_{t\in (s-r^2,s)}\fint_{\bb_r(x)}\vert \nabla v(t,\cdot)\vert^2\lesssim \fint_{s-4r^2}^s \dd s'\fint_{\bb_{2r}(x)}\dd y\,\vert \nabla v(s',y)\vert^2+\fint_{\bb_{2r}(x)}\vert h\vert^2 + r^2\fint_{\bb_{2r}(x)} f^2.  \label{Blowgrad} \end{align}
\end{lemma}

Next, we recall some localized energy estimates for parabolic and elliptic systems, we refer to \cite[Lemma 1 \& Lemma 2]{gloria2015corrector} for a proof.
\begin{proposition}[Localized energy estimates]\label{lem:GOlm1}
Let $a$ be a time-independent coefficient field satisfying \eqref{UniformElliptiIntro} and $v,f,g$ and $h$ related through the parabolic equation on $(0,\infty)\times\mathbb{R}^d$
$$\left\{
    \begin{array}{ll}
        \partial_{\tau}v-\nabla\cdot a\nabla v=f+\nabla\cdot g & \text{in $(0,\infty)\times\mathbb{R}^d$},\\
        v(0)\equiv h.& 
    \end{array}
\right.$$
Defining the exponential weight $\eta_R:=R^{-d}\exp(-\frac{\vert\cdot\vert}{R})$, there exists $C>0$ such that for any $R\geq \sqrt{T}>0$ we have for any $x\in\mathbb{R}^d$
\begin{itemize}
\item[(i)] In the case $f=g=0$ and $h=\nabla\cdot q$ with $q\in \LL^{2}_{\mathrm{loc}}(\mathbb{R}^d)$
\begin{equation}\label{GOlm1}
\int_{\R^d}\eta_R(\tfrac{\cdot-x}{C})\vert(T\nabla v(T,\cdot),\sqrt{T} v(T,\cdot))\vert^2+\int_{\R^d}\eta_R(\tfrac{\cdot-x}{C})\bigg\vert\int_{0}^T\dd \tau\,(\nabla v(\tau,\cdot),\frac{1}{\sqrt{T}}v(\tau,\cdot))\bigg\vert^2\lesssim \int_{\R^d}\eta_R(\tfrac{\cdot-x}{C})\vert q\vert^2.
\end{equation}
\item[(ii)] If $f,g,h\in \LL^{2}_{\mathrm{loc}}(\mathbb{R}^d)$
\begin{equation}\label{EnergyEstiAppendix:Eq1}
\sup_{t<T}\int_{\R^d}\eta_R(\tfrac{\cdot-x}{C})\vert v(t,\cdot)\vert^2+\int_{0}^T\int_{\R^d}\eta_R(\tfrac{\cdot-x}{C})\vert(\tfrac{1}{\sqrt{T}}v,\nabla v)\vert^2\lesssim \int_{\R^d}\eta_R(\tfrac{\cdot-x}{C})\vert h\vert^2+\int_0^T\int_{\R^d}\eta_R(\tfrac{\cdot-x}{C})\vert(\sqrt{T}f,g)\vert^2.
\end{equation}
In particular, if $f,g,h\in \LL^{2}(\mathbb{R}^d)$, for any $T>0$
\begin{equation}\label{FullSpaceEnergyEstimate}
\sup_{t<T}\int_{\R^d}\vert v(t,\cdot)\vert^2+\int_{0}^T\int_{\R^d}\vert(\tfrac{1}{\sqrt{T}}v,\nabla v)\vert^2\lesssim \int_{\R^d}\vert h\vert^2+\int_0^T\int_{\R^d}\vert(\sqrt{T}f,g)\vert^2.
\end{equation}
\end{itemize}
Let $M\geq 1$ and $u,f$ and $g$ related through the elliptic equation
$$\frac{1}{M}u-\nabla\cdot a\nabla u=f+\nabla\cdot g\quad\text{in $\mathbb{R}^d$}.$$
We have for all $R\geq \sqrt{M}$ and $x\in\mathbb{R}^d$
\begin{equation}\label{LocalizedMassiveTerm}
\int_{\R^d}\eta_R(\tfrac{\cdot-x}{C})\vert(\tfrac{1}{\sqrt{M}}u,\nabla u)\vert^2\lesssim \int_{\R^d}\eta_R(\tfrac{\cdot-x}{C})\vert(\sqrt{M} f,g)\vert^2.
\end{equation}
\end{proposition}

%We recall in this section regularity estimates for linear elliptic and regularity systems. The next proposition is a direct consequence of classical Schauder's theory and the moment bound \eqref{RegCoef}.\nc{Check that we get exponential moments. Good reference for Schauder theory for parabolic systems ?}. For a reference, see \cite[Chapter 4\&8]{krylov1996lectures}.
%
Finally, we recall the classical Schauder theory for elliptic and parabolic systems. We emphasize that the dependence on the H\"older norm of the coefficient field $a$ in the estimates are at most polynomial which, combined with \eqref{SmoothnessCoef}, produces estimates with random constant with stretched exponential moments. For further details, we refer to \cite[Chapter 4\&8]{krylov1996lectures}.
\begin{proposition}[Elliptic and parabolic Schauder's estimates]\label{SchauderTheoryStandard}
%
%\nc{We should state it already in an annealed form}
Let $a$ be a time-independent coefficient satisfying Assumption \ref{Gaussian} and Assumption \ref{DecayCor} (in particular \eqref{SmoothnessCoef} holds). We have the two following regularity statements:
\begin{itemize}
\item Let $u$ and $f$ related through the elliptic system
$$-\nabla\cdot a\nabla u=f\quad\text{in $\mathbb{R}^d$}.$$
We have, for any $p>d$,
\begin{equation}\label{SchauderEllipticEsti}
\bigg\vert \Big(u(x)-\fint_{\bb_1(x)} u,\nabla u(x)\Big)\bigg\vert\leq \mathcal{C}(x)\bigg(\bigg(\fint_{\bb_1(x)}\vert \nabla u\vert^2\bigg)^{\frac{1}{2}}+\|f\|_{\LL^p(\bb_1(x))}\bigg)\text{ for any $x\in\mathbb{R}^d$,}
\end{equation}
where there exists $\gamma>0$ such that
\begin{equation}\label{StandardReg:Eq1}
\sup_{x\in \mathbb{R}^d}\mathbb{E}\big[\exp(\mathcal{C}^{\gamma}(x))\big]<\infty.
\end{equation}
\item Let $u$ and $f$ (time-independent) related through the parabolic system
$$\partial_{\tau}u-\nabla\cdot a\nabla u=\nabla\cdot f\text{ on $(-\infty,0)\times \mathbb{R}^d$}.$$
Provided $f\in \cc^{0,\alpha}(\mathbb{R}^d)$ for some $\alpha\in (0,1)$, we have for any $t\in (-\infty,0)$
\begin{equation}\label{SchauderParaEsti}
\bigg\vert\Big(u(t,x)-\fint_{\cc_1(t,x)} u,\nabla u(t,x)\Big)\bigg\vert\leq \mathcal{C}(x)\bigg(\bigg(\fint_{\cc_1(t,x)}\vert\nabla u\vert^2\bigg)^{\frac{1}{2}}+\|f\|_{\cc^{0,\alpha}(\bb_1(x))}\bigg)\text{ for any $x\in\mathbb{R}^d$},
\end{equation}
where $\cc_1(t,x):=(t,0)\times \bb_1(x)$ and $\mathcal{C}(x)$ satisfies \eqref{StandardReg:Eq1}.
\end{itemize}

\end{proposition}
\setcounter{equation}{0}
\section{Detailed description of the algorithm}\label{appendix:Alg}
\setcounter{equation}{0}
At the end of the paper we recall in Algorithm \ref{alg:truealg} the algorithm given in \cite{lu2021optimal} for computing approximate results for the random elliptic equations in the form of \eqref{EquationIntro}. The idea of all quantities involved are introduced in Subsection \ref{subsec:ABC} and the algorithm constitutes replacing all correctors and multipoles by their finite-domain approximation counterparts.

\begin{algorithm}
	\caption{Optimal algorithm for the approximate solution $u^{(L)}$ in $Q_L$ \label{alg:truealg}}
	\begin{algorithmic}[1]
		\State Set $M= L^{2(1-\varepsilon)}$ for some fixed $\varepsilon \in (0,\frac{1}{2})$. For $i=1,\cdots,d$, solve for the approximate first-order corrector $\phi_{i,M,L}^{(1)}$:
		\begin{equation}\label{eqn:phiML}
			\dfrac{1}{M}\phi_{i,M,L}^{(1)}-\nabla  \cdot a \nabla \phi_{i,M,L}^{(1)} =\nabla\cdot ae_i \,  \mbox{ in }Q_{2L}, \hspace{0.3in} \phi_{i,M,L}^{(1)}=0 \, \mbox{ on }\partial Q_{2L}. 
		\end{equation}
		\State Calculate the approximate homogenized coefficients via \begin{equation}\label{eqn:algahL}
			a_{\text{hom},L}e_i=\int_{\R^d} \omega q_{i,M,L}^{(1)},
		\end{equation} where \begin{equation}\label{eqn:defqiTL} q_{i,M,L}^{(1)}:=a(e_i+\nabla\phi_{i,M,L}^{(1)})\end{equation} and $\omega(x)=\frac{1}{L^d}\hat{\omega}(\frac{x}{L})$ with some nonnegative weight function $\hat{\omega} \in C_c^\infty(\QQ_1)$ and $\int_{\R^d}\hat{\omega}=1$.
		\State Compute
		\begin{equation}\label{eqn:alguhtildeL}
			u_{\text{hom},L} =\int_{\R^d}G_{\text{hom},L}\nabla\cdot h,
		\end{equation}where $G_{\text{hom},L}$ is the whole space Green function for the constant-coefficient operator $-\nabla\cdot a_{\text{hom},L} \nabla$.
		\State Solve for approximate first-order flux correctors $\sigma_{i,M,L}^{(1)}=\{\sigma_{ijk,M,L}^{(1)}\}_{j,k}$:  \begin{equation}\label{eqn:algsigma}
			\dfrac{1}{M}\sigma_{ijk,M,L}^{(1)}-\Delta \sigma_{ijk,M,L}^{(1)} =\partial_j q_{ik,M,L}^{(1)}-\partial_k q_{ij,M,L}^{(1)} \, \mbox{ in }Q_{\frac{7}{4}L}, \hspace{0.3in} \sigma_{ijk,M,L}^{(1)}=0 \, \mbox{ on }\partial Q_{\frac{7}{4}L}. 
		\end{equation}
		\State Solve for approximate second-order correctors $\phi_{ij,M,L}^{(2)}$:  \begin{equation}\label{eqn:2ndcorapprox} 
			\dfrac{1}{M}\phi_{ij,M,L}^{(2)}- \nabla\cdot a \nabla \phi_{ij,M,L}^{(2)}= \nabla\cdot (\phi_{i,M,L}^{(1)}a-\sigma_{i,M,L}^{(1)})e_j \, \mbox{ in }Q_{\frac{3}{2}L}, \hspace{0.3in} \phi_{ij,M,L}^{(2)}=0 \,\mbox{ on }\partial Q_{\frac{3}{2}L}. 
		\end{equation}
		\State For the indices \begin{equation}\label{eqn:calJ}(i,j)\in \mathcal{J}=\{(1,2),(1,3),(2,3),(2,2),(3,3)\},\end{equation} calculate  \begin{equation}\label{eqn:cijlt}
			\xi_{ij,L}^{(2)}=-\int_{\R^d}h\cdot  \nabla \Bigl(\sum_{k=1}^3\phi_{k,M,L}^{(1)}\partial_k v_{ij,L}+(2-\delta_{ij})(\phi_{ij,M,L}^{(2)}-\dfrac{a_{\text{hom}.L,ij}}{a_{\text{hom},L,11}}\phi_{11,M,L}^{(2)})\Bigr) ,
		\end{equation} where $v_{ij,L}$ denote the second-degree $a_{\text{hom},L}$-harmonic polynomials \begin{equation}\label{eqn:harmpolL}
			v_{ij,L}=(1-\dfrac{1}{2}\delta_{ij})(x_ix_j-\dfrac{a_{hL,ij}}{a_{hL,11}}x_1^2).
		\end{equation}
		\State Obtain the artificial boundary condition $u_{\text{bc}}^{(L)}$ as\begin{equation}\label{eqn:algapproxbdry}
			u_{\text{bc}}^{(L)}=(1+\phi_{i,M,L}^{(1)}\partial_i+\phi_{ij,M,L}^{(2)}\partial_{ij})\Big(u_{\text{hom},L} + \sum_{i=1}^3(\int_{\R^d} h \cdot\nabla \phi_{i,M,L}^{(1)})\partial_i G_{\text{hom},L} +\sum_{(i,j)\in\mathcal{J}}\xi_{ij,L}^{(2)}\partial_{ij} G_\text{hom,L}\Big).
		\end{equation} 
	When $d=2$ or $\beta\le 2$, we use instead
	\begin{equation}\label{eqn:2dproxbdry}
		u_{\text{bc}}^{(L)}=(1+\phi_{i,M,L}^{(1)}\partial_i)\Big(u_{\text{hom},L} + \sum_{i=1}^3(\int_{\R^d} h\cdot\nabla \phi_{i,M,L}^{(1)})\partial_i G_{\text{hom},L} \Big).
	\end{equation} 
		\State Solve for $u^{(L)}$: \begin{equation}\label{eqn:finalapprox}
			-\nabla \cdot a \nabla u^{(L)}=\nabla \cdot h\text{ in }Q_L,\hspace{0.3in} u^{(L)}=u_{\text{bc}}^{(L)}\text{ on }\partial Q_L.
		\end{equation} 
	\end{algorithmic}
\end{algorithm}

\bibliographystyle{plain}
\bibliography{references}

\end{document}